\documentclass[Afour,sageh,times]{sagej}

\usepackage{moreverb,url}

\usepackage[colorlinks,bookmarksopen,bookmarksnumbered,citecolor=red,urlcolor=red]{hyperref}

\usepackage{caption}    
\usepackage{subcaption} 
\usepackage{makecell} 
\usepackage{bm} 
\usepackage[flushleft]{threeparttable} 
\usepackage{amssymb} 
\usepackage{amsthm}
\usepackage{braket} 
\usepackage{mathtools}
\usepackage{comment} 
\usepackage{derivative} 
\usepackage{bbm} 
\usepackage{siunitx} 
\usepackage{adjustbox} 
\usepackage{xcolor}
\usepackage{bbding} 
\usepackage[normalem]{ulem} 
\newcommand{\tr}{\mathsf{T}}  
\newcommand{\Rb}{\mathbb{R}} 
\newcommand\norm[1]{\left\lVert#1\right\rVert}
\DeclareMathOperator*{\argmax}{arg\,max} 
\DeclareMathOperator*{\argmin}{arg\,min} 
\usepackage{soul} 
\usepackage{helvet} 
\newcommand\BibTeX{{\rmfamily B\kern-.05em \textsc{i\kern-.025em b}\kern-.08em
T\kern-.1667em\lower.7ex\hbox{E}\kern-.125emX}}

\usepackage[normalem]{ulem}  
\def\volumeyear{2016}
\begin{document}

\runninghead{Aoyama \textit{et~al.}}

\title{Second-Order Constrained Dynamic Optimization}

\author{Yuichiro Aoyama\affilnum{1,3}, Oswin So\affilnum{2,$\ast$}, Augustinos D. Saravanos\affilnum{2,$\ast$}, and Evangelos A. Theodorou\affilnum{1}}

\affiliation{\affilnum{1}Georgia Institute of Technology, USA\\
\affilnum{2}Massachusetts Institute of Technology, USA\\
\affilnum{3}Komatsu Ltd., Japan\\
\affilnum{$\ast$}
Contributed to this work during their time at Georgia Tech. 
}

\corrauth{Yuichiro Aoyama, School of 
Aerospace Engineering, Georgia Institute of Technology, Atlanta, GA 30332 USA}

\email{yaoyama3@gatech.edu}

\begin{abstract}
This paper provides an overview, analysis, and comparison of second-order dynamic optimization algorithms, i.e., constrained Differential Dynamic Programming (DDP) and Sequential Quadratic Programming (SQP). Although a variety of these algorithms have been proposed and used successfully, there exists a gap in understanding the key differences and advantages, which we aim to provide in this work. For constrained DDP, we choose methods that incorporate nonlinear programming techniques to handle state and control constraints, including Augmented Lagrangian (AL), Interior Point, Primal-Dual Augmented Lagrangian (PDAL), and Alternating Direction Method of Multipliers (ADMM). Both DDP and SQP are provided in single- and multiple-shooting formulations, where constraints that arise from dynamics are encoded implicitly and explicitly, respectively.  As a byproduct of the review, we propose a single-shooting PDAL DDP that has more favorable properties than the standard AL variant, such as the robustness to the growth of penalty parameters. We perform extensive numerical experiments on a variety of systems with increasing complexity to investigate the quality of the solutions, the levels of constraint violation, and the sensitivity of final solutions with respect to initialization, as well as targets. The results show that 
single-shooting PDAL DDP and multiple-shooting SQP  are the most robust methods. For multiple-shooting formulation, both DDP and SQP can enjoy informed initial guesses, while the latter appears to be more advantageous in complex systems. It is also worth highlighting that DDP provides favorable computational complexity and feedback gains as a byproduct of optimization as is.



\end{abstract}


\maketitle

\section{Introduction}\label{sec:introduction}



Second-order dynamic optimization methods are powerful optimization techniques used for optimal control of systems with nonlinear dynamics and non-quadratic cost functions. Dynamic systems with these characteristics can be
found in robotics \citep{Morimoto2003}, aerospace \citep{PELLEGRINI2020ALmulti} and transportation systems \citep{FUREY1993SQPgas}, economics \citep{Weber2011OptctrlEco}, biology \citep{Driess2018SQPmuscle} and computational neuroscience \citep{Todorov2005iLQG}, etc.
There exist two main families of methods for dynamic optimization, namely Differential Dynamic Programming (DDP) \citep{Jacobson1970ddp} and Sequential Quadratic Programming (SQP) \citep{wilson1963simplicialSQP}. Both approaches are iterative and rely on first/second-order approximations of the dynamics and the cost computed along the trajectories corresponding to each iteration. This paper provides an in-depth overview of how state and control constraints are incorporated into second-order dynamic optimization algorithms. Such constraints appear in almost all applications of trajectory optimization and iterative optimal control methods \citep{tassa2014controllimited, Howell2019ALTRO}.


\subsection{Differential Dynamic Programming}
Using Bellman's principle of optimality \citep{bellman1966dynamic}, dynamic programming (DP) divides the original optimization problem into a sequence of smaller subproblems at each time step.
Nevertheless, DP is known to suffer from ``curse of dimensionality'' because its computational and memory demands explode as the dimension of the problem increases. DDP solves this issue by considering a local approximation around the nominal trajectory. Moreover, DDP can implicitly satisfy dynamic constraints thanks to its backward and forward nature. Furthermore, DDP provides feedback gains as a byproduct of optimization. 

In practical applications, state and control constraints are of great importance. These include actuation limits and obstacles in robotics and autonomy \citep{tassa2014controllimited}, flow constraints in transportation systems \citep{Heidari1971water, Murray1979water}, and positivity constraints in biology and computational neuroscience \citep{Todorov2005iLQG}. To handle these constraints, variants of DDP have been extensively studied in the literature.
In early work, active constraints were captured in the value and state-action $Q$ function during the backward pass of DDP \citep{MurrayDDP1979,Lin91}. In the same spirit, control-limited DDP was proposed by \cite{tassa2014controllimited}, which can strictly satisfy the box control constraints by solving a Quadratic Programming (QP) in the backward pass. This method sets the feedback gains to zero when the nominal control sequence hits the control limit. As a result, these gains are not as reliable as those of normal DDP. As an extension, \cite{Zhaoming2017} presented both state- and control-constrained DDP. This method solves a similar QP with a trust region in the forward pass to surely satisfy the constraints. Consequently, the algorithm discards the feedback gains obtained in the backward pass. Moreover, due to the trust region, a good initial guess is required to achieve a task. 

Another approach to handling constraints is applying Nonlinear Programming (NLP) techniques to DDP. This includes the penalty barrier \citep{Frisch1955logarithmic, Fiacco1968sequential}, Augmented Lagrangian (AL) \citep{Powell1969AL, Hestenes1969MultiplierAG}, Interior Point (IP) \citep{byrd1999interior, wachter2005line,Wachter2006IPOPT}, Primal-Dual Augmented Lagrangian (PDAL) \citep{Gill2012PDAL, Robinson2014PDAL}, and Alternating Direction Method of Multipliers (ADMM) \citep{Boyd2011ADMM}. We note that these techniques can be used to solve dynamic optimization for robotic applications, such as \cite{Vanroye2023FATROP}.

The combination of DDP and the $\log$ barrier penalty method was proposed in \cite{Grandia2019BarrierDDP}, where the barrier function is relaxed to allow constraint violations. The penalty function's coefficient, known as a penalty parameter, is driven to zero in the original formulation but is fixed in this DDP approach. \cite{Almubarak2022} uses a similar approach with an exact barrier function. 
Despite these approximations, they are shown to work sufficiently well in practice.
AL DDP is the most widely used among these combinations of NLP with DDP for inequality constraints \citep{Plancher2017ALDDP,PELLEGRINI2020ALmulti} and equality constraints \citep{Kazadi2021eqDDP}. The method is quite robust in terms of cost reduction, but it can violate constraints especially in the early stage of optimization, where Lagrangian multipliers are inaccurate or penalty parameters are not large enough. Moreover, it can take many iterations to achieve strict feasibility.  
To alleviate this problem, the researchers proposed extensions in which algorithms switch from AL DDP to other methods when a trajectory approaches convergence \citep{Lantoine2008hybrid, Howell2019ALTRO, aoyama2021}.
IP DDP \citep{Pavlov2021IPDDP} optimizes control variables, Lagrangian multipliers, and slack variables as the original IP method using DDP. PDAL DDP \citep{Jallet2022PDALDDP} was proposed most recently. This method is similar to the AL DDP, but the Lagrangian multipliers are also optimized using DDP in contrast to the AL variant. 

Both AL and PDAL use penalty parameters that are increased during optimization. In AL, large penalty parameters are known to interrupt optimization, whereas PDAL is more robust to them as described in \cite{Robinson2014PDAL}. The PDAL DDP is presented in multiple-shooting formulation, which we elaborate on in the next paragraph. Using the Alternating Direction Method of Multipliers (ADMM), several variations of constrained DDP have been presented such as in \cite{2017SindwaniADMMDDP, zhou2020accelerated}, which split the problem into smaller subproblems that are solved sequentially. Distributed ADMM-based DDP algorithms have also been proposed for handling constrained multi-agent control problems by utilizing the parallelizable nature of ADMM \citep{Saravanos2023ADMMDDP, huang2023decentralized}.

In DDP, there exist single- and multiple-shooting formulations, which are named after the work by \cite{BOCK1984multishooting}. Single-shooting DDP is a normal DDP whose decision variables are the control sequence of the system. In this formulation, the constraints of dynamics are implicitly satisfied. On the other hand, in the multiple-shooting variant, the dynamics are handled as equality constraints or residuals that can be violated. This property allows users to initialize the algorithm with good initial guesses in both the state and control trajectories. There exist two types of multiple-shooting DDPs. In the first type, the constraints from dynamics are encoded as part of the objective and captured by the state-action $Q$ function of DDP (\cite{Jallet2022PDALDDP, PELLEGRINI2020ALmulti, Jallet2022implicit}). This type has both state and control of the system as decision variables. In the second type, the constraints are not absorbed in the cost but handled separately as residuals and the decision variables are control sequence as the single-shooting (normal) DDP. \cite{Giftthaler2018familymultiDDP} first introduced this formulation with residuals using a linear update law of dynamics in both the state and the control. \cite{Mastalli2020FDDP} further analyzed and improved the method with a nonlinear update law to solve complex tasks with high-dimensional dynamics, such as the dynamic maneuver of robots with contacts. \cite{Mastalli2022FDDPctrllim} introduced a control limit to the method. Although these methods can successfully handle complex dynamics, state constraints, such as obstacles, are not presented in contrast to the first type.

\subsection{Sequential Quadratic Programming}\label{sec:intro_SQP}
SQP was introduced as a constrained optimization technique in \cite{wilson1963simplicialSQP}, showing its power in many fields including robotics \citep {Yunt2006SUMT, yunt2007combined, yunt2011augmented, kuindersma2016optimization,posa2014direct}, aerospace \citep{kenway2014multipoint,kamyar2014aircraft} and chemical engineering \citep{lucia2013multi}.
The approach relies on sequentially solving QP subproblems with quadratically approximated objectives and linearized constraints, generating new nominal trajectories. 
For dynamical systems, SQP encodes the constraints that arise from linearized dynamics at every time step as equality constraints. 
SQP also has single-shooting and multiple-shooting formulations. The multiple-shooting variant has state and controls as decision variables. The single-shooting variant has some variations. It can be achieved by solving the same subproblem as the multiple-shooting variant and propagating the state using the system's dynamics. A technique called condensing \citep{DIEHL2002577realtimeMPC} can also be used to eliminate the state variables from the subproblems. There exists a significant amount of SQP variants that have been proposed in the literature recently for nonlinear optimal control and model predictive control (MPC) \citep{kouzoupis2018recent}.
One of the main difficulties of the SQP approach was its computational demands to solve QP subproblems. To alleviate this problem \cite{pantoja1991sequential} proposed a stage-wise version of SQP under control inequality constraints.  In the same spirit of DDP, this approach solves smaller subproblems at every time step, rather than solving a problem for a whole time horizon. Another approach for improving the time complexity is to exploit the special sparse structure of the constraints arising from dynamics while solving QP. The growth of computational complexity with respect to time horizon can be reduced from cubic to linear, which is as good as DDP. This approach was proposed with active-set \citep{GLAD1984SQPactive},
IP \citep{Steinbach1994efficientIPSQP}, and barrier \citep{
Wang2010linearMPCefficientT} methods. The same reduction is achieved by Riccati-based recursion for Linear Quadratic Regulator (LQR) problems
\citep{dohrmann1997efficient, Rao1998IPMPC, JORGENSEN2004efficientMPC}. Note that this technique is available only for multiple-shooting SQP because the sparse structure is lost by condensing \citep{Diehl2009efficientMPCandEstim}. 
It is worth noting that the recursion provides feedback gains as byproducts, as in DDP, although the gains have not been actively used. In the single-shooting variant, rather than working on the complexity, \cite{Singh2022} incorporates the DDP-style closed loop rollout, achieving a faster convergence speed. Recently, \cite{Jordana2023stagewiseSQP} proposed an approach similar to the multiple-shooting SQP, showing its superiority over single- and multiple-shooting DDPs \citep{Mastalli2020FDDP} in tasks without state constraints. This work uses LQR and ADMM to solve QP subproblems efficiently and to enforce constraints. Inspired by \cite{Stellato2020OSQP}, the matrices used in QP are not updated every iteration to improve computational efficiency.


\subsection{Other Methods}
The most straightforward approach for dynamic optimization is collocation methods, which discretize the problem in time and treat (possibly nonlinear) dynamics as equality constraints \citep{Kelly2017collocation}. The discretization process and the time step, where equality constraints are applied, are known as transcription and collocation points, respectively. Any solver can be used to solve the transcribed problem. Since the method does not have any requirements for the representation of dynamics, such as an integration scheme, it is less restrictive than those we present in this paper. However, a naive implementation does not scale well with the size of the problem. This paper focuses on methods that can take advantage of the problem structure in dynamic optimization.

\subsection{Contribution}
Although DDP and SQP have a wealth of literature separately, there exist only few works comparing the two approaches. In the unconstrained case - where SQP coincides with Newton's method with line search - the early works by \cite{Liao1992advantagesDDP,Murray1984DDPandNewton} had made a comparison of the two approaches. DDP was noted to have the same quadratic convergence properties as Newton's method, but with the advantage of solving linear equations whose size remains constant w.r.t. the time horizon. However, given the variety of constrained DDP methods available, a \textit{thorough} comparison of how these methods compare with SQP-based methods is still missing. A comparison between AL DDP and an SQP-based solver SNOPT \citep{GillSNOPT2002} has been made in \cite{Howell2019ALTRO}, where SQP was shown to converge slower in wall clock time. \cite{2017SindwaniADMMDDP}, compared ADMM DDP with control-limited DDP \cite{tassa2014controllimited} and SQP under only control constraints. \cite{Zhaoming2017} also compared their constrained DDP and SQP (SNOPT) using state and control constrained robotic tasks under a time budget, showing the superiority of DDP. However, the details of SQP, including single- or multiple-shooting etc., are not presented.  

To better understand the modern landscape
of algorithms and how they relate to each other and the optimization literature,
we compare these algorithms from derivations to performance in this paper. In addition, we propose single-shooting PDAL DDP that inherits the merit of PDAL over AL, that is, the robustness to the large penalty parameter, and add it to the comparison. 

Our thorough numerical experiments with simple to complex dynamical systems reveal the difference of the algorithms in terms of the quality of solutions, the levels of constraint violation, iterations for convergence, and the sensitivity of the final solutions with respect to initialization. The results indicate that DDP frequently shows its capability to find better local minima, whereas SQP generally performs better in satisfying constraints. It is also shown that although both DDP and SQP have multiple-shooting formulations that can enjoy informed initial guesses, the SQP variant tends to work more reliably in complex systems. Furthermore, our analysis shows that single-shooting DDPs can offer lower computational complexity, especially when the system is underactuated
and the problem has a long time horizon. This computational efficiency is one of the main motivations for users to choose DDP.

We summarize the different methods we discuss in Table \ref{tab:intro_comparison}. In the table, {\small{\fontfamily{phv}\selectfont{Single}}} and {\small{\fontfamily{phv}\selectfont{Multi.}}} represent the single and multiple shooting formulation.  $\checkmark$ and $-$ indicate whether the property is satisfied or not. For the satisfaction of the constraints, $\checkmark$ indicates that the constraints may be violated during optimization but eventually satisfied upon convergence. In contrast, $\checkmark \checkmark$ means that the trajectory is always feasible even before convergence. 

This work is organized as follows: section 
\ref{sec:non-lin-prog} includes an overview of nonlinear programming optimization techniques for static problems, including $\log$ barrier, AL, IP, PDAL, ADMM, and SQP. In section \ref{sec:DDP}, we review unconstrained DDP. Section \ref{sec:Constrained DDP} provides the derivation of constrained DDP techniques, and section \ref{sec:dynamic_SQP} gives the SQP approach for dynamic optimization. Section \ref{sec:time_complexity} analyzes time complexity of the algorithms. Section \ref{sec:experiments} presents the results of numerical experiments to compare dynamic optimization algorithms. Finally, the conclusion is drawn in section \ref{sec:conclusions}.

\begin{table*}[ht]
\footnotesize\sf\centering
\caption{Overview of the landscape of second-order constrained dynamic optimization algorithms. A problem with a dynamical system that has $n$ state, $m$ control with $N$ time horizon is considered. For algorithmic complexity, we consider inversion operations as these are more difficult to reduce compared to matrix multiplications. The complexity of unconstrained DDP is $2n^{3} + 7/2 n^{2}m+2nm^{2}+1/3m^{3}$, and its breakdown can be found in \cite{yakowitz1991convergence}.}
\label{tab:intro_comparison}
\begin{center}
\begin{threeparttable}
\begin{tabular}{c|c|cc|c|cc|c|cc}
& Barrier DDP & \multicolumn{2}{c|}{AL DDP}     & IP DDP     & \multicolumn{2}{c|}{PDAL DDP}   & ADMM DDP   & \multicolumn{2}{c}{SQP}        \\ \hline
Const. Satisfaction  &
\checkmark \checkmark  &
\multicolumn{2}{c|}{\checkmark} &
\checkmark \tnote{a} & \multicolumn{2}{c|}{\checkmark} & \checkmark & \multicolumn{2}{c}{\checkmark} \\ \hline
Single or Multi. &
Single &
\multicolumn{1}{c|}{Single} & Multi. &
Single \tnote{b}
& \multicolumn{1}{c|}{Single \tnote{c}} & Multi & Single & \multicolumn{1}{c|}{Single} & Multi. \\ \hline
Feedback Gain        & \checkmark  & \multicolumn{2}{c|}{\checkmark} & \checkmark & \multicolumn{2}{c|}{\checkmark} & \checkmark & \multicolumn{1}{c|}- & \checkmark \tnote{d}         \\ \hline
\makecell{Matrix Inversion\\ Complexity}&
$Nm^{3}$ &
\multicolumn{1}{c|}{$Nm^{3}$}     & $N(n+m)^{3}$ &   
$Nm^{3}$ & 
\multicolumn{1}{c|}{$Nm^{3}$}     & $N(n+m)^{3}$ & 
$Nm^{3,} $&
\multicolumn{1}{c|}{$N^{3}m^{3, } \tnote{e}$} & $N(n+m)^{3}, Nm^{3,} \tnote{f}$    \\ \hline
\end{tabular}
\begin{tablenotes}
  \item[a] There exist infeasible $(\checkmark)$ and feasible $(\checkmark\checkmark)$ formulations in \cite{Pavlov2021IPDDP}.
  We modify the infeasible one due to its better performance.
  \item[b] We modify the line-search filter in section \ref{subsec:IPDDP}.
  \item[c] Proposed in this paper.
  \item[d] When proper method is used as mentioned in section \ref{sec:intro_SQP} and \ref{subsec:SQP_multi}.
  \item[e] Although there exist multiple variants, here we use single-shooting with condensing. It cannot enjoy the reduction available in the multiple-shooting variant but can be reduced from $N^{3}m^{3}$ to $N^{2}m^{3}$ \citep{kouzoupis2018recent}. 
  \item[f] Originally $N^{3}(n+m)^{3}$, but the dependency on $N$ can be reduced to linear with the sparse structure as mentioned in section \ref{sec:intro_SQP} and \ref{subsec:SQP_multi}.
  \end{tablenotes}
\end{threeparttable}
\end{center}
\end{table*}

\section{Nonlinear Programming Preliminaries}\label{sec:non-lin-prog}
In this section, we review NLP methods for static optimization problems that are relevant to DDP and SQP variations for constrained dynamic optimization.
\subsection{Notation}
All methods presented here use the penalty parameter $\mu = {1}/{\rho}>0$ to penalize constraint violations.
Small values of $\mu$, corresponding to large $\rho$, will penalize constraint violations more, except for the $\log$ barrier methods, where small $\mu$ makes the barrier function closer to an indicator function.
We use $\mu$ and $\rho$ interchangeably, depending on the literature in which they were introduced.
In addition, we use both a subscript and $\nabla$ to denote derivatives depending on the situation, e.g.,
$\frac{\partial{f(x,u)}}{\partial x} = f_{x}(x,u) = \nabla_{x} f(x,u)$. For univariate functions, we omit the subscript of $\nabla$, that is, $\nabla_{x} g(x)$ = $\nabla g(x)$. In the following section, the optimization variable $x$ is a vector in $\Rb^n$. The identity matrix of size $w$ is denoted by $I_{w}$. The zero matrix of size $a$ by $b$ is denoted by $O_{a,b}$. When it is square, only one subscript is provided, i.e., $O_{a} \in \mathbb{R}^{a \times a}$. $\oslash$ and $\odot$ denote element-wise division and multiplication, respectively.  
 
\subsection{Log Barrier methods}\label{subsec:log barrier}
Consider the constrained optimization problem
\begin{align}\label{eq:nonlinear_intro_inequality}
    \min \ &f_{0}(x), \quad 
    {\text {s.t.}} \quad g(x) \leq 0, 
\end{align}
where $f_{0}: \Rb^n \to \Rb$ and $g : \Rb^n \to \Rb^w$. Barrier methods solve the problem above by minimizing a sequence of new objectives that are a sum of the original objective and a barrier function associated with the constraints. Barrier functions are parameterized by a scalar nonnegative penalty (barrier) parameter $\mu$ \cite{Frisch1955logarithmic,Fiacco1968sequential}. Here, we consider a logarithmic barrier function and a new objective denoted by $\mathcal{P}$ as
\begin{align}{\label{eq:log_barrier_penalty_function}}
    \mathcal{P}(x;\mu) = 
    \begin{cases}   
    f_{0}(x) - \mu \sum_{i=1}^{w}\log(-g_{i}(x)), & g(x) < 0,\\
    \infty, & \text{else.}
    \end{cases}
\end{align}
Observe that $\mathcal{P}$ increases as $g(x)$ gets closer to the boundary of constraints and becomes infinity when $g(x) \geq 0$. Thus, this method can only handle inequality constraints.
The function $\mathcal{P}(x;\mu)$ is minimized over $x$ iteratively with a fixed $\mu$ through Newton's method \citep{Nocedal2006numerical}. A new $x$ is obtained as 
\begin{align}\label{eq:update_x_alpha}
 x^{\rm{new}} = x + \alpha \delta x^{\ast}, 
\end{align}
where $\delta x^{\ast}$ is the solution of Newton's method and $\alpha \in (0,1]$ is a step size which reduces the penalty function. By ensuring that the cost is finite, the constraints remain satisfied for all iterations.
After each update of $x$, $\mu$ is reduced to make the minimizer of $\mathcal{P}$ closer to the true minimizer of \eqref{eq:nonlinear_intro_inequality}.

The Hessian of the penalty function $\mathcal{P}$ used in Newton's method is given by
\begin{align*}
    \nabla_{xx} \mathcal{P}(x;\mu) &= \nabla_{xx} f_{0}(x) - \mu\sum_{i}^{w}\frac{\nabla_{xx} g_{i}(x)}{g_{i}(x)} \\
    & \quad + \mu\sum_{i}^{w}\left[\frac{\nabla g_{i}(x)[\nabla g_{i}(x)]^{\tr}}{[g_{i}(x)]^2}\right],
\end{align*}
which is required to be Positive Definite (PD).
In the Hessian, the first term $\nabla_{xx} f_{0}$ is PD when the objective $f_{0}(x)$ is convex, and the third term is Positive SemiDefinite (PSD) by construction.
However, the second term may not be.
By eliminating the second term, we guarantee a PD approximation of the Hessian of $\mathcal{P}$,
which corresponds to the Gauss-Newton (GN) approximation of the Hessian \citep{Nocedal2006numerical}.  

\subsection{Augmented Lagrangian methods}\label{subsec:AL}
Consider an optimization problem similar to \eqref{eq:nonlinear_intro_inequality}, but with equality constraints explicitly included
\begin{align}\label{eq:non_linear_general_problem}
    \min \ f_{0}(x), \quad 
    {\text {s. t.}} \quad g(x) \leq 0, \quad h(x) = 0, 
\end{align}
where $h : \Rb^n \to \Rb^W$.
The Powell-Hestenes-Rockafellar (PHR) augmented Lagrangian \citep{Powell1969AL,Hestenes1969MultiplierAG,Rockfellar1974AL} is given as follows,
\begin{align}\label{eq:AL_4_general_form}
    \mathcal{L}_{\rm{A}}(x; \lambda, \nu) &= f_{0}(x) + 
    \sum_{i=1}^{w}
    \frac{\rho_{{\rm{I}},i}}{2}\norm{\left[g_{i}(x) + \frac{\lambda_{i}}{\rho_{{\rm{I}},i}}\right]_{+}}^2 \\\notag
    &\quad +\sum_{j=1}^{W}
    \frac{\rho_{{\rm{E}},j}}{2}\norm{h_{j}(x) + \frac{\nu_{j}}{\rho_{{\rm{E}},j}}}^2,
\end{align}
where $\lambda \in \mathbb{R}^{w}$ and $\nu \in \mathbb{R}^{W}$ are Lagrange multipliers associated with the inequality and equality constraints, respectively. Here, the penalty parameters $\rho_{\rm{I}}$ and $\rho_{\rm{E}}$ can be scalar or vectors of size $\mathbb{R}^{w}$ and $\mathbb{R}^{W}$, respectively - the latter is used for generality. We denote $[\cdot]_{+}$ as the projection to the nonnegative orthant, i.e.,
\begin{align}\label{eq:positive_projection}
\left[g_{i}(x) + \lambda_{i}/\rho_{{\rm{I}},i}\right]_{+} 
 = \max \{ 0, \ g_{i}(x) + \lambda_{i}/{\rho_{{\rm{I}},i}}\}.
\end{align}
The AL method has an inner and an outer loop.
In the inner loop, $\mathcal{L}_{\rm{A}}$ is minimized over $x$ with fixed penalty parameters and Lagrangian multipliers, 
until $\norm{\nabla \mathcal{L}_{\rm{A}}} \leq \epsilon_{\rm{AL}}$ is achieved, where  $\epsilon_{\rm{AL}}$ is a prescribed tolerance.
Since the inner loop problem is a minimization problem of $\mathcal{L}_{\rm{A}}$ on a single variable $x$, DDP can be easily used.  
In the outer loop, the Lagrangian multipliers are updated based on the constraint satisfaction and optimality conditions. The gradient and Hessian of $\mathcal{L}_{\rm{A}}$ are given by
\begin{align*}
     &\nabla{\mathcal{L}}_{\rm{A}} 
     = \nabla f_{0}(x) + [\nabla g(x)]_{\mathcal{A}}^{\tr}P_{\rm{I}}[g(x) + \lambda \oslash  \rho_{I}]_{+}\\
     &\hspace{10mm} 
     + [\nabla h(x)]^{\tr}P_{\rm{E}}[h(x) + \nu \oslash \rho_{\rm{E}}],\\
     &\nabla_{xx}{\mathcal{L}}_{\rm{A}}\\
     =& 
     \nabla_{xx}f_{0}(x) + [\nabla g(x)]_{\mathcal{A}}^{\tr}P_{\rm{I}}[\nabla g(x)]_\mathcal{A} \\
     &+ P_{\rm{I}}[g(x) + \lambda\oslash\rho_{I}]_{+}[\nabla_{xx} g(x)]_{\mathcal{A}}+ [\nabla h(x)]^{\tr}P_{\rm{E}}[ \nabla h(x)]\\
     & \quad + P_{\rm{E}}[\nabla h(x) + \nu \oslash \rho_{\rm{I}}]\nabla_{xx} h(x),\\
     & {\text{with}} \quad P_{\rm{I}} = {\rm{diag}}[\rho_{\rm{I}}], \quad P_{\rm{E}} = {\rm{diag}}[\rho_{\rm{E}}].
\end{align*}
%
 %
 The subscript $\mathcal{A}$ denotes the projection onto the coordinates corresponding to the active constraints, that is, elements of the indices where \eqref{eq:positive_projection} is positive. The update laws for multipliers $\lambda_i, \nu_j$ are 
\begin{align}\label{eq:AL_multiplier_update}
    \lambda^{{\rm{new}}}_{i} &= \max \{ \lambda_{i} + [\rho_{{\rm{I}},i}] g_{i}(x),0 \}, \quad  i = 1 \cdots w, \\\notag \nu^{{\rm{new}}}_{j} &= \nu_{j} + \rho_{{\rm{E}},j} h_{j}(x), \quad  j = 1 \cdots W, 
\end{align}
which are obtained by comparing the gradient of $\mathcal{L}_{\rm{A}}$ and that of the (normal) Lagrangian $\mathcal{L} = f_{0}(x) + \lambda^{\tr}g(x) + \nu^{\tr}h(x)$. The penalty parameters $\rho$ are increased if the constraint satisfaction after the inner loop is not satisfactory, determined by a tolerance for constraint satisfaction $\eta_{\rm{I}}$ as below.
\begin{align}\label{eq:AL_rho_update}
    \rho_{\rm{I},i}^{\rm{new}} = \begin{cases}
    \min \{ \rho^{\rm{max}}, \ \beta\rho_{\rm{I},i} \}, \quad g_{i}(x) \geq \eta_{\rm{I},i},\\
    \rho_{\rm{I},i}, \hspace{25mm} \text{otherwise,}
    \end{cases}
\end{align}
with $\beta>1$. The parameter $\rho_{\rm{E}}$ follows the similar law where it is updated when $\norm{h(x)} \geq \eta_{\rm{E}}$.
The tolerance for the inner loop $\epsilon_{\rm{AL}}$ starts at a moderate value and decreases as the optimization progresses to avoid local minima.
This is because convergence with small penalty parameters and suboptimal multipliers does not lead to an optimal solution of the original problem in \eqref{eq:nonlinear_intro_inequality}. Moreover, too many iterations with these parameters may get the algorithm ``trapped'' at poor local minima, e.g., a small objective with a significant constraint violation. When inner loop minimization is successful, that is, the constraint violation is sufficiently small after the inner loop, $\epsilon_{\rm{AL}}$ is reduced to allow for more inner loop iterations. Otherwise, the tolerance is kept, reset, or conservatively reduced \citep{Andrew1997AL}. The constraint satisfaction tolerance $\eta$ is updated in a similar manner. 

\subsection{Interior Point methods}\label{subsec:IP}
Consider the problem in \eqref{eq:nonlinear_intro_inequality}. 
Introducing the slack variable $s\in \mathbb{R}^{w}$, the problem is reformulated as
\begin{align}{\label{eq:non_linear_general_problem_slack}}
    \min_{x} &\ f_{0}(x), \ 
    \text{s.t.} \quad g(x) + s = 0, \ s \geq 0.  
\end{align}
 The KKT conditions for the problem in \eqref{eq:non_linear_general_problem_slack} are given by
\begin{subequations}\label{eq:IP_KKT}
\begin{align}
    \nabla f_{0}(x) + [\nabla g(x)]^{\tr}\lambda = 0,\label{eq:IP_KKT_grad}\\
    [g(x)]_{i} \lambda_{i} = 0, \label{eq:IP_KKT_cmp} \, \Leftrightarrow \, s_{i}\lambda_{i} =  0\\
    s \geq 0, \quad \lambda \geq 0, \\ 
    g(x) + s = 0.
\end{align}
\end{subequations}
Note that \eqref{eq:IP_KKT_grad} is derived from the original problem \eqref{eq:nonlinear_intro_inequality}, not from the one with slack variables \eqref{eq:non_linear_general_problem_slack} \citep{Nocedal2006numerical}. The complementary slackness \eqref{eq:IP_KKT_cmp} is relaxed by the parameter $\mu$, i.e., $s_{i}\lambda_{i} = \mu (>0)$. The parameter biases $s$ and $\lambda$ toward the feasible region, i.e., $s_{i}, \lambda_{i} \geq 0$, which is known as the central path. The relaxed condition of \eqref{eq:IP_KKT} can also be obtained by adding constraints for $s$ to the objective, forming a modified problem
\begin{align}\label{eq:IP_seen_as_relaxed_barrier}
    \min_{x, s}\ f_{0}(x) - \mu\sum_{i}^{w}\log s_{i}, \
    {\text{s.t.}} \ g(x) +s = 0,
\end{align}
similar to \eqref{eq:log_barrier_penalty_function}. The equality constraint $g(x)+s=0$ can be violated during the optimization process as long as $s \geq 0$ unlike the log barrier method \eqref{eq:log_barrier_penalty_function}. Thus, \eqref{eq:IP_seen_as_relaxed_barrier} can be seen as a relaxation of it. 
By applying Newton's method, the primal-dual system is obtained as
\begin{align}\label{eq:IP_KKT_matrix}
    &\begin{bmatrix}
    \nabla_{xx}\mathcal{L} & O_{n,w} & \nabla{g(x)}^{\tr}\\
    O_{w,n} & \bar{\Lambda} & S\\
    \nabla g(x) & I_w & O_{w}
    \end{bmatrix}
    \begin{bmatrix}
    \delta x\\
    \delta s\\
    \delta \lambda
    \end{bmatrix}
    = \\\notag
    &\hspace{30mm} -\begin{bmatrix}
    \nabla f_{0}(x) + \nabla g(x)^{\tr}\lambda \\
    \bar{\Lambda} s - {\mu} e\\
    g(x) + s
    \end{bmatrix},
\end{align}
where $S = {\rm{diag}}[s]$, $\bar{\Lambda} = {\rm{diag}}[\lambda]$, and $e = [1, \cdots, 1]^{\tr} \in \mathbb{R}^{w}$. To guarantee that the direction obtained by solving \eqref{eq:IP_KKT_matrix} is a descent direction, the matrix on the LHS, called primal-dual matrix, must have $n+w$ positive, $w$ negative, and no zero eigenvalues. When this condition is not met, 
the matrix can be modified, the details of which are found in \cite{Nocedal2006numerical,Wachter2006IPOPT}.
As in AL, optimization is performed iteratively by solving \eqref{eq:IP_KKT_matrix} and updating the decision variables via
\begin{align*}
    x^{\rm{new}} &= x + \alpha_{s}\delta x, \quad
    s^{\rm{new}} = s + \alpha_{s}\delta s,\\
    \lambda^{\rm{new}} &= \lambda + \alpha_{\lambda}\delta\lambda, \quad
\end{align*}
where the maximum step sizes $\alpha_{s}$ and $\alpha_{\lambda}$ are given by the fraction-to-the-boundary rule,
\begin{align}\label{eq:IP_fraction_to_boundary}
    \alpha_{s} &= \max \{\alpha_{s}\in(0,1]:s^{\rm{new}} \geq (1-\tau)s\}, \\\notag 
    \alpha_{\lambda} &= \max\{\alpha_{\lambda}\in (0,1]:\lambda^{\rm{new}} \geq (1-\tau){\lambda}\},
\end{align}
where $\tau (\leq 0.995)$ is a constant. $\alpha_\lambda$ is used directly to update $\lambda$. For $s$ and $x$, after determining the maximum step size, a backtracking line search is performed until both sufficient cost reduction and constraint satisfaction are achieved using a filter \citep{Fletcher1999filter} formed by a modified objective in \eqref{eq:IP_seen_as_relaxed_barrier} and constraint violation $\norm{g(x)+s}$. After the update, a new problem \eqref{eq:IP_KKT_matrix} is solved. This process is repeated until the norm on the RHS of \eqref{eq:IP_KKT_matrix} is smaller than some predetermined tolerance $\epsilon_{\rm{IP}}$, that is,
\begin{align}
    \max(\norm{\nabla \mathcal{L}}, \norm{{\Lambda}s-\mu e}, \norm{g(x)+s}) \leq \epsilon_{\rm{IP}},
\end{align}
When the condition above is satisfied, $\mu$ is reduced
for the next iteration, making the relaxed complementary slackness close to the actual one. There is an important mechanism implemented in IP-based packages such as IPOPT \cite{Wachter2006IPOPT}, which is known as feasibility restoration. When the line search cannot find an acceptable step size, even when the candidate step is smaller than a threshold, the feasibility restoration phase is invoked. In this phase, the algorithm focuses on minimizing infeasibility to find a solution that the filter can accept.

\subsection{Primal-Dual Augmented Lagrangian methods}
In this section, we first introduce PDAL \citep{Gill2012PDAL, Robinson2014PDAL} and its formulation with slack variables. Then, we outline a method for reducing the slack variables from the objective, namely PDAL.
\subsubsection{PDAL with slack variables}
Consider the problem \eqref{eq:non_linear_general_problem_slack} with constraints $h(x)=0$ as in \eqref{eq:non_linear_general_problem}. The Primal-Dual Augmented Lagrangian (PDAL) with slack variables is given by 
\begin{align}{\label{eq:PDAL}}
    &\mathcal{L}_{\rm{PD}}(x,s,\lambda,\nu;\lambda_{\rm{e}}, \nu_{\rm{e}}, \mu_{I},\mu_{E})\\\notag
    = &f_{0}(x) + \lambda^{\tr}_{e}[g(x)+s] + \frac{1}{2\mu_{\rm{I}}}
    \norm{g(x)+s}^2 \\\notag 
    &+ \frac{1}{2\mu_{\rm{I}}}
    \norm{g(x)+s + \mu_{\rm{I}}(\lambda_{\mathrm{e}} - \lambda)}^2\\\notag
    &   +\nu^{\tr}_{e}h(x) + \frac{1}{2\mu_{\rm{E}}}
    \norm{h(x)}^2 + \frac{1}{2\mu_{\rm{E}}}
    \norm{h(x)+ \mu_{\rm{E}}(\nu_{e} - \nu)}^2,
\end{align}
where $\lambda_{\mathrm{e}}$ and $\nu_{e}$ are Lagrangian multiplier estimates for $\lambda$ and $\nu$ respectively. We take $\mu_{\rm{I}} = 1/\rho_{\rm{I}}$ and $\mu_{\rm{E}}= 1/\rho_{\rm{E}}$ as scalars for simplicity, but they can also be vectors of the corresponding size.
Note that PDAL penalizes not only the violation of constraints but also the deviation of multipliers from the trajectory of the minimizers. Therefore, PDAL $\mathcal{L}_{\rm{PD}}$ is minimized for all variables, including Lagrangian multipliers, in contrast to simple AL, where the inner loop performs optimization with fixed multipliers. 
We now wish to minimize $\mathcal{L}_{\rm{PD}}$ for $x, s, \lambda, \nu$. For this purpose, we first obtain the optimal $s$ denoted by $s^{\ast}$ for other variables. To keep the notation simple and compact, we only consider inequality constraints during the derivation of $s^{\ast}$, reducing $\mu_{\rm{I}}$ and $\mu_{\rm{E}}$ to $\mu$. Since the equality constraint terms are not affected by $s$, the complete form with inequalities is easily recovered after computing $s^{\ast}$. Completing the square of $\mathcal{L}_{\rm{PD}}$ in terms of $s$ yields
\begin{align}\label{eq:PDAL_quadratic_s}
\notag
\mathcal{L}_{\rm{PD}}(x,s,\lambda;\mu)
&= f_{0}(x) + \frac{1}{\mu}\norm{g(x)+s + \mu \left(\lambda_{\mathrm{e}} 
 -\frac{\lambda}{2}\right)}^2 \\
 & \quad \ \ 
 - \mu\norm{\lambda_{\mathrm{e}} - \frac{\lambda}{2}}^{2}+ \frac{\mu}{2}\norm{\lambda_{\mathrm{e}}-\lambda}^2.
\end{align}
From this form, $s^{\ast}$ is obtained as 
\begin{align}
&s^{\ast} = -G^{\dagger}(x,\lambda) + [G^{\dagger}(x,\lambda)]_{+}.\\ 
&\text{with} \quad  
\label{eq:PDAL_Gdagger}
 G^{\dagger}(x,\lambda) = g(x) + \mu\left(\lambda_{\mathrm{e}} - {\lambda}/{2}\right).
\end{align}
Plugging this back in \eqref{eq:PDAL} (without equality constraint terms) gives the PDAL without $s$ as
\begin{align*}
    \mathcal{L}_{\rm{PD}}(x,\lambda; \mu)
    &= f_{0}(x) + \lambda_{\mathrm{e}}^{\tr}
    \left[[G^{\dagger}(x,\lambda)]_{+}-\mu\left(\lambda_{\mathrm{e}} - \frac{\lambda}{2}\right)\right]\\\notag
    & \ \ + \frac{1}{2\mu}\norm{[G^{\dagger}(x,\lambda)]_{+}- \mu\left(\lambda_{\mathrm{e}} - \frac{\lambda}{2}\right)}^2 \\
    & \ \ \  +\frac{1}{2\mu}\norm{ [G^{\dagger}(x,\lambda)]_{+} - \frac{\mu}{2}{\lambda}}^2.
\end{align*}

\subsubsection{Minimization of PDAL}\label{subsubsec:min_of_PDAL}
As in the case of AL, PDAL method also has an inner and an outer loop. In the inner one, $L_{\rm{PD}}$ is minimized over $x$, $\lambda$, and $\nu$ with fixed $\lambda_{\mathrm{e}}$, $\nu_{e}$ and $\mu$. In the outer loop, $\lambda_{\mathrm{e}}$ and $\nu_{e}$ are updated. The parameter $\mu$ ($\rho$) is monotonically decreased (increased) when the constraint violation is not satisfactory. Newton's method minimizes PDAL in the inner loop as other methods. 
In the optimization process, systems of equations 
\begin{align}\label{eq:PDAL_system_matrix}
    &\begin{bmatrix}
    \nabla_{xx}\mathcal {L}_{\rm{PD}} & - [g_x]^{\tr}_{\mathcal{A}}  &  -h_{x}^{\tr}\\
    - [g_{x}]_{\mathcal{A}} & 
    M_{\rm{I}}
    & O_{w,W}\\
    -h_{x} & O_{W,w} & \mu_{\rm{E}}I_{w}
    \end{bmatrix}
    \begin{bmatrix}
    \delta x \\ \delta \lambda \\ \delta \nu
    \end{bmatrix}\\\notag
    & \hspace{10mm}= -\begin{bmatrix}
    \nabla f_{0}(x) + h_{u}^{\tr}[2\pi_{\rm{E}}-\nu] + [g_{x}]_{\mathcal{A}}^{\tr}[2\pi_{\rm{I}}-\lambda]_{+}\\
    - \mu_{\rm{I}}\left[\pi_{\rm{I}}-\frac{\lambda}{2}\right]_{+}+\mu_{\rm{I}}\frac{\lambda}{2}\\
    -\mu_{\rm{E}} [\pi_{\rm{E}}-\nu]
    \end{bmatrix}, \\\notag
&{\text{with}}\\\notag
&\nabla_{xx}\mathcal{L}_{\rm{PD}} 
    =\nabla_{xx} f_{0}(x) + [2\pi_{\rm{I}}-\lambda]_{+}[g_{xx}]_{\mathcal{A}} + \frac{2}{\mu_{\rm{I}}}[ g_{x}]^{\tr}_{\mathcal{A}}[g_{x}]_{\mathcal{A}}\\\notag
    & \hspace{18mm}+ [2\pi_{\rm{E}}-\nu]h_{xx}+\frac{2}{\mu_{\rm{E}}} h_{x}^{\tr}[h_{x}],\\\notag
& \pi_{\rm{I}} =  \frac{g(x)}{\mu_{\rm{I}}} + \lambda_{\mathrm{e}}, \ \pi_{\rm{E}} =  \frac{h(x)}{\mu_{\rm{E}}} + \nu_{e},\
M_{\rm{I}}=\frac{\mu_{\rm{I}}}{2}\big[[I]_{\mathcal{A}} + I\big],
\end{align}
is iteratively solved. Here, the active constraints denoted by $\mathcal{A}$ is similar to the case of AL in \eqref{eq:positive_projection}, but with 
\begin{align}\label{eq:positive_projection_PDAL}
[g(x)]_{\mathcal{A}} = \max{ \{g(x) + \mu(\lambda_{\mathrm{e}} - {\lambda}/2), 0\}}.
\end{align}
The matrix on the LHS becomes numerically unstable as $\mu$ in the denominators becomes small. This instability can be alleviated by the transformation given in the original work \cite{Robinson2014PDAL}. Considering the gradient of the Lagrangian $\mathcal{L} = f_{0} + \lambda_{\mathrm{e}}^{\tr}g + \nu_{\mathrm{e}}^{\tr}h$, and that of $\mathcal{L}_{\rm{PD}}$, $[2\pi_{\rm{I}}-\lambda]_{+}$ and $2\pi_{\rm{E}}-\nu$ can be seen as a new estimate of Lagrangian multipliers, which gives the update law of multipliers in the outer loop as 
\begin{align}\label{eq:multiplier_update_PDAL}
    \lambda_{\mathrm{e}}^{\rm{new}} &= [2\pi_{\rm{I}} - \lambda]_{+} = \max \{2\pi_{\rm{I}} - \lambda, 0\}, \ 
    \nu_{e}^{\rm{new}} &=
    2\pi_{\rm{E}} - \nu.
\end{align}
Observe the similarity of $\pi_{\rm{I}}$, $\pi_{\rm{E}}$ and updated multiplies in AL in \eqref{eq:AL_multiplier_update}. Also, observe that $\lambda_{\mathrm{e}}$ is nonnegative even though this condition is not explicitly considered in \eqref{eq:PDAL_quadratic_s}.
The same update law as that in AL can be used for penalty parameters. The exit criterion of the inner loop is determined by checking the norm on the right-hand side of \eqref{eq:PDAL_system_matrix}, that is, $\norm{\nabla \mathcal{L}_{\rm{PD}}}$, $\norm{-[g(x) + \mu_{I}(\lambda_{\mathrm{e}}-\frac{\lambda}{2})]_{+}+\mu_{I}\frac{\lambda}{2}}$, $\norm{-[h(x) + \mu_{\rm{E}}(\nu_{e}-\nu)]}$.
\subsection{Alternating Direction Method of Multipliers}\label{subsec:ADMM}
Consider the following optimization problem
\begin{align}\label{eq:admm_problem}
\min f_{0}(x) + g_{0}(z), \quad \text{s. t.} \ \ 
Ax + Bz =c.
\end{align}
Here, we have two sets of variables $x$, $z$. ADMM iteratively minimizes the objective by minimizing the Augmented Lagrangian 
\begin{align*}
\mathcal{L}_{\rm{A}} = f_{0}(x) + g_{0}(z) + \lambda^{\tr}(Ax+B-c) \\\notag 
+ \frac{\rho}{2}\norm{Ax+Bz-c}^{2},
\end{align*}
where $\lambda$ is a Lagrangian multiplier. The optimization process is performed by repeating the following three updates. 
\begin{subequations}
\begin{align}
x^{l+1} &= \argmin_{x} \mathcal{L}_{\mathrm{A}}(x,z^{l},\lambda^{l}) \label{eq:admm_update_primal}
\\
z^{l+1} &= \argmin_{z} \mathcal{L}_{\mathrm{A}}(x^{l+1},z,\lambda^{l}) \\
\lambda^{l+1} &= \lambda^{l} + \rho\underbrace{(Ax^{l+1} +Bz^{l+1} -c)}_{r^{l+1}},
\end{align}
\end{subequations}
where $l$ in the superscripts of variables indicates $l$ th iteration.
The violation of the constraint given by $r^{l}$ above is also known as the primal residual.
The dual residual $d$ is derived from the first order optimality condition of \eqref{eq:admm_update_primal}
\begin{align*}
d^{l+1} = \nabla f(x^{l+1}) + A^{\tr}\lambda^{l+1} = \rho A^{\tr}B(z^{l+1}-z^{l}).
\end{align*}
As in other methods, the penalty parameter $\rho$ can be updated during optimization. However, the update is not based on the constraint violation, but on the relationship of the primal and dual residuals. When the primal feasibility is greater than the dual counterpart, $\rho$ is increased to make the relative significance of constraint violation higher in $\mathcal{L}_{\rm{A}}$. On the other hand, when the dual residual is higher, $\rho$ is decreased to prioritize the optimality of the original objective \citep{He2000ADMMvaryrho}.

\subsection{Sequential Quadratic Programming}
Consider \eqref{eq:nonlinear_intro_inequality}, again.
SQP transforms this problem to a QP with linearized constraints around current $x$, yielding 
\begin{align}{\label{eq:SQP_QP}}
    \min_{\delta x}\Big[f_{0}(x) + [\nabla f_{0}(x)]^{\tr} \delta x + \frac{1}{2}\delta x^{\tr}H\delta x\Big],\\\notag
    {\text{s.t.}} \ g(x) + \nabla g(x)\delta x \leq 0,
\end{align}
where $H$ is Hessian of Lagrangian $\mathcal{L}$ for \eqref{eq:nonlinear_intro_inequality}, i.e., $\mathcal{L} = f_{0}(x) + \lambda^{\tr}g(x)$, and $H = \nabla_{xx}\mathcal{L}$.
In practice, an approximation of $H$ is used instead of the exact one. $H$ is required to be PD as in other methods. \eqref{eq:SQP_QP} is known as a QP subproblem, whose solution is used to update $x$ with a step size $\alpha$ in \eqref{eq:update_x_alpha}. The updated $x$ leads to a new QP subproblem. SQP repeats solving the QP subproblem and updating $x$ sequentially. This process is repeated until the following KKT conditions are met:
\begin{align}\label{eq:SQP_KKT}
    \nabla f_{0}(x) + [\nabla g(x)]^{\tr}\lambda = 0,\\\notag
    g_{i}(x)\lambda_{i} \leq 0,\\\notag
    g_{i}(x) \leq 0, \quad \lambda_{i} \geq 0, \ i = 1 \cdots w.
\end{align}
To find an appropriate step size $\alpha$, AL is used as a merit function that achieves cost reduction and constraint satisfaction. Since SQP solves QP under linearized constraints, it may violate the original constraints if it attempts a large $\alpha$. The detailed derivation of the SQP based on \cite{GILL1986NPSOL, GillSNOPT2002} is provided in the appendix \ref{sup_sec:SQP}. 

\section{Unconstrained Differential Dynamic Programming}\label{sec:DDP}
This section provides a brief review of the derivation and implementation of unconstrained DDP. More details can be found in \cite{Jacobson1970ddp}. Consider the discrete-time optimal control problem
\begin{equation}
\label{eq:unconstrained_optimalcontrol}
\begin{split}
\min_{{U}}&\hspace{0.8mm}J({X},{U})=\min_{{U}}\big[{\sum_{k=1}^{N-1} l({x}_k,{u}_k)}+\Phi({x}_N)\big]\\
&\text{s.t.}\ {x}_{k+1}={f}({x}_k,{u}_k),\ k=1,...,N-1.
\end{split}
\end{equation}
where ${x}_k\in\mathbb{R}^n$, ${u}_k\in\mathbb{R}^m$ denote the state and control input of the system at the time instant $t_k$, respectively, and ${f}:\mathbb{R}^n\times \mathbb{R}^m\rightarrow\mathbb{R}^n$ corresponds to the transition dynamics function. The scalar-valued functions $l$, $\Phi$, $J$ denote the running, terminal, and total cost of the problem, respectively. We also let
\begin{align}\label{eq:XU_long_vec}
    X &= [x_{1}^{\tr}, \cdots , x_{N}^{\tr}]^{\tr} \in \mathbb{R}^{nN}, \\\notag
    U &= [u_{1}^{\tr}, \cdots , u_{N-1}^{\tr}]^{\tr} \in \mathbb{R}^{m(N-1)},
\end{align}
be the state and control trajectory over the horizon $N$. 
The cost-to-go at time $k=i$, i.e., the cost starting from $k=i$ to the end of the time horizon $N$ is given by 
\begin{align}\label{eq:DDP_cost_to_go}
    {J}_{i}(X_{i},U_{i}) = {J}_{i}(x_{i},U_{i}) = \sum_{k=i}^{N-1} l(x_k,u_k) + \Phi(x_{N}),
\end{align}
with $X_{i} =[{x}_{i}^{\tr},\dots,{x}_{N}^{\tr}]^{\tr}$, $U_{i} =[{u}_{i}^{\tr},\dots,{u}_{N-1}^{\tr}]^{\tr}$. In the above equation, the first equality holds because $X_{i}$ is recovered by $x_{i}$ and $U_{i}$. The value function is defined as the minimum cost-to-go at each state and time step via
\begin{equation}\label{value-function-definition}
V_k({x}_k):=\min_{{U}_k}J({x_{k}},{U}_{k}).
\end{equation}
Note that the value function is a function of $x_{k}$ due to the minimization with respect to $U_{k}$. 
Based on this, Bellman's principle of optimality gives the following rule:
\begin{equation}\label{eq:bellman}
    V_k(x_k) =  \min_{{u}_k} [l({x}_k,{u}_k) + V_{k+1}({x}_{k+1})].
\end{equation}
The DDP algorithm finds locally optimal solutions to \eqref{eq:unconstrained_optimalcontrol} by expanding both sides of \eqref{eq:bellman} around given nominal trajectories, $\bar{X}$, $\bar{U}$. Specifically, let us define the $Q$ function as the argument of $\min$ on the RHS of \eqref{eq:bellman},
\begin{equation}\label{eq:Qfunction}
Q_k(x_k,u_k) = l(x_k,u_k) + V_{k+1}(x_{k+1}).
\end{equation}
If we take quadratic expansions of both sides of \eqref{eq:Qfunction}, then the LHS expansion around $\bar{x}_{k}$, $\bar{u}_{k}$ gives
\begin{align}\label{eq:Q_expanded}
  & Q_k({x}_k,{u}_k)
  \approx Q_k(\bar{x}_{k}, \bar{u}_{k})+{Q}_{{x},k}^{{T}}\delta{x}_k+{Q}_{{u},k}^{{T}}\delta{u}_k \\\notag & \ +\textstyle{\frac{1}{2}}({\delta{x}}_k^{\tr}{Q}_{{xx},k}\delta{x}_k+2{\delta{x}}_k^{\tr}{Q}_{{xu},k}\delta{u}_k+{\delta{u}}_k^{\tr}{Q}_{{uu},k}\delta{u}_k).
\end{align}
where $\delta x_k:=x_k -\bar{x}_{k}$, $\delta u_k:=u_{k} - \bar{u}_{k}$ are deviations from the nominal sequences. The RHS is expanded in the same way about $\bar{x}_{k}$, $\bar{u}_{k}$, and $\bar{x}_{k+1}$, where $\delta \bar{x}_{k+1}$ is eliminated by plugging in 
quadratic approximation of dynamics:
\begin{align*}
  & \delta x_{k+1} 
  \approx f(\bar{x}_{k}, \bar{u}_{k})+{f}_{{x}}\delta{x}_k+{f_{u}}\delta{u}_k \\\notag & \ +\textstyle{\frac{1}{2}}({\delta{x}}_k^{\tr}{f}_{{xx}}\delta{x}_k+2{\delta{x}}_k^{\tr}{f}_{{xu}}\delta{u}_k+{\delta{u}}_k^{\tr}{f}_{uu}\delta{u}_k).
\end{align*}
Mapping the terms of both sides of the expanded \eqref{eq:Qfunction} gives
\begin{align*}
    & {Q}_{{x},k}={l}_{{x}}+{f_x}^{\tr}{V}_{{x},k+1}, \ \
    {Q}_{{u},k}={l_{ u}}+{f_u}^{\tr}{{V}_{{x},k+1}}, \\
    &{Q}_{{xx},k} ={l}_{{xx}}+{f_x}^{\tr}{V}_{{xx},k+1}{f_x} + V_{x,k+1} \cdot f_{xx}, \\
    &{Q}_{{uu},k} = {l_{{uu}}}+{f_u}^{\tr}{{V}_{{xx},k+1}}{f_u} + V_{x,k+1} \cdot f_{uu},\\
    &{Q}_{{xu},k} ={l_{{xu}}}+{f_x}^{\tr}{{V}_{{xx},k+1}}{f_u} + V_{x,k+1} \cdot f_{xu},
\end{align*}
where $\cdot$ in the second-order terms is tensor contraction along the first dimension.
After plugging \eqref{eq:Qfunction} and \eqref{eq:Q_expanded} into \eqref{eq:bellman}, we can explicitly optimize the value function with respect to $\delta{u}_{k}$ by taking a partial derivative and setting it to zero, obtaining the locally optimal control update
\begin{align}\label{eq:delta-u-star}
     \delta {u}^{\ast}_k&={\kappa}_{k}+{K}_{k}\delta {x}_{k}, \\\notag
    \text{with}\ \kappa_{k} &= -{Q}^{-1}_{{uu}}{Q_{{u}}},\ {K}_{k} = -{Q}^{-1}_{{uu}}{Q_{{ux}}},
\end{align}
where $\kappa_k$ and $K_k$ are known as feedforward and feedback gains, respectively.
Note that we have dropped the time indices for $Q$ to lighten the notation. To ensure convergence, $Q_{{uu}}$ must be regularized when it is not PD \citep{yakowitz1991convergence}, which is achieved with
\begin{align}\label{eq:DDP_regularization}
    Q_{uu}^{\rm{reg}} = Q_{uu} + \tau I_{m}.
\end{align}
This is equivalent to adding a cost that penalizes a large $\delta u_{k}$. Observe that $\delta {u}_{k}^{\ast}$ is computed using $V_{x,k+1}$ and $V_{xx,k+1}$. To propagate these back in time, we plug the minimizer of \eqref{eq:bellman}, i.e., $u_{k}^{\ast}$ back to the right-hand side of quadratically expanded \eqref{eq:bellman}, and
map the terms, obtaining 
\begin{align}\label{eq:value_update_with_gain}
    V_{x,k} &= Q_{{x}} + {K}^{\tr}Q_{{uu}}{\kappa}
    +{K}^{\tr}Q_{{u}} + Q_{{ux}}^{\tr}{\kappa}, \\\notag
    V_{xx,k}&= Q_{{xx}} + {K}^{\tr}Q_{{uu}}{K}
    +{K}^{\tr}Q_{{ux}} + Q_{{ux}}^{\tr}{K},
\end{align}
where regularized $Q_{uu}$ is captured in $\kappa$ and $K$. These equations are propagated backward in time with a terminal condition $V_{N}(x_{N}) = \Phi(x_{N})$, which is known as a backward pass. Then, a new state and control sequence is determined by propagating dynamics forward in time, typically with a backtracking line search. This propagation is called a forward pass. In the line search, a trial control sequence is applied to the system, generating a new state sequence 
\begin{align}\label{eq:DDP_update_control}
\bar{u}_{k}^{\rm{new}} = \bar{u}_{k} + \alpha \kappa + K\delta x_{k},\quad  \bar{x}_{k+1}^{\rm{new}} = f(\bar{x}^{\rm{new}}_{k},\bar{u}^{\rm{new}}_{k})
\end{align}
and candidate cost. Starting from step size $\alpha =1$, $\alpha$ is decreased until cost reduction is achieved. A pair of new state and control trajectories that achieves cost reduction is used new nominal pair for the next iteration. DDP repeats the backward and the forward pass until some convergence criteria are satisfied. 
In practice, when the cost is not reduced with small $\alpha$, the gains in the current iteration are discarded, and a new backward pass with larger $\tau$ is invoked.
As it can be seen in \eqref{eq:delta-u-star}, too large $\tau$ vanishes the information of $Q_{uu}^{-1}$ in the gains. Therefore, a proper choice of $\tau$ is important. An efficient scheduling technique is found in \cite{TassaDDP2012}.

Although the original DDP is introduced with second-order expansion of dynamics, in this work we consider only the first-order expansions as it is computationally cheaper and tends to be more numerically stable \citep{li2004iterative}. In fact, many practitioners neglect the second-order terms, resulting in the so-called iLQR family of algorithms: \cite{TassaDDP2012, Giftthaler2018familymultiDDP, Boutselis2021Liegroup}.

\section{Constrained Differential Dynamic Programming}\label{sec:Constrained DDP}
In this section, we present the main constrained variations of DDP which emerge through combinations with NLP techniques.
For all algorithms, we demonstrate the impact of constraints on the objective function $J$ and then examine how the $Q$ functions are modified compared to the unconstrained DDP. For notational brevity, we use the concatenated variable $y_{k} = [u_{k}^{\tr}, x_{k}^{\tr}]^{\tr}$, which allows us to write
\begin{align*}
Q_{y,k} = [Q_{u,k}, Q_{x,k}], \ 
Q_{yy,k} = \begin{bmatrix}
Q_{uu,k} & Q_{ux,k} \\ Q_{xu,k} & Q_{uu,k}
\end{bmatrix}.
\end{align*}
Adding inequality constraints $g\leq0$ to \eqref{eq:unconstrained_optimalcontrol}, we consider
\begin{align}\label{eq:constrained_OCP_dyn_sys}
&\min_{U}{J}(X,U) = \min_{U} \Big[\sum_{k=1}^{N-1} l(x_{k},u_{k})\Big] + \Phi(x_{N}), \\ \notag
&{\text{s.t.}} \ x_{k+1} = f(x_{k}, u_{k}), \ g(x_{k}, u_{k}) \leq 0, \quad k = 1, \cdots N-1,  \\\notag
                    & \quad \ g(x_{N}) \leq 0.
\end{align}
Note that in the final time step, the constraint is a function of only $x_{N}$.
To simplify our argument, we keep the dimension of $g_{N}(x_{N}) \in \mathbb{R}^{w}$, which is the same as that of $g(x_{k}, u_{k})$. 


\subsection{Log Barrier DDP}
First, we show how problem \eqref{eq:constrained_OCP_dyn_sys} can be addressed through combining DDP and $\log$ barrier method. Following \eqref{eq:log_barrier_penalty_function}, by incorporating the inequality constraints in the objective, we have a modified problem with cost $\hat{J}$ as
\begin{align}\label{eq:log_barr_DDP_problem}
\notag
\min_{U}\hat{J}(X,U) &= \min_{U} \sum_{k=1}^{N-1}\big[ l(y_{k}) - \mu \sum_{i=1}^{w}\log [ -g(y_{k})]\big] \\
& + \Phi(x_{N})- \mu \sum_{i=1}^{w}\log [ -g(x_{N})],
\end{align}
which modifies the value function as 
\begin{align*}
    {V}_{k}(x_{k}) = \min_{u_{k}}\big[l(y_{k}) - \mu\sum_{i}^{w}\log [-g(y_{k})] + {V}_{k+1}(x_{k+1})\big].
\end{align*}
On the RHS, the second term is added compared to \eqref{eq:bellman}. Let the argument for the $\min$ operator be $\hat{Q}$. The derivatives of $\hat{Q}$ are
\begin{align*}
    \hat{Q}_{y} &= Q_{y} - \mu \sum_{i=1}^{w}\frac{g_{i,y}(y)}{g_{i}(y)},\\
    \hat{Q}_{yy} &= Q_{yy} - \mu \sum_{i=1}^{w}\frac{g_{i,yy}(y)}{g_{i}(y)} + \mu \sum_{i=1}^{w}\frac{g_{i,y}(y)g_{i,y}(y)^{\tr}}{g_{i}(y)^{2}},
\end{align*}
where the second term of $\hat{Q}_{uu}$ can be omitted for the GN approximation as in section \ref{subsec:log barrier}. This approximation was rederived by augmenting the barrier term as an additional element of the state in dynamics in \cite{Almubarak2022}. In this work, only the first-order derivatives of constraints contribute to the second-order derivatives of $Q$. The approximation is also indirectly regularizing the Hessian, making the optimization problem more well-conditioned for DDP.
There exist techniques that relax the $\log$ barrier to facilitate optimization, accepting constraint violation \cite{Grandia2019BarrierDDP}. In this paper, however, we use the exact (not relaxed) barrier function because it can keep trajectories always feasible even before convergence. This is a unique algorithm property that distinguishes the method from others. As mentioned in section \ref{subsec:log barrier}, $\mu$ needs to be reduced as optimization progresses. However, as reported in previous work, a single small value of $\mu$ is sufficient in most cases \cite{Almubarak2022}.

\subsection{Augmented Lagrangian DDP}\label{sec:AL_DDP}
In this section, we first introduce two constrained DDP formulations, that is, single- and multiple-shooting DDPs. Then, we derive AL DDP in both formulations.
\subsubsection{Single and multiple shooting constrained DDP}
Single-shooting DDP is the normal DDP explained in section \ref{sec:DDP}, where the decision variables are only control variables. In this formulation, equality constraints arise from the dynamics are implicitly satisfied. On the other hand, in multiple-shooting DDP, the constraints from dynamics can be violated and treated as residual in dynamics, or constraint violation penalized in the cost. Since the algorithm can violate the dynamics, computing an initial trajectory of the state is easy in contrast to the case of single-shooting. In a reaching task of a vehicle, for example, one can draw its trajectory by linearly interpolating the initial point to the target, or one could use sampling-based algorithms, e.g., rapidly exploring random tree \citep{lavelle1998RRT}. We test the multiple-shooting DDP by feeding these initial guesses later in section \ref{sec:experiments}. 

As mentioned in section \ref{sec:introduction}, two types of multiple-shooting DDP exist. In the first type, the decision variables are a pair of current control and next state, i.e., $u_{k}$ and $x_{k+1}$, which are related by equality constraints from the dynamics. In this formulation, the infeasibility of equality constraints is penalized as part of the cost. The linear update law updates all variables, including state and control. Inequality constraints, such as obstacles, can easily be added as part of the cost. However, in the second type, the infeasibility of the equality constraints, also known as residual or defect, is not part of the cost. Instead, they are captured in $Q$ functions in the backward pass through dynamics and reduced in the forward pass, with the linear \citep {Giftthaler2018familymultiDDP} or nonlinear \citep{Mastalli2020FDDP} update law. In these methods, in the second type, decision variables are only control variables,  which differs from the methods of the first type. Nevertheless, they still have the property of multiple-shooting because they can accept dynamically infeasible trajectories. They show their power in solving problems with complex dynamics, such as humanoid robots with contacts. However, state constraints such as obstacles have not been presented. We use the first type as a representative of multiple-shooting DDP because we are more interested in environments with constraints.

\subsubsection{Inequality constrained DDP with single shooting}
Here, we derive the AL DDP with the single shooting-method, where only the inequality constraints $g$ are considered. Following \eqref{eq:AL_4_general_form}, by adding penalty terms from constraints to the objective, we have from \eqref{eq:constrained_OCP_dyn_sys} that
\begin{align}\label{eq:AL_DDP_cost}
    &\min_{U}{\hat{J}}(X,U)\\\notag
    =& \min_{U} \sum_{k=1}^{N-1}\left[ l_{k}(y_{k}) +\sum_{i=1}^{w}\frac{\rho_{k,i}}{2}\norm{\left[ g_{i}(y_{k}) + \frac{\lambda_{k,i}}{\rho_{k,i}}\right]_{+}}^{2} \right]\\\notag
    &\quad
    + \Phi(x_{N})+ \sum_{i=1}^{w}\frac{\rho_{N,i}}{2}\norm{\left[g(x_{N})+ \frac{\lambda_{N,i}}{\rho_{N,i}}\right]_{+}}^2,
\end{align}
for fixed penalty parameter $\rho_{k} \in \mathbb{R}^{w}$ and Lagrangian multiplier $\lambda_{k} \in \mathbb{R}^{w}$. In the inner loop, DDP is used to solve \eqref{eq:AL_DDP_cost} with modified $Q$ functions whose derivatives are given below.
\begin{align}
    \notag
    \hat{Q}_{y} &= Q_{y} + [g_{y}(y_{k})]_{\mathcal{A}}^{\tr}P_{k}\left[g(y_{k}) + \lambda_{k} \oslash \rho_{k} \right]_{+},\\\label{eq:ALDDP_Qs}
    \hat{Q}_{yy} &= Q_{yy} + [g_{y}(y_{k})]_{\mathcal{A}}^{\tr}P_{k}[g_{y}(y_{k})]_{\mathcal{A}} \\\notag
    & \hspace{12mm} + P_{k}\left[g(y_{k})+\lambda_{k} \oslash \rho_{k} \right]_{+}[g_{yy}(y_{k})]_{\mathcal{A}},
\end{align}
where $P_{k} = {\rm{diag}}[\rho_{k}]$. In the outer loop, the multipliers and penalty parameters are updated by \eqref{eq:AL_multiplier_update} and \eqref{eq:AL_rho_update}. Note that in our implementation, we vary $\rho$ for the constraints but keep the same throughout the time horizon, that is, $\rho_{k} =\rho$ for all $k$. This is because even if a large constraint violation exists at time step $k$ in an iteration, this may not be the case in the next iteration at the same time step. Rather, the same constraint is more likely to be violated.
The tolerance for the constraint satisfaction, $\eta \in \mathbb{R}^{w}$, is used to determine whether the satisfaction of the constraints is sufficient. Let $i_{\rm{s}}$ and $i_{\rm{f}}$ be indices of sufficient and insufficient constraint satisfaction. i.e.,
\begin{align*}
 i_{\rm{s}} = \{i| \ g_{i}(x_{k},u_{k}) \leq \eta_{i}, \ \forall k =1\cdots N\}, \quad 
 i_{\rm{f}} = i - i_{\rm{s}},
\end{align*}
where the subtraction in the second equation is for sets. Initialized by $\eta_{0}$, using these, $\eta$ is updated by
\begin{align*}
\eta_{i_{\rm{s}}}^{\rm{new}} =\max (\eta_{i_{\rm{s}}}
/{{\rho}^{{\beta}_{\eta}}_{i_{\rm{s}}}}, \ \eta_{\rm{min}} ), \quad 
\eta_{i_{\rm{f}}}^{\rm{new}} = \max({\eta_{0}/{\rho^{\alpha_{\eta}}_{i_{\rm{f}}}},\ \eta_{\rm{min}}}),
\end{align*}
where $\alpha_{\eta}, \beta_{\eta} \in (0,1)$.
The tolerance of the inner loop $\epsilon_{\rm{AL}}$ in  \ref{subsec:AL} is updated when the inner loop is successful. The inner loop is considered successful when it satisfies the following conditions,
\begin{align*}
&g_{i^{\dagger}}(x_{k^{\dagger}},u_{k^{\dagger}}) \leq \eta_{i^{\dagger}}, \\
\text{with} \ {k}^{\dagger}, {i}^{\dagger} &= \argmax_{k \in [1,\cdots N], i\in [1,\cdots w]} 
 g_{i}(x_{k},u_{k}).
\end{align*}
i.e., the largest constraint violation is below the specified tolerance. We also reduce the $\epsilon_{\rm{AL}}$ even when the inner loop is not successful to let the inner loop run more as optimization proceeds as in AL for static problems.

\subsubsection{Inequality and equality constrained DDP with multiple shooting}\label{subsubsec:DDP multi}
In this section, we derive multiple shooting AL DDP.
We augment $u_{k}$, $x_{k+1}$, having $\tilde{u}_{k}\in\mathbb{R}^{m+n}$ and define an operator $\Pi$ that extracts $x_{k+1}$ from $\tilde{u}_{k}$
\begin{align*}
    \tilde{u}_{k} = \begin{bmatrix}
    u_{k}^{\tr} & x_{k+1}^{\tr}
    \end{bmatrix}^{\tr},
    \quad \Pi(\tilde{u}_{k}) = [O_{1, m}, e]
    \tilde{u}_{k}  = x_{k+1},
\end{align*}
with $e = [1,\cdots,1] \in \mathbb{R}^{1\times n}$. Equality constraints from dynamics are given as follows.
\begin{align*}
    h(x_{k},\tilde{u}_{k}) = x_{k+1}-f(x_{k},u_{k}), \ k=1, \cdots, N-1.
\end{align*}
Replacing constraints from dynamics with these and using augmented control, the problem in \eqref{eq:constrained_OCP_dyn_sys} is modified as
\begin{align}\label{eq:OCP_dyn_sys_equality}
&\min_{\tilde{U}}{\hat{J}}(X,\tilde{U}) = \min_{\tilde{U}} \Big[\sum_{k=1}^{N-1} l(x_{k},\tilde{u}_{k})\Big] + \Phi(x_{N}), \\ \notag
&{\text{s.t.}} \  h(x_{k},\tilde{u}_{k}) = 0, \quad g_{k}(x_{k}, \tilde{u}_{k}) \leq 0, \quad k = 1, \cdots N-1,\\\notag
&\hspace{8mm} g_{N}(x_{N}) \leq 0.
\end{align}
As in the single shooting case, we handle the constraints as  a part of the cost.  Here, the running cost, value function, and constraints are functions of $\tilde{u}_{k}$, which gives derivatives as
\begin{align*}
    l_{\tilde{u}} = \begin{bmatrix}
    l_{u} \\ O_{n, 1}\\
    \end{bmatrix}, & \
    \frac{{\partial V(x')}}{\partial \tilde{u}} = 
    \Big(\frac{\partial x'}{\partial \tilde{u}}\Big)^{\tr}\Big( \frac{\partial V(x')}{\partial x'}\Big)
    = \begin{bmatrix}
    O_{m}\\V_{x}(x') 
    \end{bmatrix},
    \\
    \frac{\partial^{2} V(\Pi(\tilde{u}))}{\partial \tilde{u}^{2}} &= \frac{\partial}{\partial{\tilde{u}}}
    \begin{bmatrix}
    O_{m,1}\\V_{x}(x')
    \end{bmatrix}
    = \begin{bmatrix}
    O_{m,m} & O_{m,n} \\
    O_{n,m} & V_{xx}(x')
    \end{bmatrix}, \\ 
    h_{\tilde{u}} &= [h_{u}, h_{x'}],\ 
    g_{\tilde{u}} = [g_{u}, g_{x'}],  
\end{align*}
 where we drop time index $k$ and use $x' = x_{k+1}$ for simplicity. 
The $Q$ function for this problem is now 
\begin{align*}
    \hat{Q}(x_{k},\tilde{u}_{k}) &= l(x_{k}, u_{k}) + V(x_{k+1}) \\\notag 
    & \ + \sum_{i}^{w}\frac{[\rho_{\rm{I}}]_{k,i}}{2}\norm{\Big[g_{i}(x_{k},u_{k}) + \frac{\lambda_{k,i}}{\rho_{\rm{I}}}\Big]_{+}}^{2} \\\notag
    &\quad \quad +\sum_{j}^{n}\frac{[\rho_{\rm{E}}]_{k,j}}{2}\norm{h_{j}(x_{k},u_{k}) + \frac{\nu_{k,j}}{\rho_{\rm{I}}}}^{2}.
\end{align*}
Using the augmented control as $\tilde{y}_{k} = [\tilde{u}_{k}^{\tr}, x_{k}^{\tr}]^{\tr}$, we have 
\begin{align}\label{eq:Q_with_z_tilde}
{Q}_{\tilde{y}} = [(l_{\tilde{u}}+ V_{\tilde{u}})^{\tr}, l_{x}^{\tr}]^{\tr}, \
{Q}_{\tilde{y}\tilde{y}} =
\begin{bmatrix}
l_{\tilde{u}\tilde{u}}+ V_{\tilde{u}\tilde{u}} & l_{\tilde{u}x} \\l_{x\tilde{u}} & l_{xx}
\end{bmatrix}.
\end{align}
The derivatives of $\hat{Q}$ is given by
\begin{align*}
\hat{Q}_{\tilde{y}} &= Q_{\tilde{y}}
+[g_{\tilde{y}}]_{\mathcal{A}}^{\tr}P_{\rm{I}}[g(y_{k})+\lambda_{k}\oslash {\rho_{\rm{I}}}]_{+} \\
&  \hspace{20mm} + [h_{\tilde{y}}]P_{\rm{E}}[h(y_{k})+\nu_{k}\oslash {\rho_{\rm{E}}}],\\
\hat{Q}_{\tilde{y}\tilde{y}} &= Q_{\tilde{y}\tilde{y}} +
 [g_{\tilde{y}}]_{\mathcal{A}}^{\tr}P_{\rm{I}}[g_{\tilde{y}}]_{\mathcal{A}} +  P_{\rm{I}}[g(y_{k})+\lambda_{k}\oslash {\rho_{\rm{I}}}]_{+} [g_{\tilde{y}\tilde{y}}]_{\mathcal{A}} \\&\hspace{8mm}+ 
 [h_{\tilde{y}}]^{\tr}P_{\rm{E}}h_{\tilde{y}} + P_{\rm{E}}[h(y_{k})+\nu_{k}\oslash {\rho_{\rm{E}}}]h_{\tilde{y}\tilde{y}},
\end{align*}
with $P_{\rm{I}} = {\rm{diag}}[\rho_{\rm{I}}]$, $P_{\rm{E}} = {\rm{diag}}[\rho_{\rm{E}}]$. Here, the GN approximation can be used to modify $\hat{Q}_{\tilde{y}\tilde{y}}$ as follows
\begin{align}\label{eq:AL_DDP_gauss_newton_approx}
\hat{Q}_{\tilde{y}\tilde{y}} \approx  Q_{\tilde{y}\tilde{y}} +
 [g_{\tilde{y}}]_{\mathcal{A}}^{\tr}P_{\rm{I}}[g_{\tilde{y}}]_{\mathcal{A}} + 
 [h_{\tilde{u}}]^{\tr}P_{\rm{E}}h_{\tilde{u}}.
 \end{align}
By changing $u$ to $\tilde{u}$, DDP can be used for optimization. One key difference is the update law in the forward pass, where both state and control variables are updated by the following linear update law
\begin{align}\label{eq:ALDDP_multiple_variable_update}
    \tilde{u}_{k}^{\rm{new}} = \tilde{u}_{k} + \kappa + K \delta x_{k}, \quad 
    x_{k+1}^{\rm{new}} = \Pi(\tilde{u}_{k}^{\rm{new}}).
\end{align}
The rest of the parameters are updated in the same manner as in the single-shooting case.

\subsection{Interior Point DDP}\label{subsec:IPDDP} In this section, we derive and modify IP DDP based on \cite{Pavlov2021IPDDP}. We consider the single-shooting DDP variation with inequality constraints. In particular, we minimize the Lagrangian over control $u$ and maximize it over multiplier $\lambda$ to compute the optimal value of the original problem \eqref{eq:constrained_OCP_dyn_sys}. Introducing slack variables, and removing the equality constraints from dynamics, we get a new objective 
\begin{align*}
&\qquad \qquad \qquad \quad \min_{U}\max_{\Lambda} \hat{J}(X,U,\Lambda)\\ 
&\text{with} \ \hat{J} =
\sum_{k=1}^{N-1}\left[l(x_{k},u_{k}) + \lambda_{k}^{\tr}g(y_{k})\right]
 + \Phi(x_{N}) + \lambda_{N}^{\tr}g(x_{N}), \\\notag
&{\text{s.t.}} \ g(x_{k}, u_{k}) + s_{k}= 0, \ s_{k} \geq 0, \ \lambda_{k}\geq 0, \  k = 1\cdots N-1, \\\notag
& \ \quad \ g(x_{N}) + s_{N} = 0, \hspace{4mm} s_{N} \geq 0, \ \lambda_{N}\geq 0,  
\end{align*}
where $\Lambda$ is a sequence of multipliers similar to \eqref{eq:XU_long_vec}. Here, we cannot simply add terms to $Q$ due to additional decision variables $\lambda_{k}$ and $s_{k}$. Instead, we define
\begin{align*}
\hat{Q}(x_{k}, u_{k}, \lambda_{k}) = l(x_{k},u_{k}) +\lambda_{k}^{\tr}g(x_{k},u_{k}) + V(x_{k+1}).
\end{align*}
The existence of constraints modifies the derivatives of $Q$, and $\lambda$ introduces new derivatives as 
\begin{align*}
    \hat{Q}_{y} &= Q_{y} + g_{y}^{\tr}\lambda_{k}, \ 
    \hat{Q}_{\lambda} = g(x_{k},u_{k}), \ \hat{Q}_{\lambda \lambda} = O_{w},\\
    \hat{Q}_{yy} &= Q_{yy} + \sum_{i=1}^{w}\lambda_{k,i}[g_{i,yy}], \
    \hat{Q}_{\lambda y} = g_{y}(x_{k},u_{k}). 
\end{align*}
The constraint term in $\hat{Q}_{yy}$ can be excluded for better conditioning, that is, $\hat{Q}_{uu} \approx Q_{uu}$.
In order to derive the backward pass, the optimality condition for $\hat{Q}$ under constraints is considered. 
Partial derivative of quadratic approximation of $\hat{Q}$ with respect to $\delta u_{k}$, first order expansion of complementary slackness, and that of slack variable give
\begin{align}\label{eq:IP_DDP_KKT}
    \hat{Q}_{uu}\delta u_{k} + \hat{Q}_{ux}\delta x_{k} + \hat{Q}_{u\lambda}\delta \lambda_{k} &= -\hat{Q}_{u},\\\notag
    \bar{\Lambda} \delta s_{k} + S\delta \lambda_{k} &= -\bar{\Lambda} s + \mu e,\\\notag
    g_{u}\delta u_{k} + g_{x}\delta x_{k} + \delta s_{k} &= - g(\bar{x}_{k}, \bar{u}_{k}) -s_{k},
\end{align}
where $S = {\rm{diag}}[{s}_{k}]$, $\bar{\Lambda}=
{\rm{diag}}{[\lambda_{k}]}$ and $e$ is given in \eqref{eq:IP_KKT_matrix}.
Solving the system above, we obtain the deviation of decision variables with gains as 
\begin{align}  \label{eq:IP_DDP_delta_variables}
    \delta u_{k} &= \kappa + K\delta x_{k}, \
    \delta s_{k} = \kappa_{s} + K_{s}\delta x_{k},\\\notag 
    \delta \lambda_{k} &= r + R\delta x_{k}.
\end{align}    
By plugging $\delta u_{k}$ and $\delta \lambda_{k}$ into the quadratic expansion of $\hat{Q}$ and mapping terms of $\delta x_{k}$ s with $V_{k}$, derivatives of the value function are obtained by
\begin{align}\label{eq:IPDDP_value_recursion_full}
    V_{x,k+1} &= \hat{Q}_{x} + K^{\tr}\hat{Q}_{u} + R^{\tr}\hat{Q}_{\lambda} + K^{\tr}\hat{Q}_{uu}k \\\notag
    & \quad +\hat{Q}_{x\lambda}r  
     + \hat{Q}_{xu}\kappa + K^{\tr}\hat{Q}_{u\lambda}r + R^{\tr}Q_{\lambda u}\kappa,\\\notag
    V_{xx,k+1} &= \hat{Q}_{xx} + K^{\tr}\hat{Q}_{uu}K + \hat{Q}_{x\lambda}R + R\hat{Q}_{\lambda x}  \\ \notag
    &\quad + \hat{Q}_{xu}K + K\hat{Q}_{xu} +K^{\tr}\hat{Q}_{u\lambda}R + R^{\tr}\hat{Q}_{\lambda u}K.
\end{align}
We note that the approach described here differs from the one in the original work, which excludes the gains of $\lambda_{k}$ from the recursion of the value function and excludes constraints in the final time step. We provide further details about this difference in Appendix \ref{sec:app_IP}.
In the forward pass, $u_{k}, s_{k}$ and $\lambda_{k}$ are updated using the gains in \eqref{eq:IP_DDP_delta_variables} and based on \eqref{eq:IP_fraction_to_boundary} as
\begin{align}\label{eq:IP_DDP_update_lawconda}
\notag
    u^{\rm{new}}_{k} &= u_{k} + \alpha_{s}\kappa + K\delta x_{k}, \
    s^{\rm{new}}_{k} = s_{k} + \alpha_{s}\kappa_{s} + K_{s}\delta x_{k},\\
    \lambda^{\rm{new}}_{k} &= \lambda_{k} + \alpha_{\lambda}r + R\delta x_{k}.
\end{align}
Observe that $u_{k}, s_{k}$ and $\lambda_{k}$ are updated with different step sizes $\alpha_{s}$, $\alpha_{\lambda}\in(0,1]$. The values of $\alpha_{s}$ and $\alpha_{\lambda}$ are determined by the line search filter method so that they satisfy \eqref{eq:IP_fraction_to_boundary}. Although both merit function and filter-based line search approaches can be used in the IP method, the authors of the original work use a filter-based approach similar to the one of the popular IPOPT package \citep{Wachter2006IPOPT}.
Here, a filter is formed by the sum of the violation of the constraint on the trajectory on one axis $\sum_{k=1}^{N}\norm{s_{k} +g_{k}}_{1}$, and the barrier cost $J - \mu\sum_{k=1}^{N}\sum_{i=1}^w\log[s_{k}]_{i}$ on the other axis.
The upper bound of the line search for $\alpha_{s}$, denoted by $\alpha_{s}^{\text{max}}$ and the step size $\alpha_{\lambda}$ are given as 
\begin{align}\label{eq:DDP_fraction_to_boundary}
    \alpha_{s}^{\text{max}} &\coloneqq \max \Set{ \alpha \in (0,1] | s_{k}^{\text{new}} \geq (1-\tau )s_{k} }, \\\notag 
    \alpha_{\lambda} &\coloneqq \max \Set{ \alpha \in (0,1] | \lambda_{k}^{\text{new}} \geq (1-\tau )\lambda_{k} }, \quad  
\end{align}
for $k = 1\cdots N-1$. The step size for $\lambda$ is obtained without the line search because it does not affect the filter.  Unfortunately, applying this rule to DDP is not as straightforward as applying a similar rule to IP for static problems due to the existence of feedback terms in DDP.
Because search directions are defined via the linear feedback equations in DDP as \eqref{eq:IP_DDP_update_lawconda}, $\alpha_{\lambda}$ is affected by $\alpha_{s}$ through $\delta x_{k}$. In the original work, without setting the max. values of each $\alpha$ s, the line search is performed with a common parameter for $\alpha_{s}$ and $\alpha_{\lambda}$, making the line search unnecessarily conservative. This is because each step size is affected by that of the paired variable. To alleviate this problem, we propose applying the rule with different $\alpha$ for $s$ and $\lambda$ as in the original formulation of the IP method. We first find $\alpha_{s}^{\rm{max}}$ using dynamics. Since there is no closed-form solution for $\alpha_{s}^{\rm{max}}$, we rely on a line search to find $\alpha_{s}^{\rm{max}}$. Next, the line search filter is performed, setting $\alpha_{s}^{\rm{max}}$ as the upper bound. Inside of this line search, the inner line search for $\alpha_{\lambda}$ is also performed, using $\delta x_{k}$ generated with $\alpha_{s}$ from the outer line search. The IP DDP algorithm used in section \ref{sec:experiments} is implemented with this modified line search filter. Finally, we note that the feasibility restoration mechanism explained in \ref{subsec:IP} is not implemented with DDP.

\subsection{Primal-Dual Augmented Lagrangian DDP}
\label{subsec:PDAL DDP}
This section shows the derivation of single- and multiple- shooting PDAL DDP.
\subsubsection{Single-shooting}
Following \eqref{eq:PDAL}, and using $\rho = \frac{1}{\mu}$, we have PDAL, which is the objective of PDAL DDP as
\begin{align*}
    &\min_{U, \Lambda}{\hat{J}}(X,U,\Lambda) \\
    = &\min_{U,\Lambda} \sum_{k=1}^{N-1} \left[
    l(x_{k},u_{k}) + \lambda_{\mathrm{e},k}^{\tr}\Big[[G^{\dagger}_{k}]_{+}-P_{k}^{-1}\Big(\lambda_{\mathrm{e},k} - \frac{\lambda_{k}}{2}\Big)\right.\\
    & \ \left.+\sum_{i=1}^{w}\frac{\rho_{k,i}}{2}\norm{\big[G^{\dagger}_{k,i}\big]_{+} - \frac{1}{\rho_{k,i}}\Big([\lambda_{\mathrm{e},k}]_{i} - \frac{\lambda_{k,i}}{2} \Big)}^{2} \right.\\
    & \ \ \left.+\sum_{i=1}^{w}\frac{\rho_{k,i}}{2}\norm{[G^{\dagger}_{k,i}]_{+} - \frac{\lambda_{k,i}}{2\rho_{k,i}}}^{2}
    \right] \\
    &\quad  + \Phi(x_{N}) + 
    \lambda_{\mathrm{e},N}^{\tr}\Big[[
    G^{\dagger}_{N}]_{+}-P_{N}^{-1}\Big(\lambda_{\mathrm{e},N} - \frac{\lambda_{N}}{2}\Big)\Big] \\
    &\quad \ +\sum_{i=1}^{w}\frac{\rho_{N,i}}{2}\norm{[G^{\dagger}_{N,i}]_{+} - \frac{1}{\rho_{N,i}}\left([\lambda_{\mathrm{e},N}]_{i} - \frac{\lambda_{N,i}}{2} \right)}^{2}\\
    &\quad \ \ +\sum_{i=1}^{w}\frac{\rho_{N,i}}{2}\norm{[G^{\dagger}_{N,i}]_{+} - \frac{\lambda_{N,i}}{2\rho_{N,i}}}^{2},\\
&\text{with} \ G^{\dagger}_{k} =  g(x_{k},u_{k}) + P_{k}^{-1}(\lambda_{\mathrm{e},k}-{\lambda_{k}}/{2}),\\
& \quad \quad P_{k} = {\rm{diag}}[\rho_{k}].
\end{align*}
This objective modifies $Q$ and its derivatives as follows.
\begin{align}\label{eq:Quu_PDALDDP}
\notag
    &\hat{Q}_{y}= 
    Q_{y}  +  [g_{y}]_{\mathcal{A}}^{\tr}[2\pi_{\rm{I}}-\lambda_{k}]_{+}, \\
    \notag
    & \hat{Q}_{\lambda} = -\frac{1}{2}P^{-1}\big[[2\pi_{\rm{I}} - \lambda_{k}]_{+} - \lambda_{k}\big],\ \hat{Q}_{\lambda y} = -[g_{y}]_{\mathcal{A}}, \\
    &\hat{Q}_{yy} = 
    Q_{yy} + [2\pi_{\rm{I}}-\lambda_{k}]_{+}[g_{yy}]_{\mathcal{A}} + 2[g_{y}]^{\tr}_{\mathcal{A}}P[g_y]_{\mathcal{A}}, \\
    \notag
    &\hat{Q}_{\lambda \lambda} = \frac{P^{-1}}{2}\big[[I_{w}]_{\mathcal{A}} + I_{w}\big] = \frac{{\rm{diag}}[\mu_{\rm{I}}]}{2}\big[[I_w]_{\mathcal{A}} + I_w\big],\\
    \notag
    &\text{where} \hspace{8mm} \pi_{\rm{I}} = P g(x_{k},u_{k}) + \lambda_{\rm{e}},
\end{align}
and projection in \eqref{eq:positive_projection_PDAL}. $[I_{w}]_{\mathcal{A}}$ is a modified identity matrix whose $i$ th diagonal element is 1 if $g_{i}$ is active (positive after projection in \eqref{eq:positive_projection_PDAL}) and zero otherwise.
First-order optimality condition for quadratic approximation of $\hat{Q}$ gives
\begin{align}\label{eq:PDAL_DDP_single_opt_matrix}
    \begin{bmatrix}
    \hat{Q}_{uu} & \hat{Q}_{u\lambda}\\
    \hat{Q}_{\lambda u} & \hat{Q}_{\lambda \lambda}
    \end{bmatrix}
    \begin{bmatrix}
    \delta u_{k} \\ \delta \lambda_{k}
    \end{bmatrix}
    = 
    -\begin{bmatrix}
    \hat{Q}_{u}\\ \hat{Q}_{\lambda}
    \end{bmatrix}
    - \begin{bmatrix}
    \hat{Q}_{ux} \\ \hat{Q}_{\lambda x}
    \end{bmatrix} \delta x_{k}
\end{align}
The matrix on the left-hand side is ill-conditioned as $\rho_{\rm{I}}$ ($\mu_{\rm{I}}$) gets large (small) due to the last term in $\hat{Q}_{yy}$, which is avoided by a similar transformation mentioned in section \ref{subsubsec:min_of_PDAL}. See appendix \ref{sec:app_pdal_transform} for details. 
Using this transformation, we have a transformed symmetric system
\begin{align*}
    &\begin{bmatrix}
    H_{\rm{s}} & -[g_{u}]_{\mathcal{A}}^{\tr}\\ -[g_{u}]_{\mathcal{A}}^{\tr} &  -M_{\mu}
    \end{bmatrix}
    \begin{bmatrix}
    \delta u_{k} \\ -\delta \lambda_{k}
    \end{bmatrix}
    = \\
    & -\begin{bmatrix}
    Q_{u} + [g_{u}]_{\mathcal{A}}^{\tr}\lambda \\ \hat{Q}_{\lambda}
    \end{bmatrix}
    - 
    \begin{bmatrix}
    Q_{ux} + [2\pi_{\rm{I}}-\lambda_{k}]_{+}[g_{ux}]_{\mathcal{A}} \\ \hat{Q}_{\lambda x}
    \end{bmatrix}
    \delta x_{k},
\end{align*}
with $H_{\rm{s}} = Q_{uu} + [2\pi_{\rm{I}}-\lambda]_{+}[g_{uu}]_{\mathcal{A}}$ and $M_{\mu} = \frac{\rm{diag}[\mu]}{2}\big[[I_w]_{\mathcal{A}} + I_{w}\big]$. $H_{\rm{s}}$ here implies Hessian since it can be seen as the second derivative of Lagrangian formed by the objective $Q$, multiplier $[2\pi_{\rm{I}}-\lambda]_{+}$, and constraints $g$. 

When the constraints are not active, the system
gives the same solution as normal uncostrained DDP for $\delta u_{k}$. For $\delta \lambda_{k}$, the solution makes $\lambda_{k}^{\rm{new}}$ zero. We use a similar update law proposed in \cite{Jallet2022PDALDDP}, which makes $\lambda_{k}$ strictly zero for inactive constraints as below.
\begin{align*}
    u^{\rm{new}}_{k} &= u_{k} + \alpha \kappa_{k} + K_{k}\delta x_{k}, \\ 
    \lambda^{\rm{new}}_{k} &= \big[\lambda_{k}  -[\lambda_{k}]_{\mathcal{A}^{c}} + \alpha [r_{k}]_\mathcal{A} +  [R_{k}]_{\mathcal{A}}\delta x_{k}\big]_{+},
\end{align*}
where $\mathcal{A}^{c}$ denotes the inactive set of constraints. The positive orthant projection for $\lambda$ is ensuring that $\lambda$ is nonnegative. The gains lead to recursion for the value function as below.
\begin{align*}
V_{x,k} &= \hat{Q}_{x} + K^{\tr}\hat{Q}_{u} + R^{\tr}\hat{Q}_{\lambda} + K^{\tr}\hat{Q}_{uu}\kappa  \\\notag
    &+ R^{\tr}\hat{Q}_{xu}r + \hat{Q}_{xu}\kappa + \hat{Q}_{x\lambda}r  + K^{\tr}\hat{Q}_{u\lambda}r + R^{\tr}\hat{Q}_{\lambda u}\kappa,\\
    V_{xx,k} &= \hat{Q}_{xx} + K^{\tr}\hat{Q}_{uu}K + R^{\tr}\hat{Q}_{\lambda\lambda}R  +  K^{\tr}\hat{Q}_{ux}  \\\notag + & \hat{Q}_{xu}K + R^{\tr}\hat{Q}_{\lambda x} + \hat{Q}_{x\lambda} R +  K^{\tr}\hat{Q}_{u\lambda}R + R^{\tr}\hat{Q}_{\lambda u} K.
\end{align*}
After performing the DDP presented above, $\lambda_{\mathrm{e}}$ is updated by the law \eqref{eq:multiplier_update_PDAL}. The penalty parameters and tolerances for constraint satisfaction and DDP are updated in the same way as in the AL DDP provided in section \ref{sec:AL_DDP}.

A key advantage of PDAL over the standard AL method is its robustness to changes in the penalty parameter \citep{Robinson2014PDAL}. When the parameter is updated to improve constraint satisfaction, a new objective is defined, and a corresponding search direction is computed. By accounting for the interaction between the penalty parameter and the dual variables, PDAL can provide a more effective search direction than standard AL following the updates.

\subsubsection{Multiple-shooting}
Due to the limited space, we omit the form of an objective and recursion. Following the same the procedure as AL DDP and single shooting PDAL, we have
\begin{align*}
    &\hat{Q}_{\tilde{y}}= 
    l_{\tilde{y}} + h_{\tilde{y}}^{\tr}[2\pi_{\rm{E}}-\nu_{k}] +  [g_{\tilde{y}}]_{\mathcal{A}}^{\tr}[2\pi_{\rm{I}}-\lambda_{k}]_{+}, \\
  &\hat{Q}_{\lambda} = -\frac{1}{2}P_{\rm{I}}^{-1}\big[[2\pi_{\rm{I}}-\lambda_{k}]_{+}-\lambda_{k}\big],\ 
    \hat{Q}_{\nu}  = - P_{\rm{E}}^{-1}[\pi_{\rm{E}} -\nu_{k}],\\  
    &\hat{Q}_{\tilde{y}\tilde{y}} = 
    Q_{\tilde{y}\tilde{y}} + [2\pi_{\rm{E}}-\nu_{\rm{e}}]h_{\tilde{y}\tilde{y}}+2h_{\tilde{y}}^{\tr}{P_{\rm{E}}}[h_{\tilde{y}}] \\
    & \quad \quad +[2\pi_{\rm{I}}-\lambda_{\rm{e}}]_{+}[g_{\tilde{y}\tilde{y}}]_{\mathcal{A}} \quad + 2[g_{\tilde{y}}]^{\tr}_{\mathcal{A}}{P_{\rm{I}}}[g_{\tilde{y}}]_{\mathcal{A}},\\
    &\hat{Q}_{\nu\nu} = {P_{\rm{E}}}^{-1}, \ 
    \hat{Q}_{\lambda \lambda} = \frac{P_{\rm{I}}^{-1}}{2}\big[[I_{w}]_{\mathcal{A}} + I_{w}\big], \\
    & \hat{Q}_{\nu \tilde{y}} = -h_{\tilde{y}}, \quad 
    \hat{Q}_{\lambda \tilde{y}} = -[g_{\tilde{y}}]_{\mathcal{A}}, \quad 
    \\
    &\text{where}  \ P_{\rm{E}} = {\rm{diag}}{[\rho_{\rm{E}}]}, \ \pi_{\rm{E}} = P_{\rm{E}}h(x_{k},u_{k},x_{k+1}) + \nu_{\rm{e}},  \\
    & \quad\quad\quad   P_{\rm{I}} = {\rm{diag}}{[\rho_{\rm{E}}]}, 
    \ 
    \pi_{\rm{I}} = P_{\rm{I}}g(x_{k},u_{k}) + \lambda_{\rm{e}}.
\end{align*}
Using the multiple-shooting DDP with these $\hat{Q}$ followed by updating multipliers, penalty parameters, and tolerances for inner DDP and constraints, the problem is solved by multiple-shooting PDAL DDP.

\subsection{ADMM DDP}
This section introduces the ADMM-based variaton of constrained DDP following \cite{2017SindwaniADMMDDP}.
\subsubsection{Introducing Copy Variables }
Let us introduce a copy of the variables of $X$ and $U$, denoted by $X^{\mathrm{c}}$, and $U^{\mathrm{c}}$, respectively. These copy variables are intended to strictly satisfy the additional state/control constraints, except for the dynamics ones. Original variables ($X, U$) minimize the original cost only under dynamic constraints. We assume that the constraints can be divided into state- and control-dependent parts as follows. 
\begin{align*}
g(x, u) = [g^{x}(x), g^{u}(u)]^{\tr}. 
\end{align*}
To encode strict feasibility in the cost, we use an indicator function defined as:
\begin{align*}
\mathbbm{1}_{S}(x) =  \begin{cases}
0, \quad x \in S,\\
\infty, \quad x \notin S. 
\end{cases}
\end{align*}
Hence, the problem can be reformulated as follows
\begin{align*}
&\min_{X,U,X^{\mathrm{c}},U^{\mathrm{c}}} \hat{J}(X,U,X^{\mathrm{c}}, U^{\mathrm{c}}) 
\quad \text{s.t.} \ X= X^{\mathrm{c}}, \ U = U^{\mathrm{c}}, 
\end{align*}
where 
\begin{align*}
\hat{J}(\underbrace{X,U}_{Y},\underbrace{X^{\mathrm{c}}, U^{\mathrm{c}}}_{Y^{\mathrm{c}}}) &= J(Y) + \mathbbm{1}_{f}(Y) + \mathbbm{1}_{g^{u}}(U^{\mathrm{c}}) + \mathbbm{1}_{g^{x}}(X^{\mathrm{c}}),\\
\mathbbm{1}_{f}(Y)
&= \sum_{k=1}^{N-1}\mathbbm{1}_{x_{k+1}-f(x_{k}, u_{k})}(x_{k}, u_{k}, x_{k+1}), \\
\mathbbm{1}_{g^{u}}(U^{\mathrm{c}}) = \sum_{k=1}^{N-1} & \mathbbm{1}_{g^{u}\leq0}(u^{\mathrm{c}}_{k}), \
\mathbbm{1}_{g^{x}}(X^{\mathrm{c}}) = \sum_{k=1}^{N}\mathbbm{1}_{g^{x}\leq0}(x^{\mathrm{c}}_{k}). 
\end{align*}
Here, we use $Y = [X^{\tr},U^{\tr}]^{\tr} \in \mathbb{R}^{n+m}$, which is similar to $y_{k}$.
This formulation is equivalent to setting $x = [X,U]$, $z = [X^{\mathrm{c}},U^{\mathrm{c}}]$, $A=I, B=-I$, and $c=0$ in \eqref{eq:admm_problem}. The problem’s AL function is given by
\begin{align*}
\mathcal{L}_{\mathrm{A}}
(Y, Y^{\mathrm{c}}, \Lambda) &= \hat{J}(
Y, Y^{\mathrm{c}}
) \\
& + \frac{1}{2}\sum_{k=1}^{N}\sum_{i=1}^{n_{y}}\rho_{i}(y_{k,i}-y^{\mathrm{c}}_{k,i} + {\lambda_{k,i}/{\rho_{i}}})^{2},
\end{align*}
where we use the same vector $\rho$ across time.

\subsubsection{ADMM DDP}
The ADMM DDP algorithm consists of the following updates that happen in a sequential manner. The variables $X,U$ are first updated by solving the following minimization subproblem:
\begin{align}\label{eq:ADMM_DDP_primal}
X^{\rm{new}}, U^{\rm{new}} = \argmin_{X,U} \mathcal{L}_{\mathrm{A}}(X,U, X^{\mathrm{c}}, U^{\mathrm{c}}, \Lambda)
\end{align}
Since DDP can strictly satisfy constraints from dynamics, the problem is solved by DDP with modified $Q$ functions and its derivatives:
\begin{align}
\label{eq:PDALDDP_Qs}
\hat{Q} &= Q + \frac{1}{2}\sum_{i=1}^{n_{y}}\rho_{i}(y_{i}-y^{c}_{i} + \lambda_{i} / \rho_{i})^{2} \\
\notag
\hat{Q}_{y} &= Q_{y} + \mathrm{diag}[\rho](y-y^{c} + \lambda \oslash \rho), \\
\label{eq:Quu_ADMMDDP}
\hat{Q}_{yy} &= Q_{yy} + \rm{diag}[\rho].
\end{align}
Using the $Y^{\rm{new}}$ obtained by DDP, the copy variables are updated by 
\begin{align*}
 Y^{\mathrm{c},\rm{new}}= \argmin_{Y^{c}}\mathcal{L}_{\mathrm{A}}(Y^{\rm{new}}, Y^{\mathrm{c}}, \Lambda), 
\end{align*}
which can be decomposed for each time instant as  
\begin{align}\label{eq:ADMM_DDP_copy}
y^{\mathrm{c},\rm{new}}_{k} = \argmin_{y_{k}^{c}} \norm{y_k^{\rm{new}}-y_{k}^{c} + \lambda_{k}\oslash \rho_{k}}^{2}, \\\notag
\text{s.t.} \ \ g(y^{\mathrm{c}}) \leq 0, \ k=1,\cdots,N.
\end{align}
This requires solving an optimization problem with a quadratic objective under constraints $g(y^{\mathrm{c}})\leq 0$. This optimization is performed without considering dynamics, and the resulting trajectory may be dynamically infeasible. 
As a special case, when the constraint has the simple form of $y^{\mathrm{c}} \leq y_{\rm{b}}$, the problem is solved by clamping. 
Finally, the multiplier is updated by 
\begin{align}\label{eq:ADMM_DDP_dual}
\lambda_{k}^{\rm{new}} = \lambda_{k} + \rho_{k}(y_{k}^{\rm{new}}-y_{k}^{c,\rm{new}}).
\end{align}
ADMM DDP repeats the three update processes in \eqref{eq:ADMM_DDP_primal}, \eqref{eq:ADMM_DDP_copy}, and \eqref{eq:ADMM_DDP_dual} until the residuals mentioned in section \ref{subsec:ADMM} become small enough. Upon convergence, the original and copy variables will reach to consensus, and as a result, the final solution will be optimal while satisfying all constraints. It is well known however, that ADMM might require many iterations until reaching high accuracy \cite{Boyd2011ADMM}. 

\subsection{Analysis of AL-based DDPs }\label{sec:AL_based_DDP_analysis}

In this section, we analyze AL-based DDPs. Specifically, we analyze the difference between ADMM DDP and others, including AL and PDAL DDPs. As in the case of ADMM for static problems \citep{Boyd2011ADMM}, ADMM DDP also requires long iterations to achieve an accurate solution. This can be understood by investigating how the $Q$ function of DDP captures information on cost and constraints. In AL and PDAL DDP, the active constraints are directly captured in the $Q$ functions of DDP. See \eqref{eq:ALDDP_Qs} for AL DDP and \eqref{eq:PDALDDP_Qs} for PDAL DDP. This information of constraints enables algorithms to satisfy them effectively while reducing the original cost. However, in ADMM DDP, the $Q$ function has information on constraints only through the distance from safe copies (see \eqref{eq:Quu_ADMMDDP}). In addition, safe copies may not be dynamically feasible, making the problem in \eqref{eq:ADMM_DDP_primal} difficult. A canonical example is when the DDP for solving problem \eqref{eq:ADMM_DDP_primal} commands a control sequence that exceeds its limits in many time steps. In this situation, staying close to the safe copies (clamped control) and completing the task conflict with each other, slowing down cost reduction and constraint satisfaction. Indeed, in our experiment in section \ref{subsec:result_one_guess}, we observe that ADMM DDP cannot handle problems where the control constraints are tight, and a control sequence needs to hit its limit in many time steps. 

Another difference can be found in the role of the penalty parameter $\rho$. In AL and PDAL DDPs, the product of $\rho$ and active constraints modifies $Q$ functions. Therefore, when the constraints are not active, $\rho$s does not affect $Q_{uu}$. In ADMM DDP, $\rho$ is added to $Q_{uu}$ regardless of the status of the constraints as if the regularizer in \eqref{eq:DDP_regularization}. This seems appealing for the conditioning of $Q_{uu}$, but too large $\rho$ can slow down optimization, as mentioned in section \ref{sec:DDP}. In AL and PDAL DDPs, too large $\rho$ also interrupts optimization, but what we would like to emphasize here is that in ADMM, $\rho$ always affects $Q_{uu}$. 

\section{SQP for dynamical systems}\label{sec:dynamic_SQP}
This section presents a concise overview and derivation of SQP for dynamical systems, considering both the single- and multiple-shooting approaches based on \cite{Gill2000SQPdynamical}.

Consider the constrained optimal control problem
\eqref{eq:constrained_OCP_dyn_sys}. Here, we have the sequence of state and control as long vectors as in \eqref{eq:XU_long_vec} and a deviated trajectory as in a similar manner as in section \eqref{sec:DDP}, i.e., $X = \bar{X} + \delta X, \ U = \bar{U} + \delta U$. Both single- and multiple-shooting SQP have linearized dynamics as constraints. The difference is whether the constraints are implicit or explicit. The constraints from the dynamics are linearized as below. 
\begin{align}\label{eq:SQP_linearized_dyn_multi}
  &\bar{x}_{k+1} + \delta x_{k+1} = f(\bar{x}_{k}, \bar{u}_{k}) + A_{k}\delta x_{k} + B_{k}\delta u_{k} \\\notag
  {\text{with}}& \ A_{k} = f_{x,k}, \quad B_{k} = f_{u,k}, \quad k=1,\cdots N-1.
\end{align}
In the single-shooting formulation, the state is updated via
the system dynamics, satisfying the equality constraints on the nominal trajectory. Therefore, the first terms of both sides of the equations cancel out. Therefore,  the deviation of state and control are tied with a matrix $F \in \mathbb{R}^{nN \times m(N-1)}$ by
\begin{align}\label{eq:SQP_single_eq_cnst} 
&\hspace{30mm}\delta {X} = F \delta U,\\\notag
&\text{where,}\\\notag
&F = 
    \begin{bmatrix}
    O_{n,m}         &  O_{n,m}         & \cdots    & O_{n,m}\\
    B_{1}           & O_{n,m}          & \cdots    & O_{n,m}\\
    A_{2}B_{1}      & B_{2}            & \cdots    & O_{n,m}\\
    \vdots          & \vdots           &\ddots     & \vdots\\
    A_{N-1}\cdots{A}_{2}{B}_{1} &
    {A}_{N-1}\cdots{A}_{3}{B}_{2} &
    \cdots & {B}_{N-1}
    \end{bmatrix}.
\end{align}
In the multiple-shooting case, however, the equality constraint from dynamics might be violated. Consequently, the first terms of both sides of \eqref{eq:SQP_linearized_dyn_multi} might not cancel out. Thus, the constraints on each time step take the following form. 
\begin{align*}
     \delta x_{1} &= x_{\rm{init}}- \bar{x}_{1}, \\
     \delta x_{k+1} &= A_{k}\delta x_{k} + B_{k}\delta{u}_{k}-\bar{x}_{k} + f(\bar{x}_{k},\bar{u}_{k}),
\end{align*}
$k=1, \cdots N-1$.
$x_{\rm{init}}$ is a given initial state where the control sequence cannot affect, and thus $\bar{x}_{1}-x_{\rm{init}}=0$. The equality constraints with nominal terms are given in matrix-vector form by 
\begin{align}\label{eq:SQP_multi_eq_cnst}
  \hat{F} & (\bar{X},\bar{U})  + \hat{F}_{Y} 
    \delta Y = 0, \\\notag \text{with} \
    \hat{F}(\bar{X}, \bar{U}) &= 
 \begin{bmatrix}
 \bar{x}_{1} -x_{\rm{init}}\\
 \bar{x}_{2}  -f(\bar{x}_{1},\bar{u}_{1})\\
 \vdots\\
 \bar{x}_{N}  -f(\bar{x}_{N-1},\bar{u}_{N-1})
 \end{bmatrix}, \\\notag
\hat{F}_{Y} & = 
\begin{bmatrix}
\hat{F}_{Y1} & \hat{F}_{Y2} 
\end{bmatrix},\\\notag
\text{where} \  
\hat{F}_{Y1} &=  
\begin{bmatrix}
    I_{n}  & O_{n}  & \cdots  &          & O_{n}\\
    -A_{1} & I_{n}  & O_{n}   & \cdots   & O_{n}\\
    O_{n}  & -A_{2} & I_{n}   & O_{n}      & O_{n}\\
    \vdots & \cdots & \ddots  &  \ddots    & \vdots\\
    O_{n}  & \cdots        & O_{n} & -A_{N-1} & I_{n}
\end{bmatrix},\\\notag
\hat{F}_{Y2}
    &= \begin{bmatrix}
    O_{n,m}       &        & & O_{n,m}\\
    -B_{1}        & O_{n,m}       & & O_{n,m} \\
    O_{n,m}         & -B_{2} &O_{n,m} & O_{n,m}\\
    \vdots &  \cdots& \ddots    & \vdots \\
    O_{n,m}& \cdots & O_{n,m}& -B_{N-1}
 \end{bmatrix},
\end{align}
with $Y = [X^{\tr}, U^{\tr}]^{\tr}\in \mathbb{R}^{nN + m(N-1)}$ and $\hat{F}\in \mathbb{R}^{nN \times [nN +m(N-1)]}$. 
\subsection{Single-shooting SQP}\label{sec:SQP_single}
From \eqref{eq:SQP_QP}, the SQP subproblem of the problem in \eqref{eq:constrained_OCP_dyn_sys} is given by 
\begin{align}
    &\min_{\delta X, \delta U} [J_{X}^{\tr}, J_{U}^{\tr}]
    \begin{bmatrix}
    \delta X \\ \delta U
    \end{bmatrix}
    + \frac{1}{2}
    \begin{bmatrix}
    \delta X \\ \delta U
    \end{bmatrix}^{\tr}
    \begin{bmatrix}
    \nabla_{XX} \mathcal{L} & \nabla_{XU} \mathcal{L}\\
    \nabla_{UX} \mathcal{L} & \nabla_{UU} \mathcal{L}
    \end{bmatrix}
    \begin{bmatrix}
    \delta X \\ \delta U
    \end{bmatrix},\label{eq:SQP_single_objective_0}\\ \notag
    &{\text{s.t.}}\
    \eqref{eq:SQP_single_eq_cnst},  \text{ \ and \ } 
    G(\bar{X}, \bar{U}) + [G_{X}, G_{U}]
    \begin{bmatrix}
    \delta X \\ 
    \delta U
    \end{bmatrix} \leq 0,
\end{align}
where the constraints 
\begin{equation*}
G(X,U) = [g(x_{1},u_{1})^{\tr}\dots g(x_{N})^{\tr}]^{\tr} \in \mathbb{R}^{wN}
\end{equation*}
and Lagrangian $\mathcal{L} = J + \lambda_{\rm{v}}^{\tr} G(\bar{X},\bar{U})$
with a vectorized multiplier $\lambda_{\rm{v}}\in \mathbb{R}^{wN}$. 
Gradients of $G(\bar{X}, \bar{U})$ are given by
\begin{align*}
    G_{X} &= {\rm{blkdiag}}[
    g_{x}(\bar{x}_{1}, \bar{u}_{1}), \cdots, g_{x}(\bar{x}_{N-1}, \bar{u}_{N-1}), g_{x}(\bar{x}_{N})],\\\nonumber
    G_{U} &=\begin{bmatrix}
    {\rm{blkdiag}}[
    g_{u}(\bar{x}_{1}, \bar{u}_{1}), \cdots, g_{u}(\bar{x}_{N-1}, \bar{u}_{N-1})]\\ O_{w, m(N-1)}
    \end{bmatrix},
\end{align*}
whose dimension are $\mathbb{R}^{wN \times nN}$ and $\mathbb{R}^{wN \times m(N-1)}$, respectively.
For the cost $J$, the derivatives are given by
\begin{align*}
    J_{X} &= [l_{x}(\bar{x}_{1}, \bar{u}_1)^{\tr}, \cdots, l_{x}(\bar{x}_{N-1}, \bar{u}_{N-1})^{\tr}, \Phi_{x}(\bar{x}_{N})^{\tr}]^{\tr}, \\\notag
    J_{U} &= [l_{u}(\bar{x}_{1}, \bar{u}_1)^{\tr}, \cdots, l_{u}(\bar{x}_{N-1}, \bar{u}_{N-1})^{\tr}]^{\tr},
\end{align*}
with $J_{X}\in\mathbb{R}^{nN}$ and $J_{U}\in\mathbb{R}^{m(N-1)}$. 
By eliminating $\delta X$ using the equality constraints \eqref{eq:SQP_single_eq_cnst}, single-shooting SQP for \eqref{eq:constrained_OCP_dyn_sys} is formulated as
\begin{align}\label{eq:SQP_single}
    &\min_{\delta U} \ [J_{X}^{\tr}F + J_{U}^{\tr}]\delta U + \frac{1}{2}\delta U^{\tr}[F^{\tr}\nabla_{XX}\mathcal{L}F \\ \notag 
    &\quad \quad \quad \quad + 2F^{\tr}\nabla_{XU}\mathcal{L} + \nabla_{UU} \mathcal{L}]\delta U\\\notag
    &{\text{s.t.}}\quad
    G(\bar{X}, \bar{U}) + (G_{X}F + G_{U})\delta U \leq 0.
\end{align}
 The $X$ trajectory is updated using $\delta U^{\ast}$ which is the solution of \eqref{eq:SQP_single}, and the dynamics of the system with an appropriate step size similar to the DDP in  \eqref{eq:DDP_update_control} but in open-loop fashion. We perform a line search with the AL merit function to determine the step size, whose detail is given in appendix \ref{sup_sec:SQP}. When an appropriate step size cannot be found in the line search, we regularize the Hessian as in the case of DDP  \eqref{eq:DDP_regularization} and resolve the QP subproblem. We note that computing the exact Hessian of \eqref{eq:SQP_single} can be expensive and may not even be worth computing when it is not PD. There exists an iterative Hessian approximation scheme, which is known as the BFGS update \citep{Broyden1970BFGS, Fletcher1970BFGS, Goldfarb1970BFGS, Shanno19070BFGS}. This update rule can approximate a PD Hessian based on the Hessian in the previous iteration.

\subsection{Multiple shooting-SQP}\label{subsec:SQP_multi}
Using equality constraints in \eqref{eq:SQP_multi_eq_cnst}, instead of \eqref{eq:SQP_single_eq_cnst}, and augmented variable $Y$, the multiple-shooting formulation is obtained as
\begin{align}\label{eq:SQP_multi}
    &\hspace{8mm} \min_{\delta Y} \ [J_{Y}^{\tr}]\delta Y + \frac{1}{2}\delta Y^{\tr}[\nabla_{YY}\mathcal{L}]\delta Y\\\notag
    &{\text{s.t.}}\quad
    \eqref{eq:SQP_multi_eq_cnst}, \text{ and }
    G(\bar{X}, \bar{U}) + (G_Y)\delta Y \leq 0, 
\end{align}
where the Lagrangian is given by 
 $
 \mathcal{L}= J + \lambda_{\rm{v}}^{\tr}G(X,U)+ \nu_{\rm{v}}^{\tr}\hat{F}(X,U), 
 $
with Lagrangian multipliers $\lambda_{\rm{v}} \in \mathbb{R}^{wN}$ and $\nu_{\rm{v}} \in \mathbb{R}^{nN}$. After solving the problem above, $Y$ (both $X$ and $U$) is updated by the following linear update law:
\begin{equation*}
Y = \bar{Y} + \alpha \delta Y,
\end{equation*}
 whose step size $\alpha$ is determined by the line search with the AL merit function as in the single-shooting method. Another update strategy that strictly satisfies the dynamics is possible, as presented in \cite{Tenny2004SQPMPC}. Here, the jacobians from the dynamics given by $\hat{F}_{Y}$ in \eqref{eq:SQP_multi_eq_cnst} has a sparse structure due to the recursion in the dynamics. Several methods that can exploit this structure and reduce computational complexity of SQP have been proposed \citep{dohrmann1997efficient, Rao1998IPMPC, JORGENSEN2004efficientMPC,
Wang2010linearMPCefficientT}. Note that the sparse structure is available only in the multiple-shooting formulation.

\section{Theoretical Time Complexity}\label{sec:time_complexity}
In this section, we compare the theoretical per-iteration time complexities of DDP and SQP. 
We note that for SQP, the reported complexity corresponds to one iteration of solving the inner QP subproblem, rather than the overall solution of the outer nonlinear problem. This is because with inequality constraints, we cannot tell the number of iterations required to solve the inner subproblem. All complexities are summarized in Table \ref{tab:intro_comparison}.

In unconstrained DDP, the computational bottleneck is inverting $Q_{uu}$. Reference \cite{yakowitz1991convergence} has a detailed breakdown of the computational complexity of unconstrained DDP, which includes matrix inversion, as well as matrix multiplication operations.  Matrix multiplication is a highly parallelizable operation, and therefore it can be optimized very effectively, while matrix inversion remains a much harder operation to optimize and accelerate. When $n \gg m$, matrix multiplications involving $Q_{xx}$ and $V_{xx}$, are theoretically more expensive than inverting $Q_{uu}$. However, these operations can be accelerated, whereas inversion of $Q_{uu}$ still requires a cubic time complexity. The inversion is performed $N-1$ times in the backward pass.
Therefore, the complexity is $\mathcal{O}(Nm^{3})$, in single-shooting and $\mathcal{O}(N(n+m)^{3})$ in multiple-shooting. The barrier, AL, and ADMM DDPs follow the same complexity because the number of decision variables (of inner loops) is the same as in standard DDP. In IP DDP, the coefficient matrix of the system of equations in 
\eqref{eq:IP_KKT_matrix} has size $(m+2w) \times (m+2w)$. This system can be reduced to a smaller one by eliminating the slack variable, whose coefficient matrix has the size of $(m+w)\times (m+w)$ \citep{Wachter2006IPOPT}. Further reduction can be performed by eliminating the multiplier, giving a system with a coefficient matrix of $m \times m$ \citep{Nocedal2006numerical}. Therefore, the time complexity with respect to the decision variables is cubic. The same analysis can be applied to PDAL DDPs.

In SQPs, the bottleneck is solving the QP subproblems, whose Hessian has the size of $(N-1)m \times (N-1)m$ in single- and $\{(N-1)m + Nn\}  \times \{(N-1)m + Nn\}$ in multiple-shooting. These large Hessians lead to $\mathcal{O}(N^{3}m^{3})$ and $\mathcal{O}(N^{3}(m+n)^{3})$ for single- and multiple-shooting, respectively for one iteration of QP subproblem. However, as mentioned in section \ref{subsec:SQP_multi}, several methods in structure-exploiting multiple-shooting SQP that can reduce the QP's complexity are available. 

For methods based on Riccati recursion \citep{Rao1998IPMPC}, the matrix to be inverted has size $m \times m$. 
Hence, the complexity is 
$\mathcal{O}(N m^{3})$,
following the same reasoning as in DDP. 
For methods that uses factorization \citep{Wang2010linearMPCefficientT}, the complexity can be reduced to 
$\mathcal{O}(N (n+m)^{3})$. In single-shooting, dependency on $N$ can reduced to not linear, but quadratic $\mathcal{O}(N^{2}n^{3})$ \citep{kouzoupis2018recent}. From a computational point of view, the single-shooting variant is useful for a problem with a short time horizon $N$ and a large state dimension $n$.

With a long time horizon $N$, DDP and multiple-shooting SQP are favorable due to the linear growth of complexity in $N$. 
Among these, single-shooting DDPs and multiple-shooting SQP with Riccati recursion are the most efficient because of the small number of decision variables. Their advantages are especially important in underactuated systems, e.g., systems with $n\gg m$, which are the most typical ones in robotics.
Although multiple-shooting SQP can achieve the same time complexity, this corresponds only to a single iteration of its QP subproblem. Assuming that they require a similar number of outer iterations, DDP remains advantageous.

\section{Numerical experiments}\label{sec:experiments}
We wish to answer the following three research questions to understand each method's relative strengths and weaknesses.
\begin{enumerate}
    \item[(R1)] How fast does each algorithm converge in terms of cost and constraints?
    \item[(R2)] How robust is each algorithm to varying initial conditions and targets?
    \item[(R3)] How effectively can we steer multiple shooting methods to avoid poor local minima via initial guess?
\end{enumerate}
To answer these three questions,
we compare the constrained DDPs and SQPs on  four different dynamical systems, that is, an inverted pendulum, a 2D quadrotor with a pendulum (quadpend) based on \cite{Singh2022}, a tadpole-like swimmer as in \cite{Tassa2007MPCDDP}, and, Franka Emika Panda robotic arm \citep{panda2022} simulated in Brax with mjx backend \cite{brax2021github, menagerie2022github}. For investigating (R1), we use all systems, while for (R2), we use the quadpend and Panda. In (R3), we focus on the quadpend. 
The state $x$ of the inverted pendulum consists of the angle and angular velocity of the pendulum as $x= [\theta, \dot{\theta}]^{\tr} \in \mathbb{R}^{2}$ and the control $u \in \mathbb{R}$ is the torque applied to the pendulum. The state of the quadpend consists of the 2D position and orientation of the quadrotor, the angle of the pendulum, and the time derivative of them, which leads to $x \in \mathbb{R}^{8}$. We control the force that the two rotors generate, and therefore $u \in\mathbb{R}^{2}$. We choose this example because of the nonlinearity of the constraints by the pendulum part. 
The swimmer is a more complex system consisting of five links, whose control is torque is applied in the four joints. The dynamics of the swimmer are found in \cite{Tassa2007MPCDDP}. The state consists of the 2D position of the nose, four joint angles and their time derivatives. Hence, the state is $x\in\mathbb{R}^{12}$. The control $u\in\mathbb{R}^{4}$ is the torques generated in the joints. The robotic arm Panda has seven joints. The state has the angles of joints and its derivatives, which gives state $x \in \Rb^{14}$. The control $u\in \Rb^{7}$ is a command of the joints. The angles and commands have limits.
Details of the dynamics of the systems and parameters are provided in Appendix \ref{sec:app_numerical}. All experiments except for Panda were performed with MATLAB \citep{MATLAB}. The Panda experiment is implemented with JAX \citep{jax2018github}.


\subsection{(R1) How fast does each algorithm converge in terms of cost and constraints?}
\label{subsec:result_one_guess}
We first explain how we evaluate the progress of optimization. Subsequently, the results are demonstrated.
\subsubsection{Progress of Optimization}
To monitor the progress of optimization, we keep track of the cost gradient over iterations. In unconstrained DDP, the gradient can be obtained by differentiating the cost-to-go at time step $k$ in \eqref{eq:DDP_cost_to_go} w.r.t. $u_{k}$ as
\begin{equation*}
     J_{u,k} = l_{u}(x_{k}, u_{k}) +  J_{u,k+1}.
\end{equation*}
The second term is computed recursively in the backward pass of DDP as
\begin{align}\label{eq:DDP_gradient_cost}
     J_{u,k+1} = J_{x,k+1}f_{u,k}, \ \ \
     J_{x,k} = l_{x, k}  + J_{x,k+1}f_{x,k}.
\end{align}
with boundary condition $J_{x,N} = \Phi_{x,N}$. We used
\begin{align*}
J_{u} = \max \{ \norm{{J_{u,1}}}_{\infty} \cdots \norm{J_{u,N-1}}_{\infty} \}
\end{align*} 
as a representative of the gradient of the cost of a trajectory. In a constrained setting, the gradient of $J$ changes to a gradient of the modified objective depending on the algorithm, as shown in the left column of Table \ref{tab:sup_method_metric}. For SQP, we use the gradient of the Lagrangian in \eqref{eq:SQP_KKT}. We note that this is intended for monitoring the optimization of one algorithm, but not for comparing the value across algorithms. We also note that the gradient of the $\log$ barrier DDP has a different property due to the barrier term. In other methods, the gradient may be used as an exit criterion, but in barrier DDP, it is not \citep{Almubarak2022}. Nevertheless, we present it for completeness.
We also monitor the values presented in the right column of the table, all of which are related to constraint violation, the residual of optimality condition for dual variables (if the method takes them into account), and penalty parameters (if the method includes any and changes them over iterations). The norm here is the infinity norm taken over all the time steps, e.g.,
\begin{align*}
\norm{[g(x_{k},u_{k})]_{+}} = \max \{\norm{g(x_{1},u_{1})]_{+}}_{\infty} \cdots \norm{[g(x_{N})]_{+}}_{\infty}\}.
\end{align*}

\begin{table*}[ht]
\footnotesize
\centering
\caption{Methods and Metrics.}
\small\sf\centering
\begin{tabular}{ llc } \hline \label{tab:sup_method_metric}
 Method & Gradient & Other metrics\\
 \hline
 $\log$ barrier DDP & $\mathcal{P}$ &  - \\
 AL DDP single & $\mathcal{L}_{\rm{A}}$ & $\norm{[g(x_{k},u_{k})]_{+}}$,  $\rho$ \\
 IP DDP single& $\mathcal{L}$ & $r_{c} = \norm{\Lambda s -\mu}$, \ $r_{s} = \norm{g(x_{k},u_{k}) + s}$,\ $\rho =  1/\mu$ \\
PDAL DDP single&  $\mathcal{L}_{\rm{PD}}$ & $\norm{[g(x_{k},u_{k})]_{+}}$, $r_{\rm{I}} = -\mu_{\rm{I}}[\pi_{\rm{I}}-\frac{\lambda}{2}]_{+} + \mu_{\rm{I}}\frac{\lambda}{2}$\\
AL DDP multi. & $\mathcal{L}_{\rm{A}}$ &  $\norm{[g(x_{k},u_{k})]_{+}}, \norm{h(x_{k},\tilde{u}_{k})}$, $\rho_{\rm{I}}, \rho_{\rm{E}}$\\
PDAL DDP multi. & $\mathcal{L}_{\rm{PD}}$ &  $\norm{[g(x_{k},\tilde{u}_{k})]_{+}}, \norm{h(x_{k},\tilde{u}_{k})}$, $r_{\rm{I}}$, $r_{\rm{E}}= \mu_{\rm{E}}[\pi_{\rm{E}}-\nu]$, $\rho_{\rm{I}}$,  $\rho_{\rm{E}}$ \\
ADMM DDP single & $\mathcal{L}_{\mathrm{A}}$&
$\norm{[g(x_{k},\tilde{u}_{k})]_{+}}, r,$ ($\rho$ is fixed.)
\\
\hline
SQP single &$\mathcal{L}$& $\norm{[g(x_{k},u_{k})]_{+}}$, $r_{c} = \norm{g(x_{k},u_{k})\odot \lambda_{k}}$\\

SQP multi. &$\mathcal{L}$& $\norm{[g(x_{k},u_{k})]_{+}}$,
$\norm{h(x_{k},u_{k})}$, $r_{c} = \norm{g(x_{k},u_{k})\odot \lambda_{k}}$\\
\hline
\end{tabular}
\end{table*}
\subsubsection{Results}
To compare the performance of each algorithm, we let all algorithms solve the same tasks, comparing cost reduction and constraint violation. The algorithms we used are as follows:
\begin{itemize}
    \item Single Shooting: $\log$ barrier, AL, IP, PDAL ADMM DDPs, and SQP.
    \item Multiple Shooting: AL with exact Hessian, AL with GN approximation \eqref{eq:AL_DDP_gauss_newton_approx}, PDAL DDPs, and SQP.
\end{itemize}
The optimization stops when either max. iteration is reached, or the regularizer in \eqref{eq:DDP_regularization} exceeds the prespecified value.

\begin{figure*}[!ht]
  \centering
  \begin{subfigure}[b]{0.45\linewidth}
 \includegraphics[trim={1cm 2cm 0.2cm 1.5cm},clip,width=\linewidth]{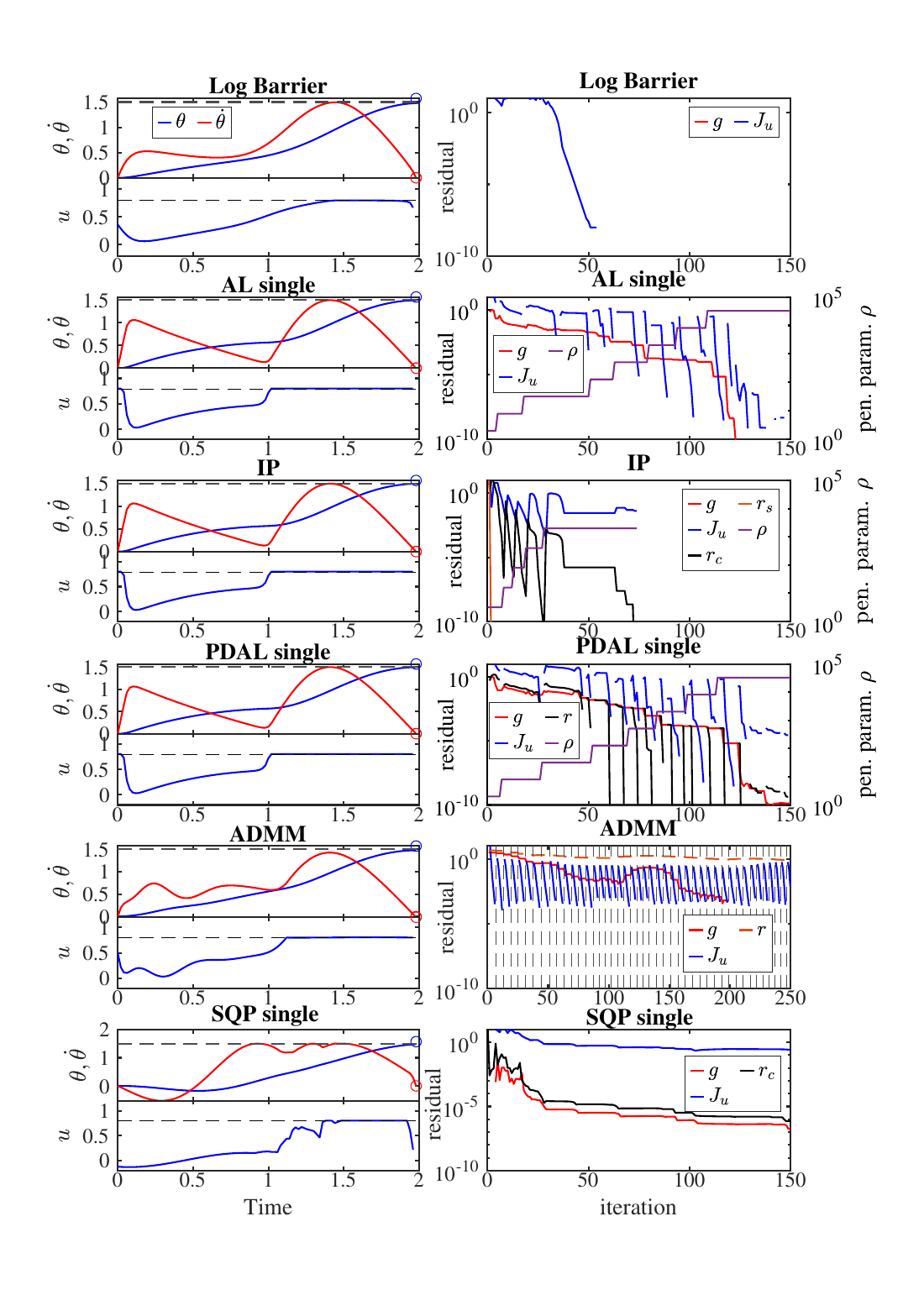}
 \end{subfigure}
 \begin{subfigure}{0.45\linewidth}
  \includegraphics[trim={1cm 2cm 0.5cm 0cm},clip,width=\linewidth]{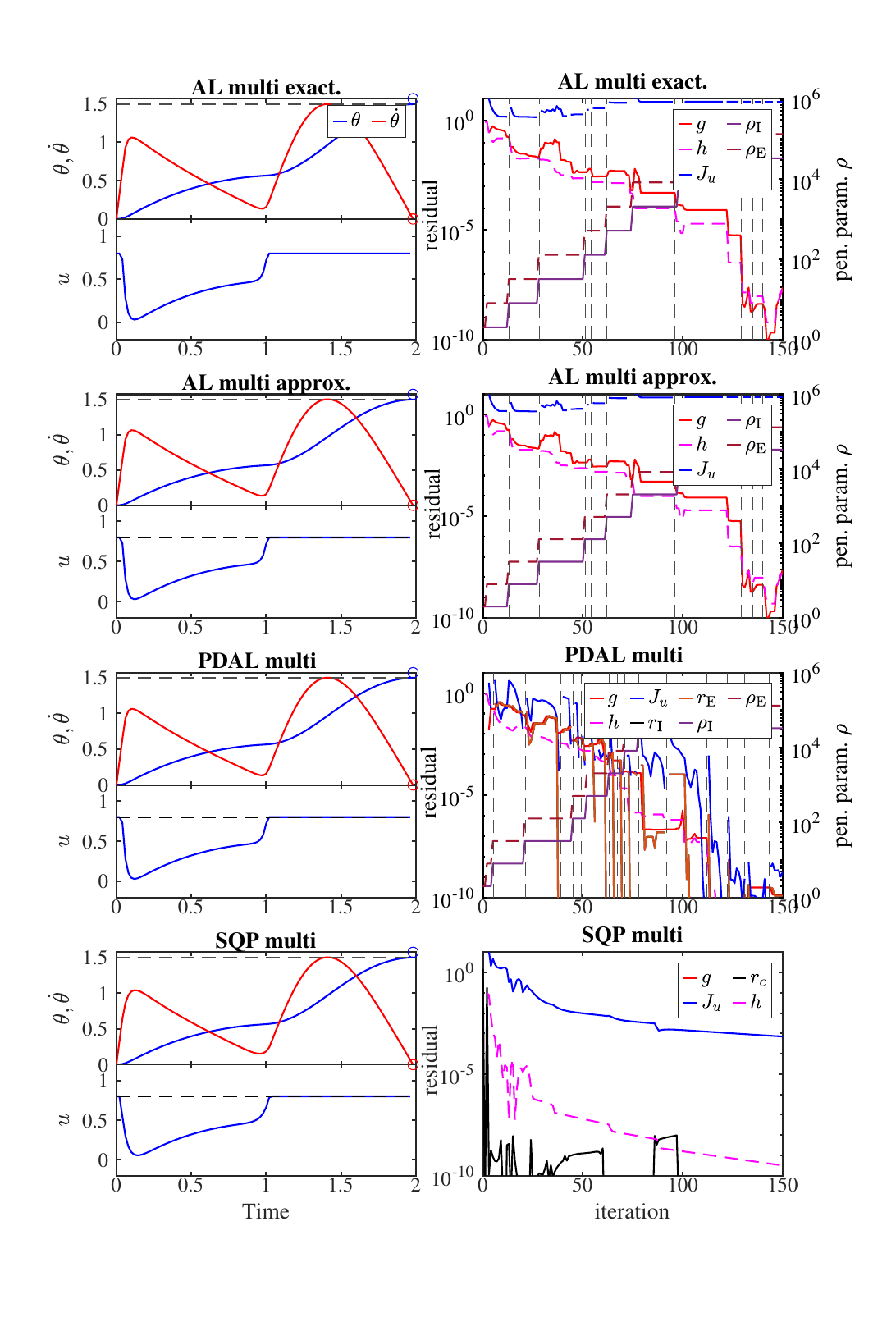}
  \end{subfigure}
  \caption{Results for inverted pendulum.}
  \label{fig:sup_pend}
\end{figure*}

\begin{figure}[tbhp]
\centering
\subfloat[Cost Reduction]
{\label{fig:inv_cost_iter}\includegraphics[trim={0.1cm 0cm 0.2cm 0.3cm},clip,width=0.5\linewidth]{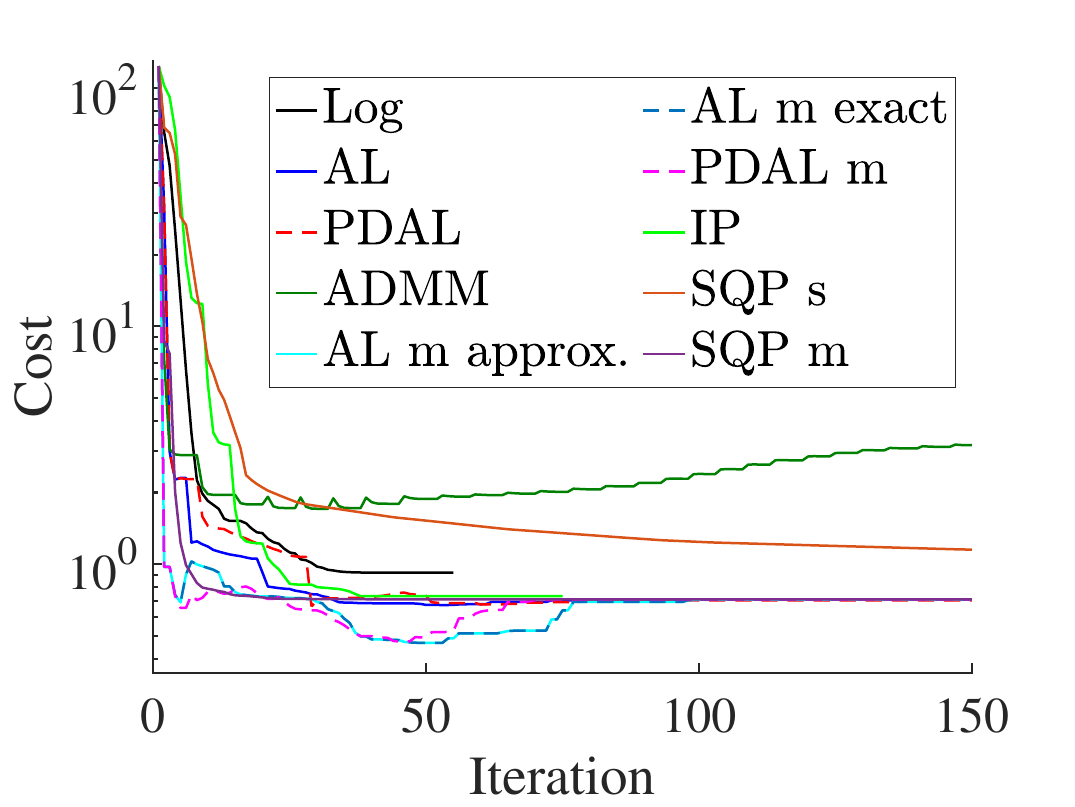}}
\subfloat[Constraint Violation]{\label{fig:inv_const_violation_iter}\includegraphics[trim={0.2cm 0cm 0.1cm 0.3cm},clip,width=0.5\linewidth]{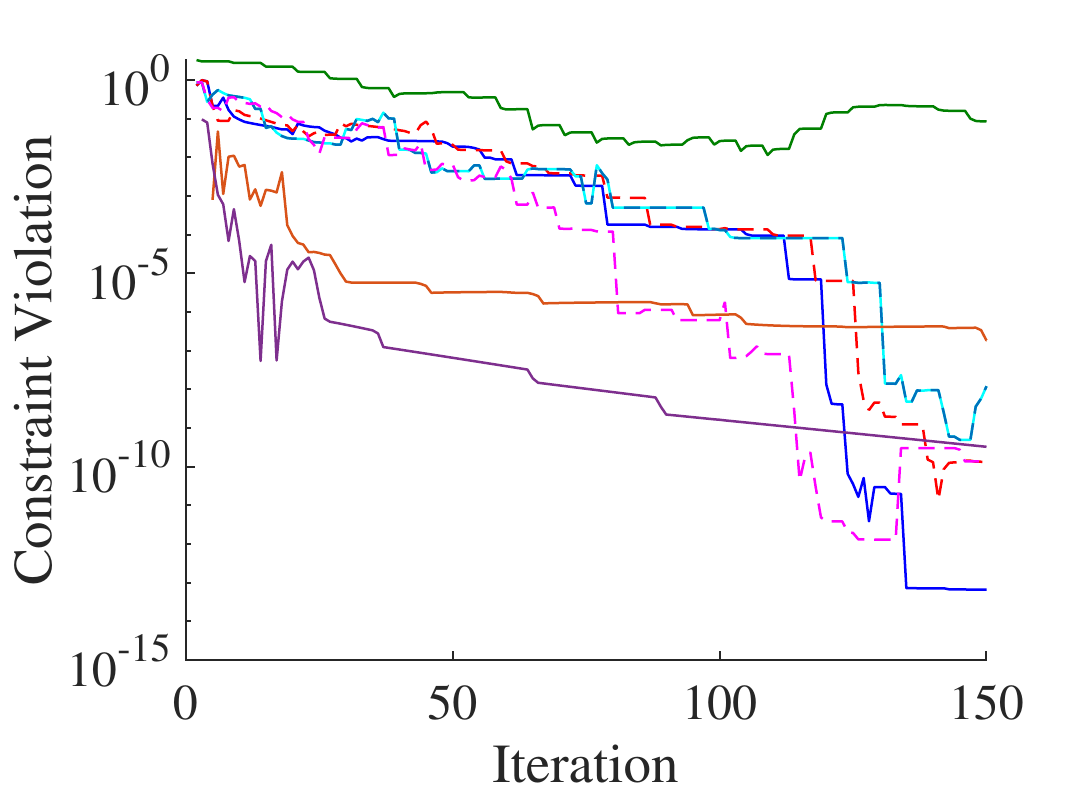}}
\caption{Comparison of algorithms with a single initial trajectory. Inverted pendulum}
\label{fig:inv_comparison_single_initial_guess}
\end{figure}

\textbf{Inverted pendulum: } In this task, the goal is to swing up the pendulum while satisfying the constraints $-0.8 \leq u_{k} \leq 0.8, -1.5 \leq \dot{\theta} \leq 1.5$.
The results are provided in Fig. \ref{fig:sup_pend}. In the right column of each figure, we present the progress of optimization with the gradient of the cost and the other metrics given in Table \ref{tab:sup_method_metric}. The vertical dotted lines in AL, PDAL, and ADMM DDP represent different outer loops. In the left column, graphs showing the evolution of $\theta$, $\dot{\theta}$ and control $u$ are illustrated. The constraints and targets are given as dotted lines and circles, respectively. Since constraints are linear in this problem, the Hessian of the constraints is zero. Thus, the exact Hessian and that with the GN approximation are identical. See AL multi exact and AL multi approx. in Fig. \eqref{fig:sup_pend}. $\log$ barrier and IP DDP terminate earlier than other methods because the descent direction cannot be found with a large regularizer in \eqref{eq:DDP_regularization}. Fig. \ref{fig:inv_cost_iter} 
shows the evolution of cost (original cost $J$ in \eqref{eq:constrained_OCP_dyn_sys}). The largest constraint violation (if any) over iterations is shown in Fig. \ref{fig:inv_const_violation_iter}. 
We observe that most methods are able to complete the task with low $(< 10^{-5})$ constraint violation.
In multiple-shooting DDPs, inequality and equality constraints exhibit comparable magnitudes, contrasting with multiple-shooting SQP, where inequality constraints consistently maintain feasibility. We postulate that this phenomenon arises from the property of SQP, that is, solving QP under constraints, rather than incorporating all elements into the $Q$ function as in DDP. Although both SQP and DDP rely on line search to find a proper step size, SQP seems to better capture the information of constraints. The $\log$ barrier method can keep the trajectory always feasible and therefore does not appear in Fig. \ref{fig:inv_const_violation_iter}. However, its cost is higher than those of the other DDP methods. This is due to the fixed penalty parameter $\mu$ in \eqref{eq:log_barr_DDP_problem}. In exchange for a simple implementation, this fixed $\mu$ makes the algorithm only approximately solve the problem, resulting in a higher cost in this experiment. In this problem, the controller needs to hit its limit over many time steps to complete the task. This makes the problem difficult for ADMM DDP, as we analyzed in section \ref{sec:AL_based_DDP_analysis}. In Appendix \ref{sec:sup_ADMM_loose_ctrl}, we further relax the control limit to $u \in[-0.9, 0.9]$ and observe that ADMM DDP is able to handle problems with less tight constraints. For SQP, the multiple shooting version performs similarly to other successful DDP methods, whereas the single shooting variant performs poorly. This difference comes from the fact that the cost function in the multiple-shooting SQP is strictly quadratic in the decision variables (state and control); on the other hand, in single-shooting, the cost is not exactly quadratic in control due to the elimination of state via linearized dynamics.

\begin{figure*}[!ht]
  \centering
  \begin{subfigure}[b]{0.45\linewidth}
 \includegraphics[trim={0cm 0cm 0cm 0cm},clip,width=\linewidth]{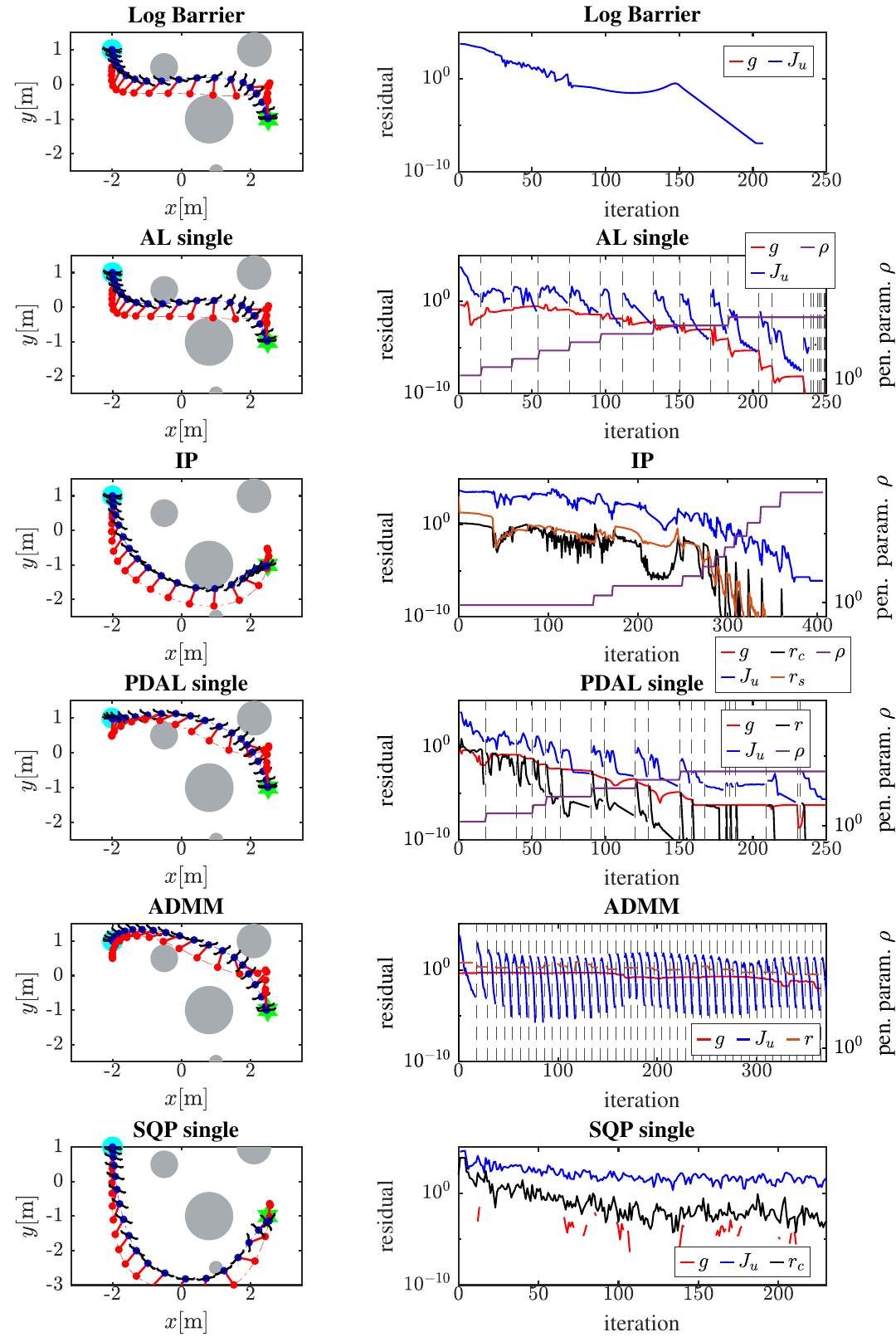}
 \end{subfigure}
 \begin{subfigure}{0.45\linewidth}
  \includegraphics[trim={0cm 0cm 0cm 0cm},clip,width=\linewidth]{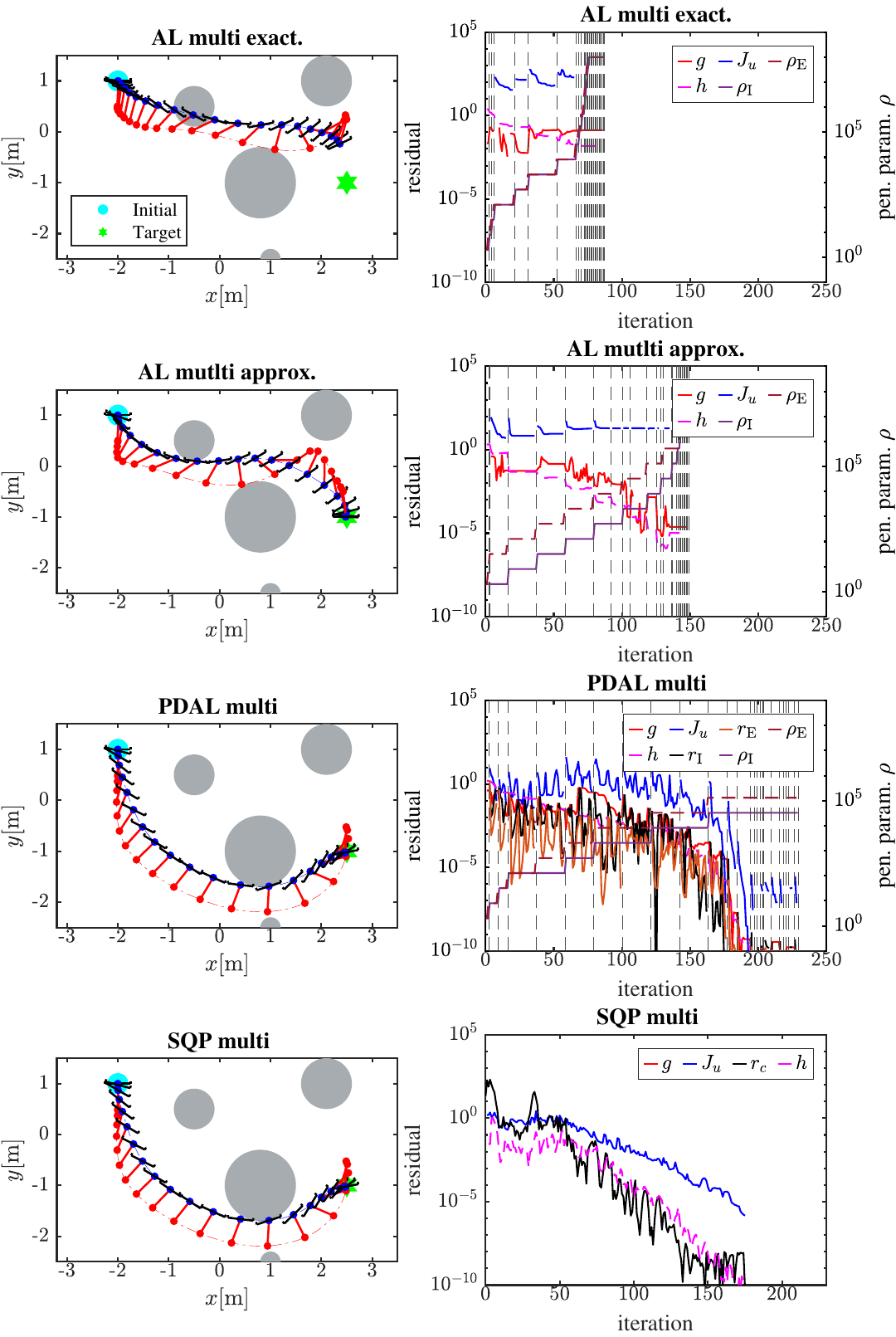}
  \end{subfigure}
  \caption{Results for quadpend.}
  \label{fig:sup_quad}
\end{figure*}

\begin{figure}[tbhp]
\centering
\subfloat[Initial trajectory. \label{fig:quad_init_last}]
{\includegraphics[trim={0.1cm 0cm 0.1cm 0.1cm},clip,width=0.5\linewidth]{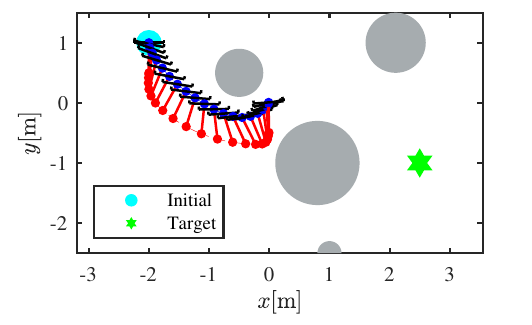}}

\subfloat[Cost Reduction]{\label{fig:quad_cost_iter}\includegraphics[trim={0.1cm 0cm 0.1cm 0.1cm},clip,width=0.5\linewidth]{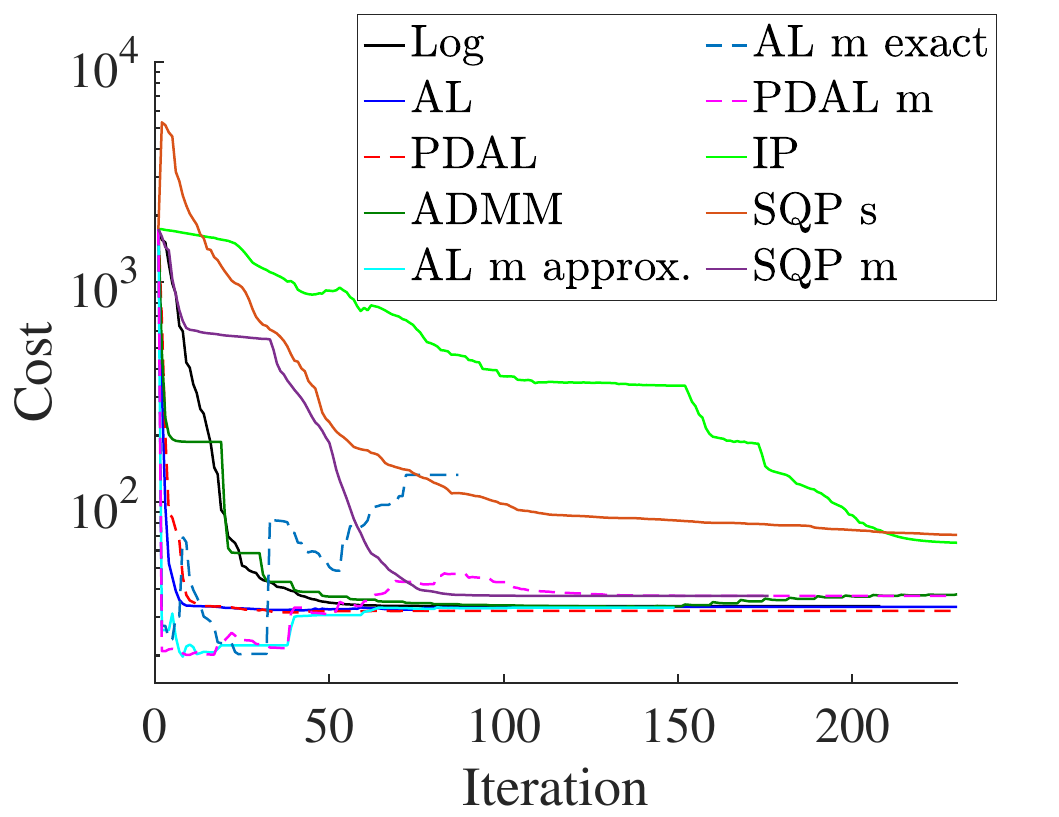}}
\subfloat[Constraint Violation]{\label{fig:quad_const_violation_iter}\includegraphics[width=0.5\linewidth]{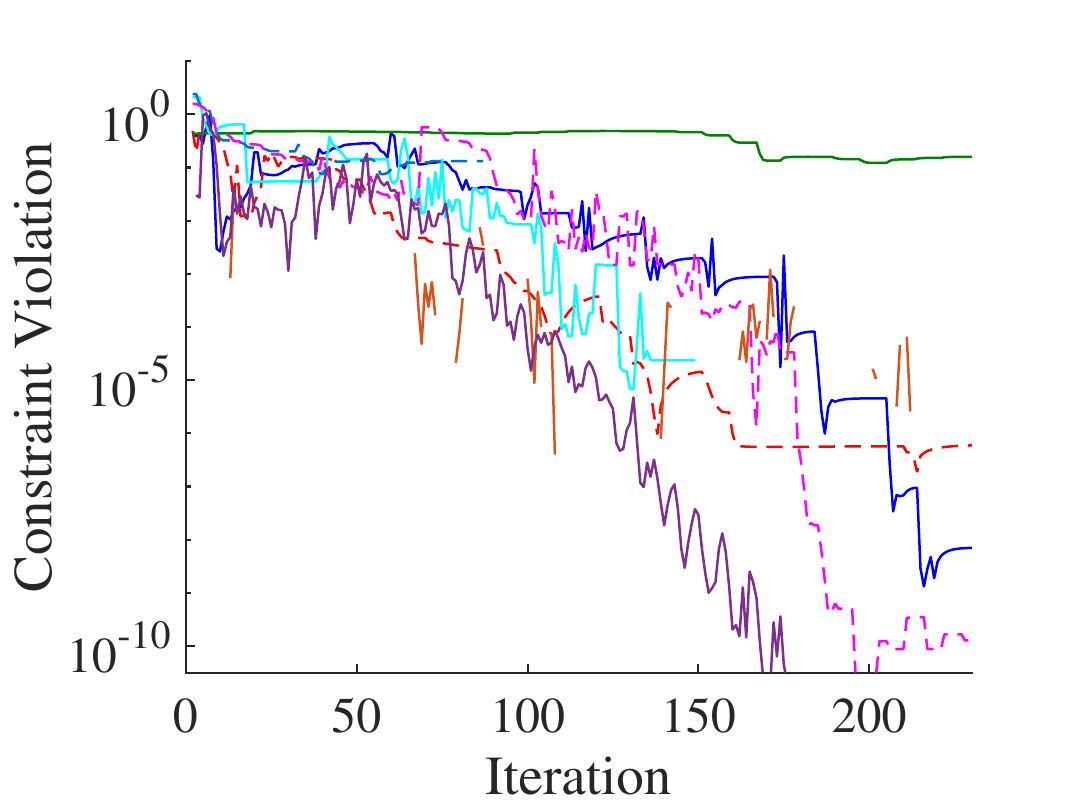}}
\caption{Comparison of algorithms with a single initial trajectory. Quadpend.}
\label{fig:quad_comparison_single_initial_guess}
\end{figure}

\textbf{Quadpend:} The quadpend is navigated to reach the target $x_{\rm{g}} = [2.5, -1, 0, {\pi}/{2}, 0, 0, 0, 0]^{\tr}$ (pendulum upright)
from the initial state $x_{0} = [-2, 1, 0, 0, 0,0,0,0]^{\tr}$ (pendulum down) while avoiding four obstacles under box control constraints. The initial trajectory is shown in Fig. \ref{fig:quad_init_last}, which is obtained with single-shooting PDAL DDP. 
We show the results of single and multiple shooting methods in Fig. \ref{fig:sup_quad}. 
In this experiment, multiple-shooting AL DDPs reached maximum iteration by reaching a plateau where the inner loops could not find a descent direction. The multiple-shooting SQP algorithm terminated due to failure in the line search to find a descent direction even after regularizing the Hessian as explained in Section \ref{sec:SQP_single}. 
Here, the GN approximation makes a difference because of nonlinear constraints. Indeed, we observe the effectiveness of the approximation as multiple-shooting AL DDP fails without the approximation, getting stuck at a poor local minimum with a large constraint violation. Other than the multiple-shooting AL DDP with exact Hessian, all algorithms can let the quadpend hit the target with the pendulum up, although the solutions vary. This is because the algorithms are local methods that use the local approximation of the cost and dynamics.

Fig. \ref{fig:quad_cost_iter} and Fig. \ref{fig:quad_const_violation_iter} compare the cost and constraint violation. 
Again, we note that in Fig. \ref{fig:quad_const_violation_iter}, some algorithms do not appear when their trajectories are feasible. AL-based DDPs (AL, PDAL, and ADMM) show rapid cost reduction with constraint violation. They hit the target first and then gradually satisfy the constraints, where ADMM shows the slowest improvement in constraint satisfaction. The other two DDP methods, IP and barrier DDP, can keep the trajectory feasible during optimization. IP DDP's slowness comes from the nature of the IP method, where the trajectory is biased to follow the central path as in \eqref{eq:IP_KKT}. In this experiment, $\log$ barrier DDP is a good option that balances cost reduction and constraint satisfaction. For SQP, as in the previous case, the multiple-shooting DDP is comparable to other DDP methods or even better in constraint satisfaction, but the shingle-shooting method performs poorly.

\begin{figure*}[!ht]
  \centering
  \begin{subfigure}[b]{0.45\linewidth}
 \includegraphics[trim={0cm 0cm 0cm 0.1cm},clip,width=\linewidth]{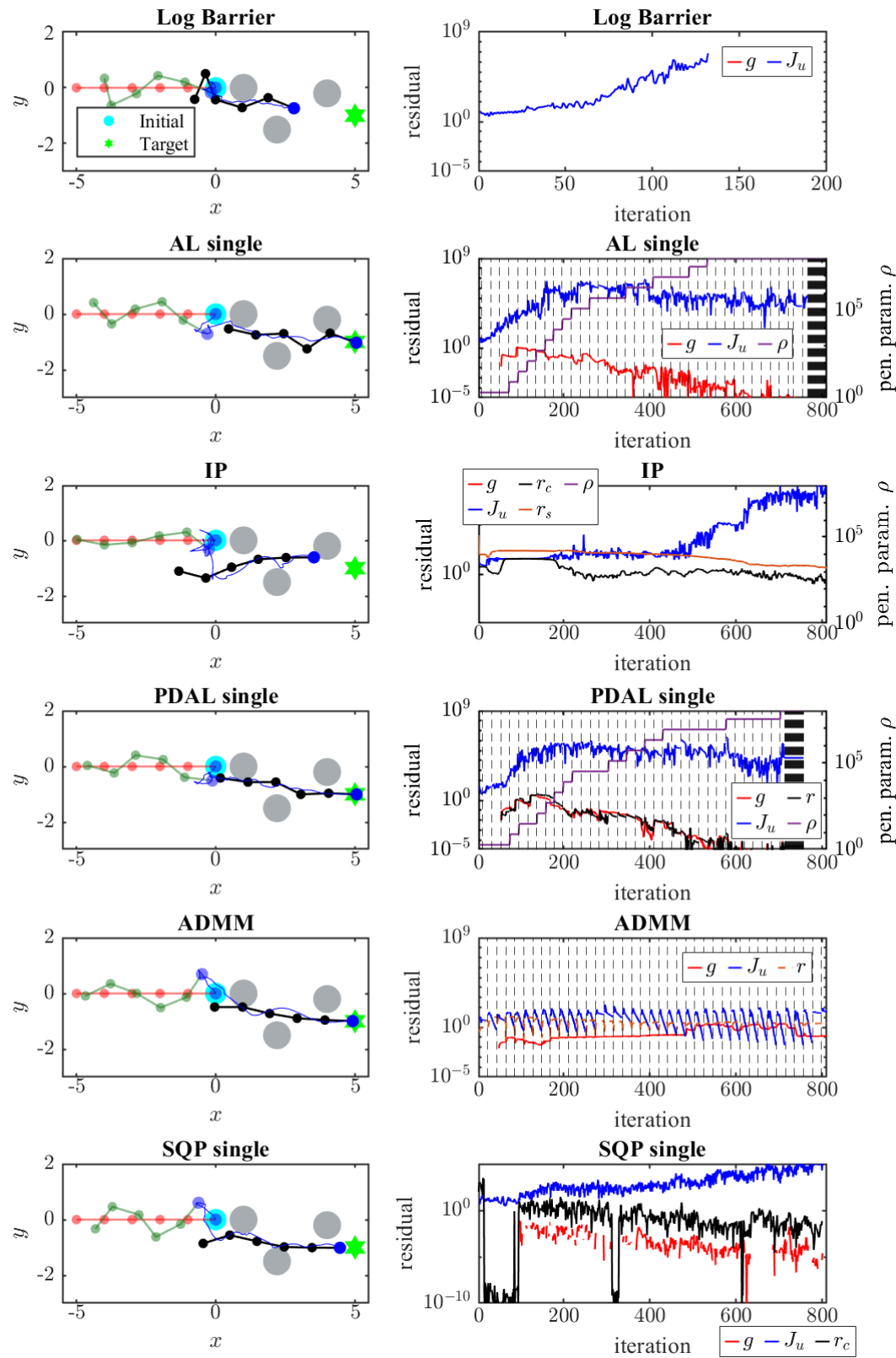}
  \end{subfigure}
 \begin{subfigure}{0.45\linewidth}
  \includegraphics[trim={0.1cm 0cm 0.5cm 0cm},clip,width=\linewidth]{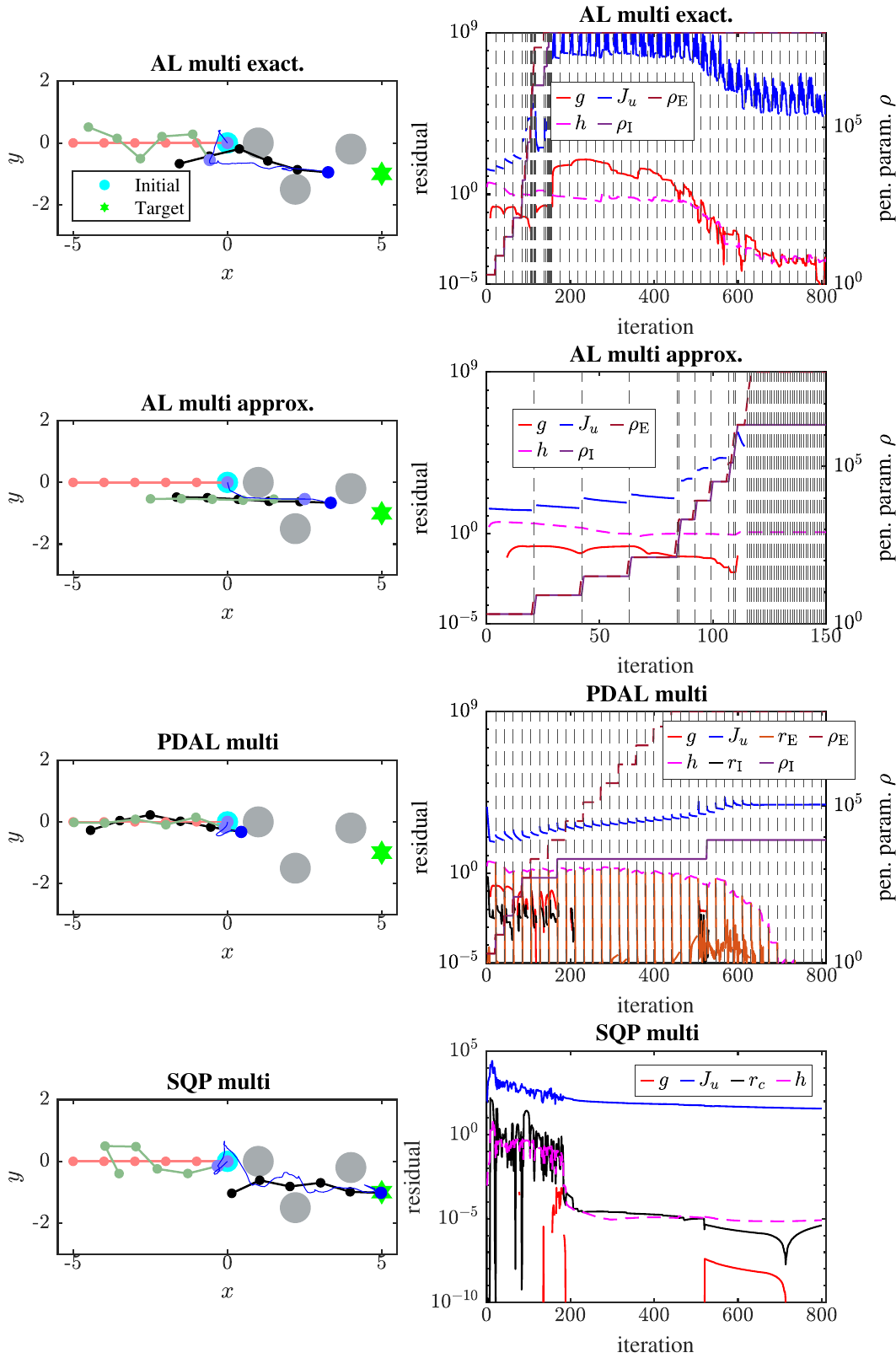}
  \end{subfigure}
  \caption{Results for swimmer. The initial, middle, and terminal configurations of the swimmer are drawn in red, green, and black.}
  \label{fig:sup_swimmer}
\end{figure*}

\begin{figure}[h!]
\centering
\subfloat[Cost Reduction]
{\label{fig:swimmer_cost_iter}\includegraphics[trim={0.2cm 0cm 0.2cm 0.3cm},clip,width=0.5\linewidth]{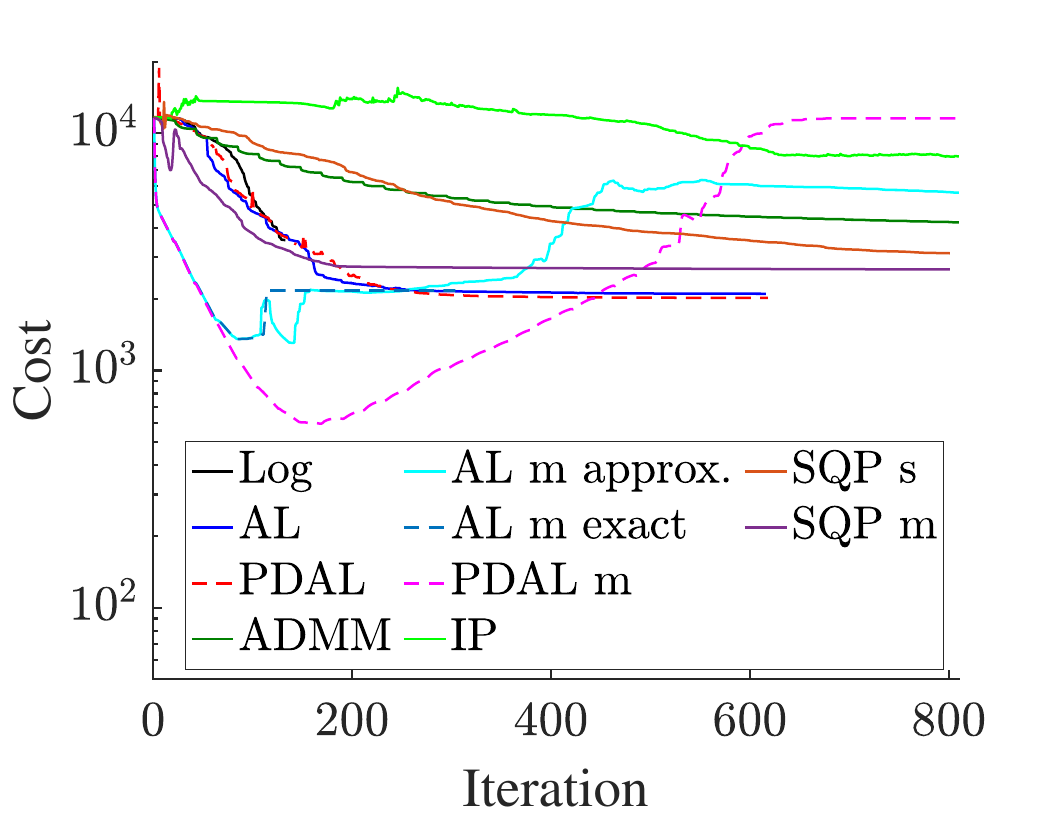}}
\subfloat[Constraint Violation]{\label{fig:swimmer_const_iter}\includegraphics[trim={0.2cm 0cm 0.1cm 0.3cm},clip,width=0.5\linewidth]{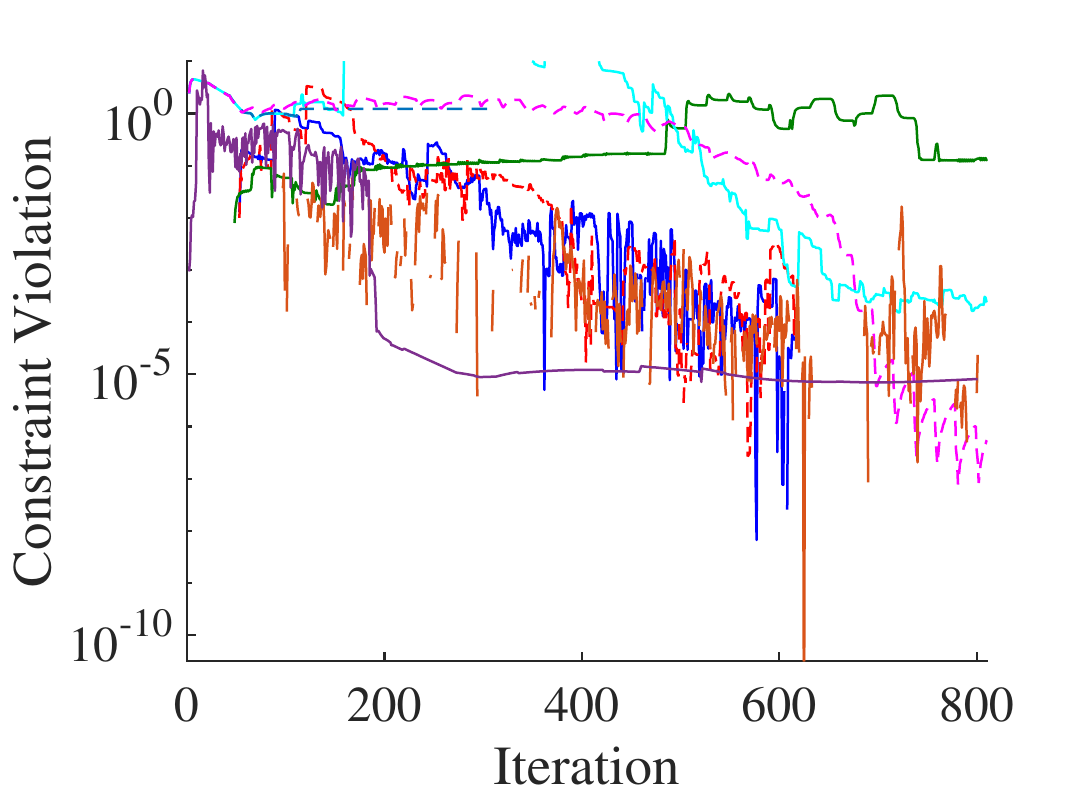}}
\caption{Comparison of algorithms with a single initial trajectory. Swimmer.}
\label{fig:swimmer_comparison_single_initial_guess}
\end{figure}

\textbf{Swimmer:} Starting from the initial position, with its nose in the origin and its tail straight, the swimmer tries to hit the target $[5, -1]$ with its nose. The initial control sequence is zeros, which keeps the swimmer in the initial place. The results are shown in Fig. \ref{fig:sup_swimmer} and Fig. \ref{fig:swimmer_comparison_single_initial_guess}, where the difference in the algorithm's performance becomes more evident than in the previous examples. The $\log$ barrier DDP method gets stuck before hitting the target due to the approximation. IP DDP shows slow progress and cannot hit the target. Although single-shooting AL and PDAL DDP can complete the task in a similar order of constraint violation $(10^{-5})$, PDAL achieves a lower cost.

In this example, multiple-shooting DDPs struggled to solve the problem. Specifically, they can hit the target with a significant equality constraint violation (especially dynamics) early in the optimization process. However, to satisfy the constraints, they start to show conservative motions and finally end up staying at the initial position with the initial state at a high cost. This is because the equality constraints are always satisfied if they do not move. To alleviate the issue, we regularize the state part of $Q_{\tilde{u} \tilde{u}}$ similar to \eqref{eq:DDP_regularization} and the technique in \citep{Jallet2022implicit}, which prevent $x$s from moving too far from the current trajectory.
The results presented here are obtained with this regularization strategy. This modification works in AL DDP with approximated Hessian, making the swimmer move forward. Multiple-shooting SQP can achieve the task with the lowest constraint violation among all methods that can complete the task. Overall, single-shooting AL, PDAL DDPs, and multiple-shooting SQP can complete the task with a similar order of constraint violation. PDAL single-shooting DDP achieves the best cost, followed by AL single-shooting DDP, and multiple-shooting SQP.

\begin{figure*}[!ht]
  \centering
  \begin{subfigure}[b]{0.45\linewidth}
 \includegraphics[trim={0cm 0cm 0cm 0.1cm},clip,width=\linewidth]{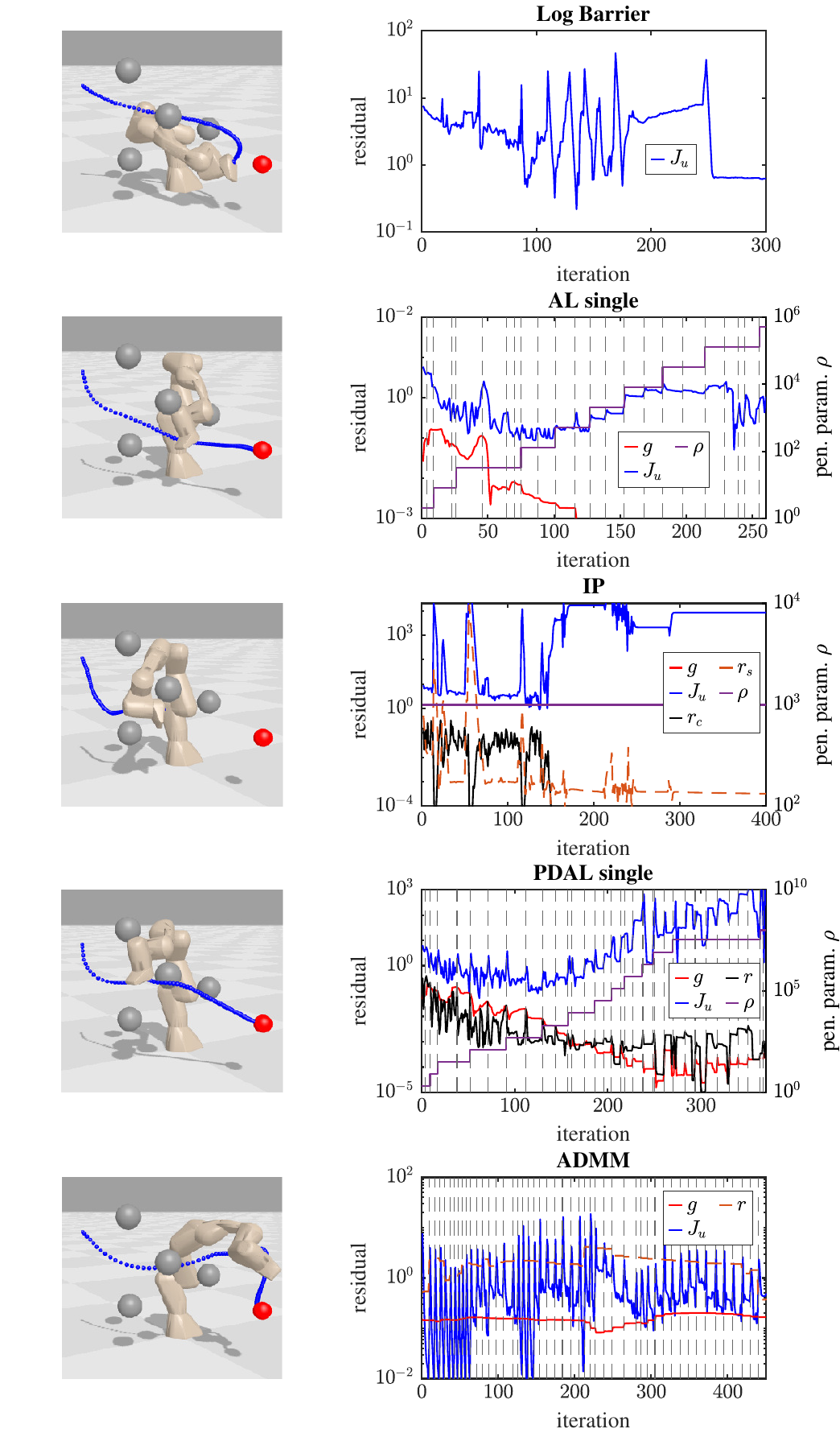}
 \end{subfigure}
 \begin{subfigure}{0.45\linewidth}
  \includegraphics[trim={0.1cm 0cm 0.2cm 0cm},clip,width=\linewidth]{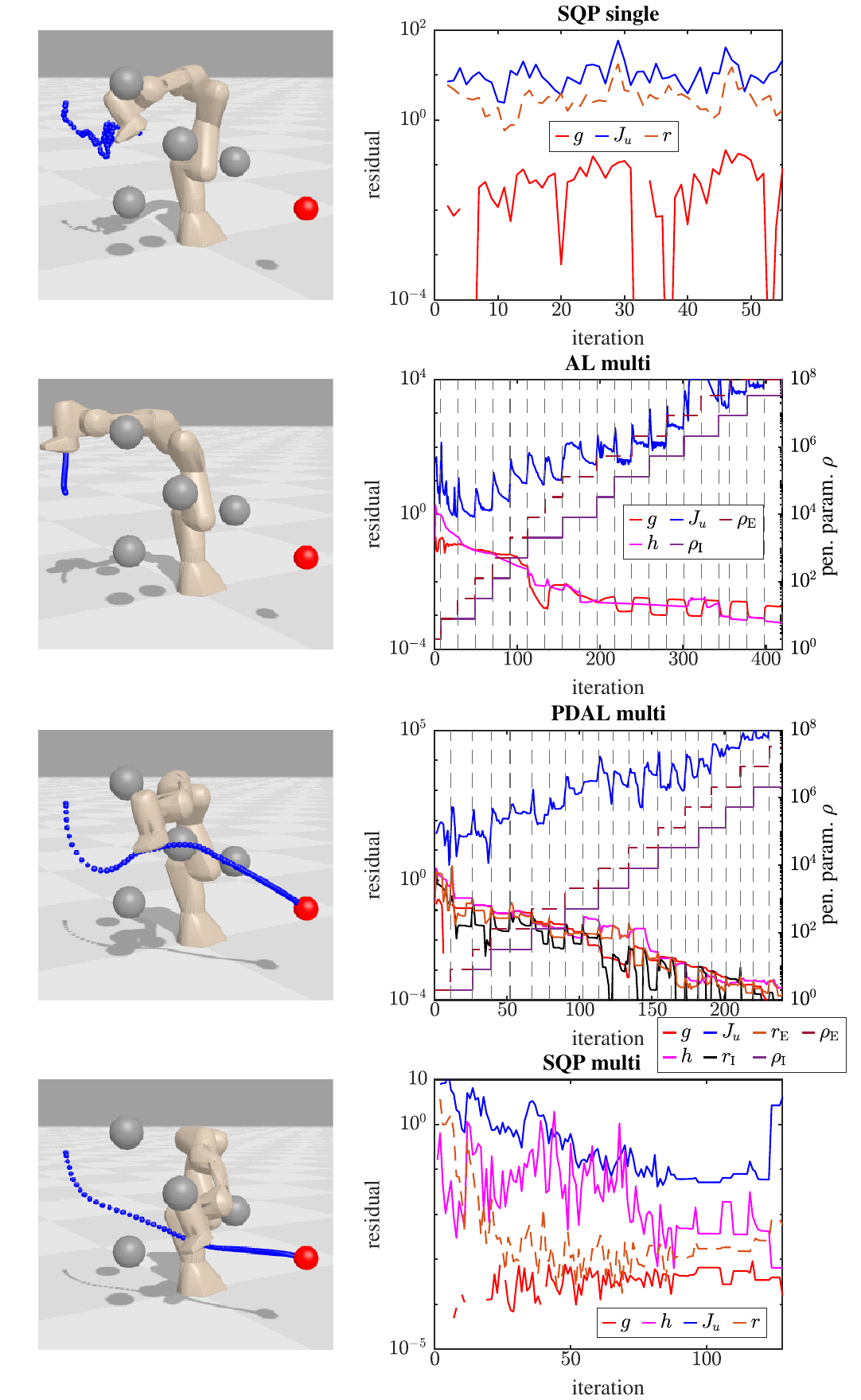}
  \end{subfigure}
  \caption{Results for panda. The target and obstacles are shown in red and gray spheres, respectively. The trajectory of the end-effector is drawn in blue spheres.}
  \label{fig:sup_panda}
\end{figure*}
\begin{figure}[!ht]
\centering
\begin{subfigure}[b]{0.48\linewidth}
\includegraphics[trim={0.3cm 0cm 0.2cm 0.3cm},clip,width=\linewidth]{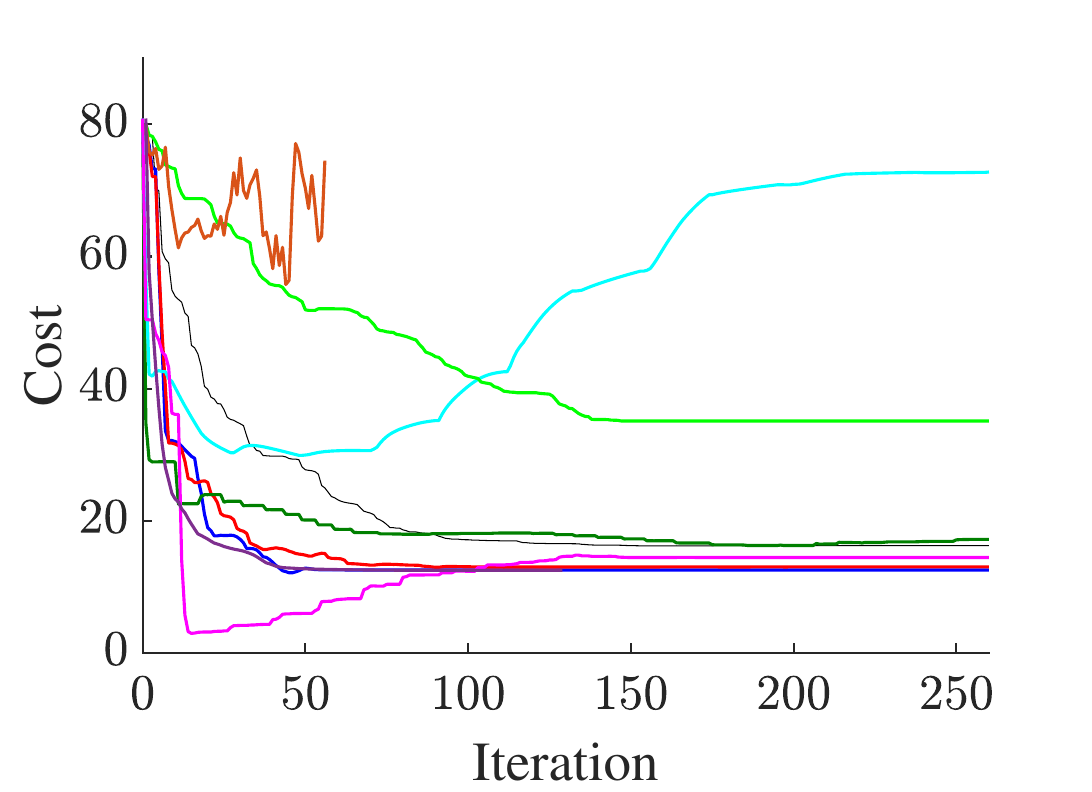}
\caption{Cost Reduction}
\label{fig:panda_cost_iter}
\end{subfigure}
\begin{subfigure}[b]{0.48\linewidth}
\includegraphics[trim={0.2cm 0cm 0.1cm 0.3cm},clip,width=\linewidth]{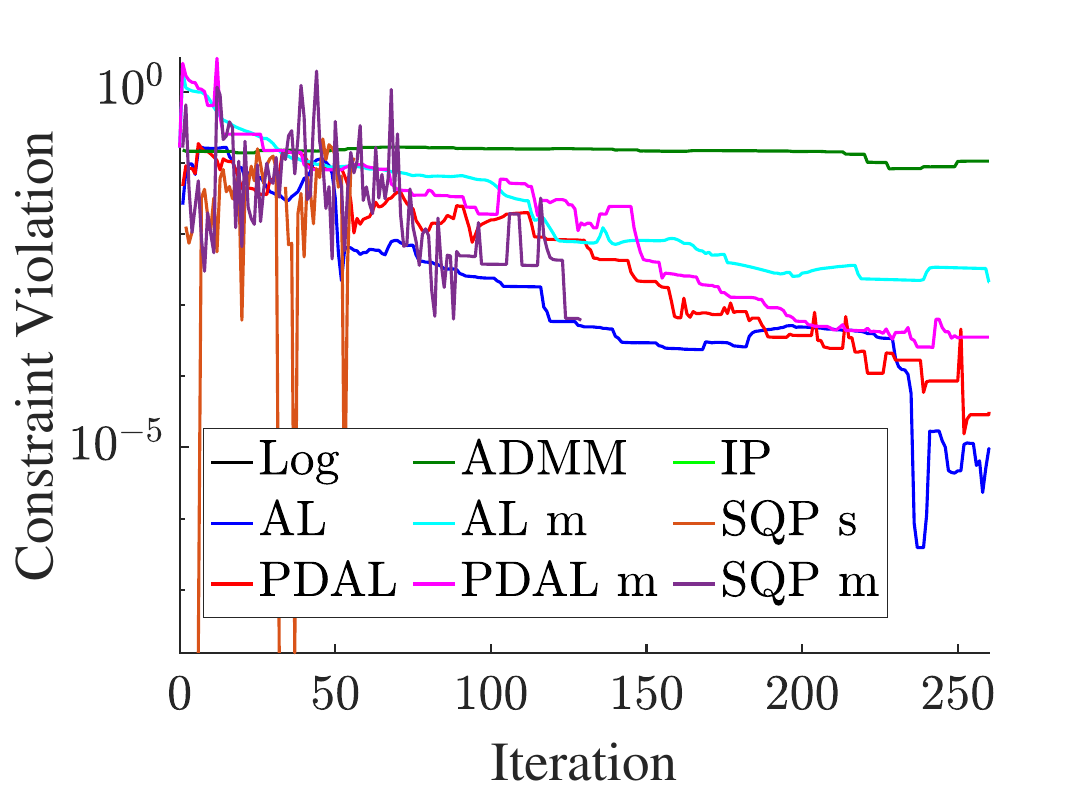}
\caption{Constraint Violation}
\label{fig:panda_const_iter}
\end{subfigure}
\caption{Comparison of algorithms with a single initial trajectory. Panda.}
\label{fig:panda_comparison_single_initial_guess}
\end{figure}

\textbf{Panda:} The task of the arm is to place the end effector on the target position and stop, while avoiding four spherical obstacles, and satisfying the joint and its command limit.
In this example, the second-order information of dynamics is not stably available from the simulator. Therefore, we drop the corresponding terms of constraints in AL and PDAL multiple-shooting DDPs. Consequently, the AL multiple-shooting method only has an approximate version. The initial position of the arm is given by the joint angles from the base link as $x_{0} = [0, 0, 0, -\pi/2, 0, \pi/2, -\pi/4]$. The initial control command maintains the initial state. The results are presented in Fig. \ref{fig:sup_panda} and Fig \ref{fig:panda_comparison_single_initial_guess}. 
The algorithms find different local solutions, as in the quadpend. Similarly to the swimmer example, $\log$-barrier and IP DDP get stuck and cannot hit the target, while other single-shooting DDPs can complete the task. An addition of the feasibility restoration mechanism could improve the performance of IP DDP. AL and PDAL show similar performance, while ADMM has high constraint violation. This is because we use a small penalty parameter to prioritize hitting the target rather than satisfying constraints in ADMM DDP. For multiple-shooting DDPs, although both AL and PDAL can reduce constraint violations over iterations, only PDAL can achieve the task. SQPs show a tendency similar to that in other experiments. The multiple-shooting one can complete the task with a similar magnitude of constraint violation as other successful DDPs, while the single-shooting one performs poorly.


\subsection{(R2) How robust is each algorithm to varying initial conditions and targets?}
\label{subsec:result_several_points}
\begin{figure}[htbp]
\centering
\includegraphics[trim={0cm 1cm 0cm 0.5cm},clip,width=\linewidth]{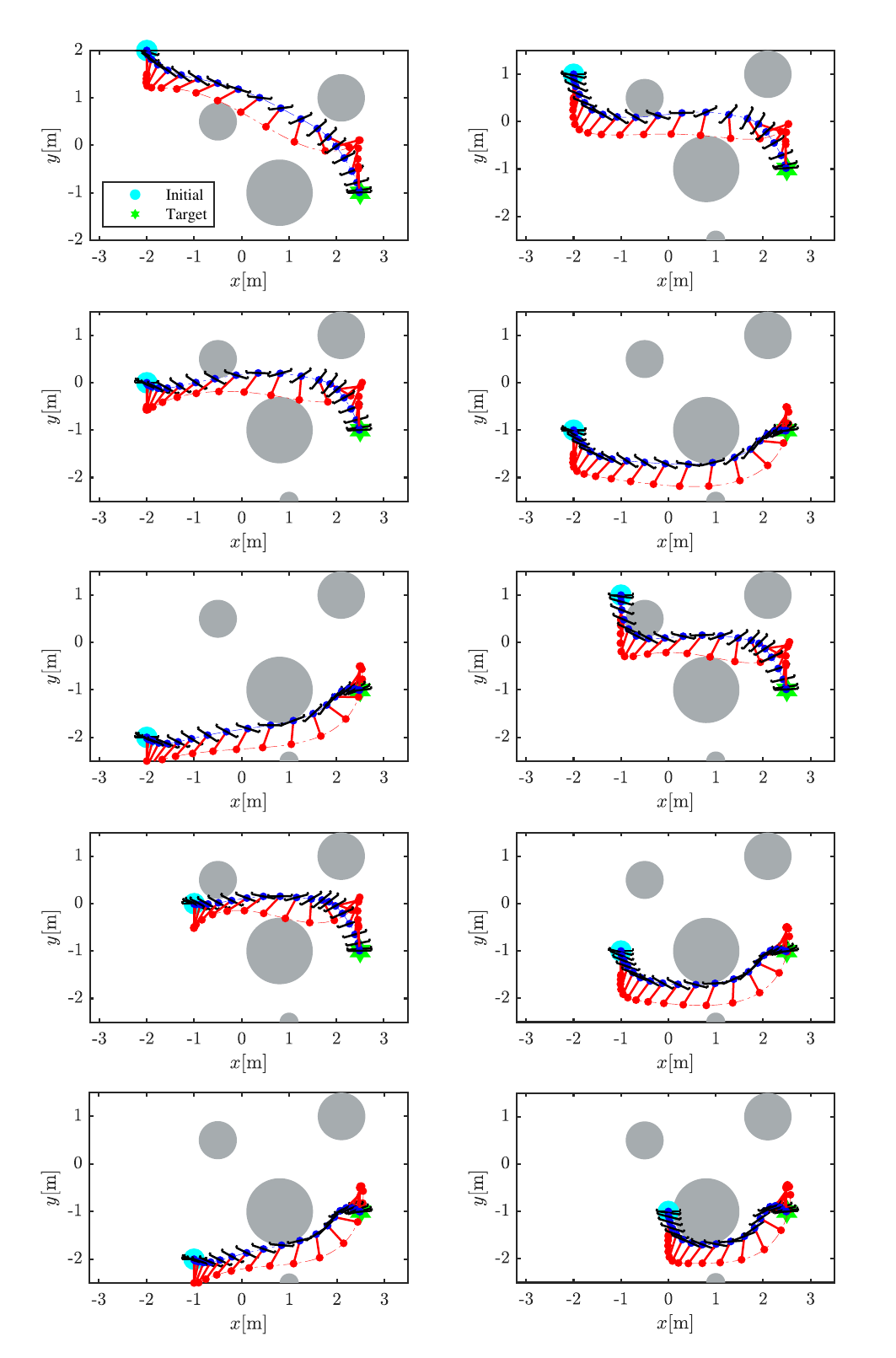}
\caption{Results of quadpend starting from multiple initial points.}
\label{fig:sup_quad_multiple_start}
\end{figure}

\textbf{Initial conditions:} To examine robustness while varying the initial conditions, we use the quadpend system and initialize the algorithms with ten different initial hovering trajectories with the pendulum down. The hovering trajectory is achieved by an initial control sequence $u_{1:N-1} = 0.5(m_{\rm{q}}+m_{\rm{p}})g_{0}[1,1]^{\tr}$, where $m_{\rm{q}}$, $m_{\rm{p}}$, and $g_{0}$ are the mass of the quadrotor part, the mass of the pendulum part, and the gravitational acceleration, respectively. 
The algorithms solve the same task as in Section  \ref{subsec:result_one_guess}. Tables \ref{tab:single_DDP} and \ref{tab:multi_DDP} present the results of the single- and multiple-shooting algorithms, respectively. A successful run is defined as the quadpend hitting the target. With $m$ and $\sigma$, we denote the mean and standard deviation of the values specified by their subscripts $J$, $I$, and $E$, which denote the cost, inequality, and equality constraint violation. When all trajectories are feasible, the symbol $\checkmark$ is used. The arrows in the table represent the preferred value. For example, we have $\uparrow$ next to the success rate because a robust algorithm can hit many targets. To provide an overview of the experiment, we show the resulting trajectories in Fig. \ref{fig:sup_quad_multiple_start} obtained by multiple-shooting SQP. 

Among single-shooting methods, PDAL DDP achieves the lowest mean and standard deviation in cost. Although its constraint violation is not as small as that of AL DDP, it is sufficiently small. The $\log$ barrier and IP DDP methods achieve strict feasibility but have higher costs than AL and PDAL. ADMM has a large constraint violation in the mean cost caused by two infeasible trajectories.
Overall, all single-shooting DDP methods are robust to the initial condition. Among the multiple-shooting methods, SQP performs best, achieving the best value on almost all items. 

\begin{figure}[tbhp]
\centering
\subfloat[Typical failure mode with hovering initial condition. \label{fig:quad_fail_through}]
{\includegraphics[trim={0.1cm 0cm 0.1cm 0.1cm},clip,width=\linewidth]{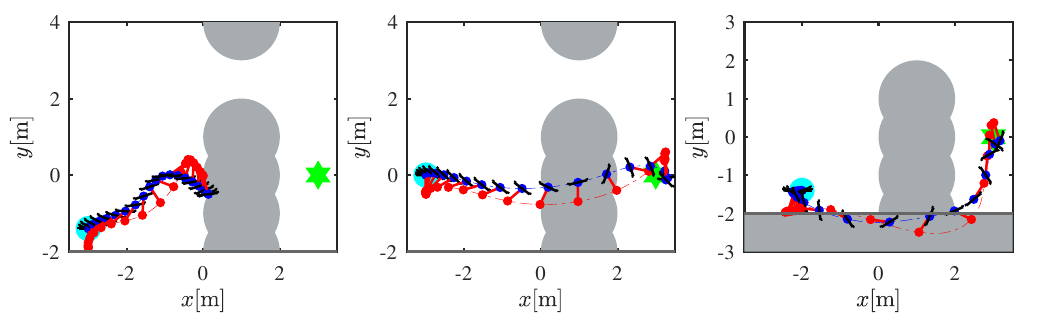}}

\subfloat[Initial and Optimal Trajectories with Multiple shooting method.]
{\label{fig:quad_warm_init}\includegraphics[trim={0.1cm 0cm 0.1cm 0.1cm},clip,width=0.7\linewidth]{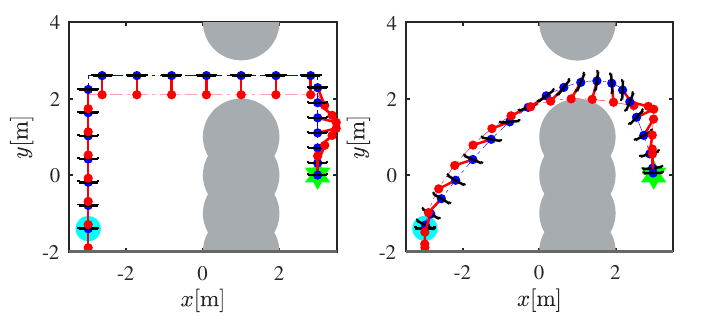}}
\caption{Failure, initial and optimal trajectories of experiments in \ref{subsec:result_warm_multi}.}
\label{fig:quad_warm_fail_init_last}
\end{figure}

\begin{figure}[htbp]
    \centering
    \includegraphics[width=0.8\linewidth]{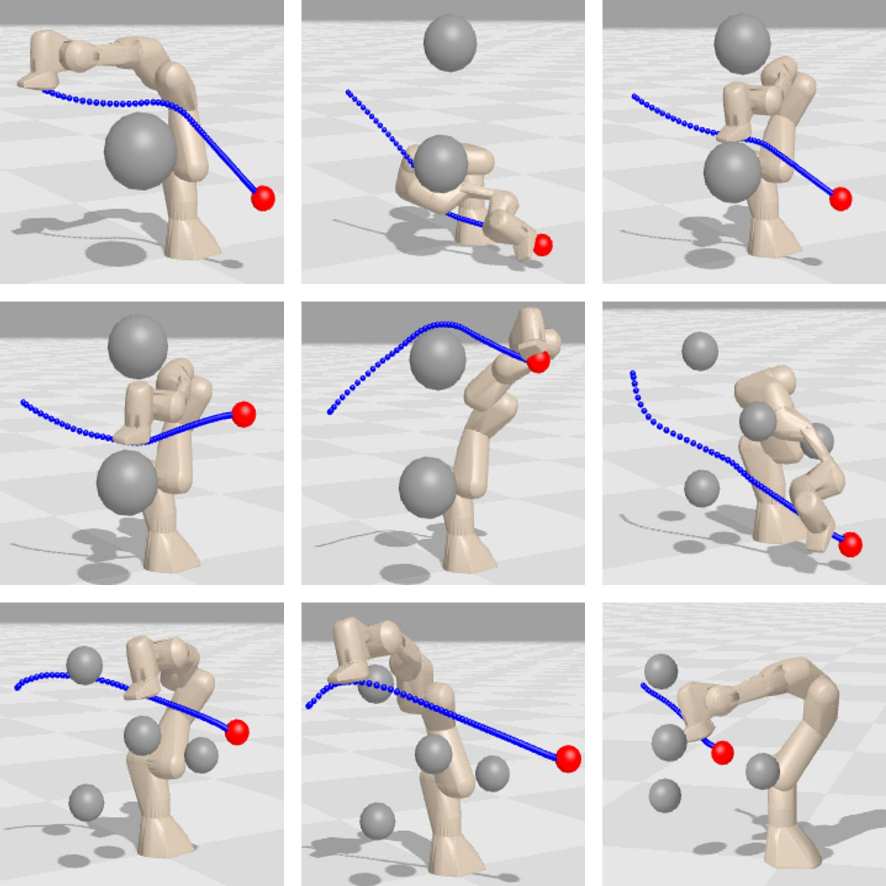}
    \caption{Results of Panda for multiple targets.}
    \label{fig:panda_trajectories}
\end{figure}

\begin{figure}[!ht]
\centering
\begin{subfigure}[b]{\linewidth}
\includegraphics[trim={0.2cm 0cm 0.2cm 0.3cm},clip,width=\linewidth]{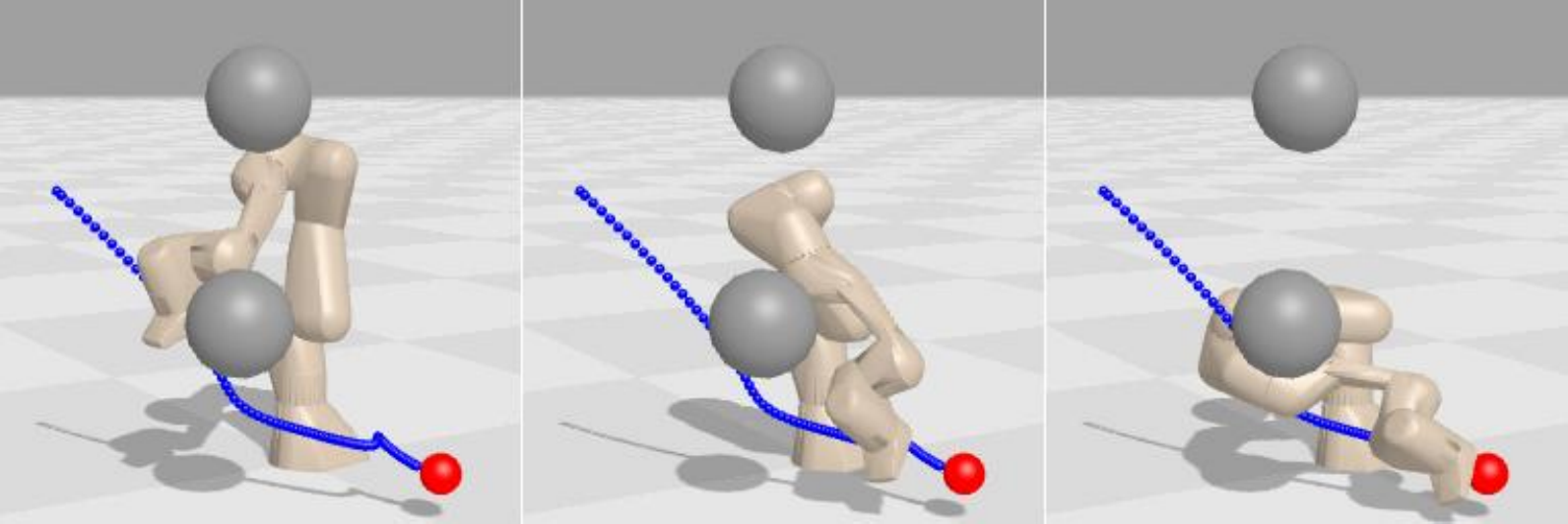}
\caption{Single-shooting PDAL DDP}
\label{fig:panda_pdal}
\end{subfigure}
\begin{subfigure}[b]{\linewidth}
\includegraphics[trim={0.2cm 0cm 0.1cm 0.3cm},clip,width=\linewidth]{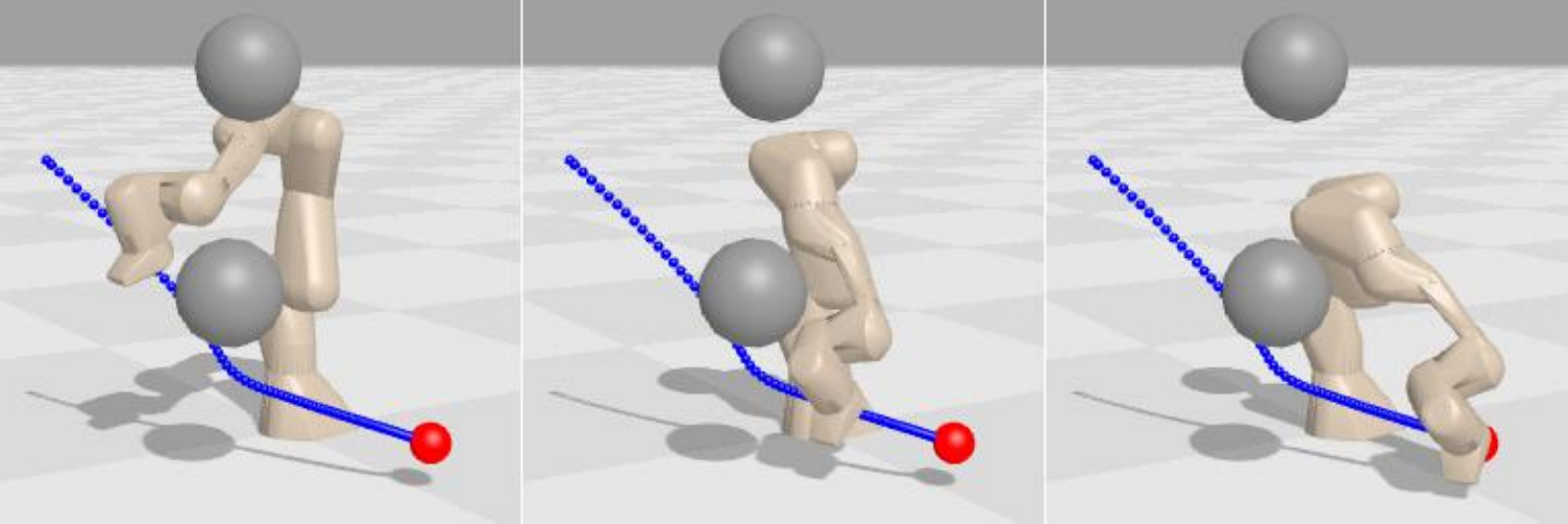}
\caption{Multiple-shooting SQP}
\label{fig:panda_sqp}
\end{subfigure}
\caption{Comparison of a failed trajectory of single-shooting PDAL DDP and a successful trajectory from multiple-shooting SQP.}
\label{fig:panda_pdal_sqp}
\end{figure}

\noindent\textbf{Different targets:}
To evaluate the robustness of the algorithms across different targets and tasks, we design three obstacle fields, containing one, two, and four obstacles, respectively. For each field, we conduct one, four, and five experiments with varying target configurations. In total, we have ten experiments. We note that one of these environments, i.e, a pair of an obstacle field and a target, is the same as the one used in the previous example. The results are provided in Table \ref{tab:single_DDP_panda} and Table \ref{tab:multi_DDP_panda}. The overview of the experiments is shown in Fig. \ref{fig:panda_trajectories} and Extension 1. 

In single-shooting methods, $\log$ barrier and IP DDPs show conservative behavior to take distance from constraints through the barrier function. As a result, their performance decreases in cluttered environments. Among the single-shooting algorithms, PDAL DDP has the highest success rate with sufficiently small constraint violation. The multiple-shooting SQP outperforms it in terms of success rate with equality constraint violation. In Fig. \ref{fig:panda_pdal_sqp},  we compare a failed trajectory from single-shooting PDAL DDP and a successful trajectory from multiple-shooting SQP in the same task. In the early time steps, both of them draw similar trajectories. However, only SQP can successfully hit the target by rotating a joint close to the base.
In contrast, a local minimum captures PDAL DDP, where it continues to bend the arm without incorporating the rotation observed in SQP. Panda can locate the end effector close to the target in this local solution, but cannot hit it. Although single-shooting DDP failed in this example, it can achieve a lower mean cost compared to multiple-shooting SQP. 

We conclude from these results that when constraint violation is the most critical factor, $\log$ barrier DDP is the best method. However, it might not be able to complete the task, especially in high-dimensional systems and cluttered environments.
If users can accept small constraint violations, single-shooting PDAL DDP or multiple-shooting SQP has a higher success rate. In our example, SQP shows a slightly better success rate, whereas PDAL DDP achieves a lower mean cost. A key difference between these two algorithms is that PDAL DDP can always satisfy dynamics, whereas SQP has small violations.

\subsection{(R3) How effectively can we steer multiple shooting methods to avoid bad local minima via the initial guess?}\label{subsec:result_warm_multi}
One key advantage of multiple-shooting methods over single-shooting ones is that they can enjoy good initial guesses, which was not demonstrated in previous experiments. In the multiple-shooting formulation, an arbitrary state sequence can be used as an initial state trajectory as explained in Sections \ref{subsubsec:DDP multi} and \ref{subsec:SQP_multi}, which can help guide the algorithm away from poor local minima. To showcase this ability, we test multiple-shooting algorithms with a new task where the quadpend flies through a narrow gap of obstacles and reaches a target behind them with ten different initial points. When initialized with hovering states and controls at a single starting point, as in the previous experiment, the performance of all methods is decremented. They cannot complete the task, i.e., hitting the target with a significant constraint violation, and getting stuck before hitting it, etc. We show typical failure trajectories in Fig. \ref{fig:quad_fail_through}.
\begin{table*}[ht]
\small\sf\centering
\caption{Comparison of single-shooting algorithms initialized with hovering controls at ten different initial states for quadpend.}
\label{tab:single_DDP}
\begin{center}
\begin{tabular}[h]{r||c c c c c c} \hline
        & \makecell{Log\\ Barrier}
        & \makecell{AL\\ single}
        & \makecell{IP\\} 
        & \makecell{PDAL\\ single}
        & \makecell{ADMM\\ single}
        & \makecell{SQP \\single}
        \\
 \hline
Success rate [\%] $\uparrow$
 & \bf{100} 
 & \bf{100} 

 & \bf{100} 
 & \bf{100} 
 & \bf{100} 
 & 60 
 \\
  \hline
 
 $m_J$ $\downarrow$   & 29.4 
 & 29.2 
 & 31.0 
 & \bf{28.5}
 & 31.5
 & 30.6
 \\
 $\sigma_{J}$ $\downarrow$  & 3.39 
& 3.16 
& 4.70 
& \bf{3.07} 
& 5.38
& 4.10
 \\
 \hline 
 $m_{\rm{I}}$ $\downarrow$
 & $\checkmark$ 
 & \num{1.42e-9} 
 & $\checkmark$ 
 & \num{1.82e-8}
 & \num{1.83e-2}
 & \num{8.26e-7}
 \\

 $\sigma_{\rm{I}}$ $\downarrow$ 
 & $\checkmark$  
 & \num{1.97e-9}
 & $\checkmark$ 
 & \num{1.37e-8}
 & \num{4.40e-2}
 & \num{1.85e-6}
 \\ \hline
 \end{tabular}
\end{center}
\end{table*}


\begin{table*}[ht]
\small\sf\centering
\caption{Comparison of multiple-shooting algorithms initialized with hovering controls at ten different initial states for quadpend.}
\label{tab:multi_DDP}
\begin{center}
\begin{tabular}[h]{r||c c c c} \hline
        & \makecell{AL \\multi exact}
        & \makecell{AL \\multi approx.}
        &  \makecell{PDAL\\ multi}
        & \makecell{SQP \\multi}\\
 \hline
 Success rate[\%] $\uparrow$

 & 80 
 & 80 
 & 90 
 & \bf{100} 
 \\
  \hline
 
 $m_J$ $\downarrow$   
 & 29.4
 & \textbf{28.3}
 & 29.7
 & 28.5
 \\
 $\sigma_{J}$ $\downarrow$  
& 4.37 
& 3.42
& 3.94
& 3.39 

 \\
 \hline 
 $m_{\rm{I}}$ $\downarrow$
 & \num{7.17e-6}
 & \num{1.61e-7}
 & \num{5.11e-7}
 & \textbf{\num{5.13e-11}}
 \\

 $\sigma_{\rm{I}}$ $\downarrow$ 
& \num{1.79e-5}
& \num{1.98e-7}
 & \num{8.98e-7}
 & \textbf{\num{6.16e-11}} \\
  \hline 
 $m_{\rm{E}}$ $\downarrow$
 & \num{1.80e-5}
 & \num{2.10e-7}
 & \num{7.35e-7}
 & \bf{\num{8.14e-7}}
 \\

 $\sigma_{\rm{E}}$ $\downarrow$ 
& \num{1.93e-6}
 & \num{2.61e-7}
 & \num{7.93e-7}
 & \bf{\num{2.30e-6}}
 \\ \hline
 \end{tabular}
\end{center}
\end{table*}

\begin{table*}[h]
\small\sf\centering
\caption{Comparison of single-shooting algorithms for ten different targets for Panda.}
\label{tab:single_DDP_panda}
\begin{center}
\begin{tabular}[h]{r||c c c c c c} \hline
        & \makecell{Log\\ Barrier}
        & \makecell{AL\\ single}
        & \makecell{IP\\} 
        & \makecell{PDAL\\ single}
        & \makecell{ADMM\\ single}
        & \makecell{SQP \\single}
        \\
 \hline
Success rate [\%] $\uparrow$
 & 60 & 80 & 50 & \textbf{90} & \textbf{90} & 60
 \\
  \hline
 
 $m_J$ $\downarrow$   &
\textbf{7.94} & 8.81 & 8.80& 8.82 & 9.41 & 10.7
 \\
 $\sigma_{J}$ $\downarrow$  &
1.50 & 2.54 & 1.99 & 2.48 & 2.67 & 3.28
 \\
 \hline 
 $m_{\rm{I}}$ $\downarrow$
 & $\checkmark$  &
 \num{9.59e-6} &
 \num{7.45e-9} &
 \num{3.76e-5} &
 \num{3.35e-2} &
 \num{3.73e-4}
 \\

 $\sigma_{\rm{I}}$ $\downarrow$ 
 & $\checkmark$  &
 \num{1.66e-5} & 
 \num{1.49e-8} & 
 \num{9.82e-5} & 
 \num{4.51e-2} & 
 \num{2.90e-4}
 \\
 \hline
 \end{tabular}
\end{center}
\end{table*}

\begin{table*}[h]
\small\sf\centering
\caption{Comparison of multiple-shooting algorithms for ten different targets for Panda.}
\label{tab:multi_DDP_panda}
\begin{center}
\begin{tabular}[h]{r||c c c} \hline
        & \makecell{AL multi}
        & \makecell{PDAL multi}
        & \makecell{SQP multi}
        \\
 \hline
Success rate [\%] $\uparrow$
 & 0 & 60 & \textbf{100}
 \\
  \hline
 
 $m_J$ $\downarrow$   &
 - & 10.3 & \textbf{9.09}
 \\
 $\sigma_{J}$ $\downarrow$  &
 - & 2.36 & 2.45
 \\
 \hline 
 $m_{\rm{I}}$ $\downarrow$
 & -  
 & \num{4.24e-4} 
 & \num{1.45e-4}
 \\

 $\sigma_{\rm{I}}$ $\downarrow$ 
 & -
 & \num{8.10e-4} 
 & \num{1.62e-4}\\
 \hline 
 $m_{\rm{E}}$ $\downarrow$
 & -  
 & \num{3.95e-4} 
 & \num{1.42e-3} 

 \\

 $\sigma_{\rm{E}}$ $\downarrow$ 
 & -
 & \num{7.25e-5}
 & \num{1.14e-3}
 \\
 
 \end{tabular}
\end{center}
\end{table*}

\begin{table*}
\small\sf\centering
\caption{Comparison of multiple shooting methods with good initial guess with ten initial points.
}
\label{tab:multiple_good_guess}
\begin{center}
\begin{tabular}[h]{r||c c c c} \hline
&\makecell{AL multi exact}
& \makecell{AL multi approx.}
& \makecell{PDAL multi}
& \makecell{SQP multi}\\
 \hline
 Success rate [\%] $\uparrow$
 & 80 
 & \bf{100} 
 & \bf{100} 
 & \bf{100} 
 \\
 \hline
 
 $m_J$ $\downarrow$
 & 48.8
 & 46.4
 & 46.6
 & \bf{44.0}
 \\
 
 $\sigma_{J}$ $\downarrow$
& 7.18 
& 5.03
& 5.12
& \bf{4.08} 
 \\
 \hline 
 $m_{\rm{I}}$ $\downarrow$
 & \num{7.44e-5} 
 & \num{5.24e-6} 
 & \num{3.10e-6} 
 & \num{3.31e-10} 
 \\

 $\sigma_{\rm{I}}$ $\downarrow$ 
 & \num{1.03e-4}
 & \num{1.04e-5} 
 & \num{8.82e-6} 
 & \num{5.19e-10} 
 \\
 \hline
 $m_{\rm{E}}$ $\downarrow$
 & \num{3.68e-5} 
 & \num{2.04e-6}  
 & \num{6.03e-7} 
 & \num{1.89e{-8}} 
 \\
 
 $\sigma_{\rm{E}}$ $\downarrow$
 & \num{8.47e-5} 
 & \num{4.30e-6}  
 & \num{7.63e-7} 
 & \num{4.44e-8} 
 \\
 \hline
 \end{tabular}
\end{center}
\end{table*}

\begin{figure}[htbp]
\centering
\includegraphics[trim={0cm 0cm 0cm 0cm},clip,width=\linewidth]{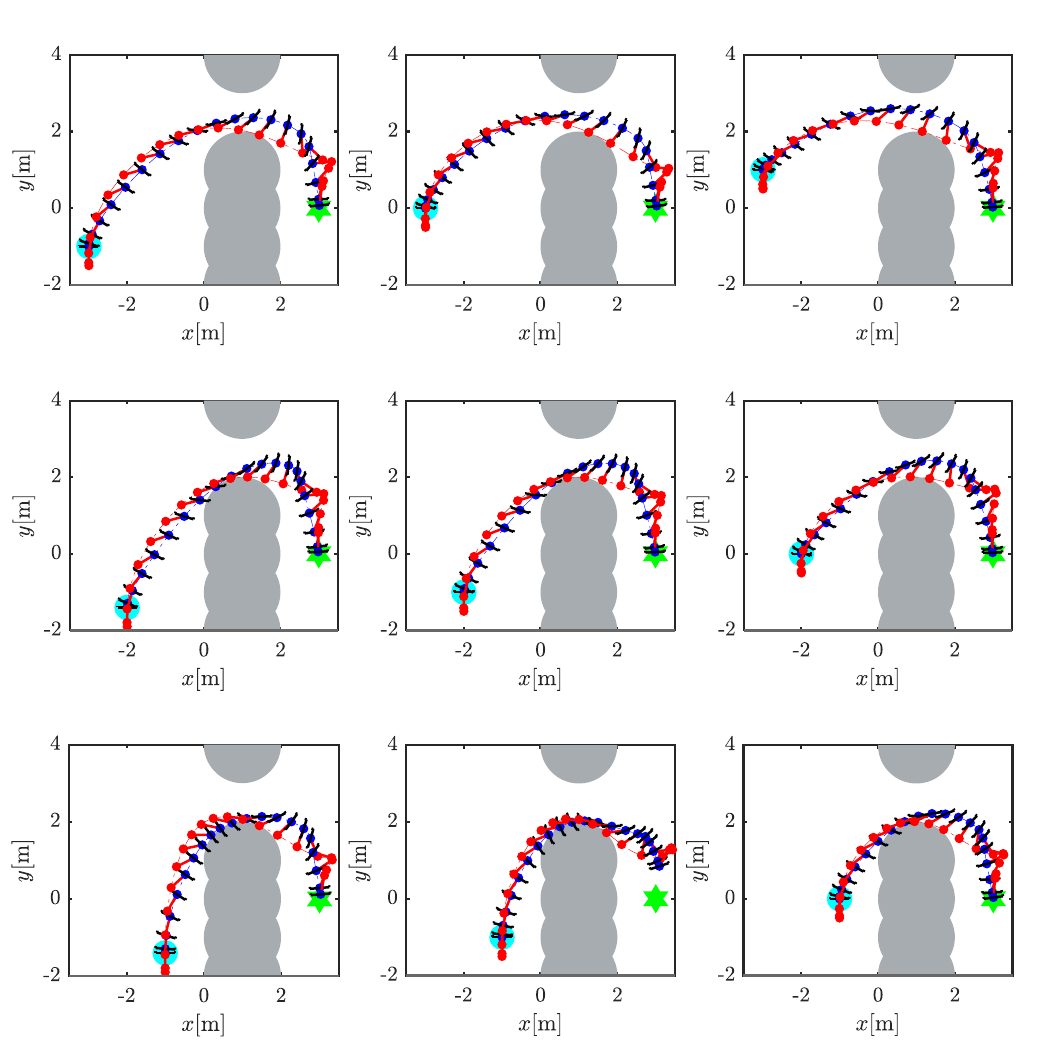}
\caption{Results of quadpend starting from multiple initial points with initial guesses by interpolation.}
\label{fig:sup_multi_warm}
\end{figure}

To circumvent this problem, we could solve a sequence of subproblems with intermediate targets and a short time horizon so that the quadpend can detour rather than get stuck at local minima. However, using the multiple shooting method with an informed initial guess, the problem should be solved without solving subproblems because the trajectory can be biased towards the right solution. To verify this idea, we initialize the multiple shooting algorithms with a state sequence, one of which is shown in the left figure of Fig. \ref{fig:quad_warm_init}. This trajectory is obtained by linear interpolation of four points such that the state trajectory does not hit obstacles. This initial guess leads to a reasonable solution presented in the right figure in
Fig. \ref{fig:quad_warm_init}. The same strategy as this interpolation is used for all ten different initial points to obtain informed initial trajectories.
The control sequence is initialized with the hover sequence presented in \ref{subsec:result_several_points}, which means that the equality constraints are violated on the initial trajectory. 

The results of this task are shown in Table \ref{tab:multiple_good_guess}.
The multiple-shooting method, except for the AL DDP with exact Hessian can complete the task. SQP is the most successful method among them, achieving the lowest cost and constraint violation.

\section{Conclusion}\label{sec:conclusions}
In this paper, we have reviewed two families of algorithms for constrained dynamic optimization: constrained DDP derived based on NLP techniques and SQP for dynamical systems. We have also discussed two distinct representations of these methods, namely, single- and multiple-shooting formulations. In addition, we derived a novel single-shooting PDAL DDP and added it to the comparison. Working towards our goal to systematize the research on second-order constrained dynamic optimization, we performed extensive benchmarking and analyzed algorithms based on criteria such as objective function minimization, task completion rate, and constraint satisfaction.
Among the different methods, the single-shooting PDAL DDP and multiple-shooting SQP algorithms stand out due to their consistent performance and robust numerical behavior across different systems and tasks. Both algorithms handle the high-dimensionality and non-convexity of trajectory optimization tasks in robotics well. When a small violation of dynamics is allowed, multiple-shooting SQP is the most stable method that can achieve the highest success rate of tasks. Another advantage of the method is that it can accept good initial guesses that help guide the optimization process to the desired trajectory.
On the other hand, when the infeasibility of dynamics is not allowed, the single-shooting PDAL DDP is suitable. 
This is because the method inherits the dynamical feasibility of the single-shooting DDP. 
 From a computational perspective, PDAL DDP has an advantage over multi-shooting SQP. This advantage originates from the fact that, in SQP, the inner constrained QP problem may typically require additional iterations to achieve convergence.  In PDAL DDP  and its backward pass, internal QPs  are unconstrained and therefore maintain their closed-form representation in which only the inversion of $ Q_{uu} $ is required. This is due to the way of how constraints are handled in PDAL DDP via the use of the Augmented Lagrangian.
The PDAL DDP has similar performance to the single-AL variant when the problem is simple, but starts to show its superiority as the problem becomes more complex. These two AL-based DDP methods can also be formulated in multiple-shooting formulations, which can take advantage of informed initial guesses. Finally, we would like to note that although  $\log$ barrier DDP struggles to complete tasks,  it can be the algorithm of choice when feasibility is prioritized.

\section*{Author Contributions}
\textbf{Yuichiro Aoyama}: Conceptualization, Formal analysis, Investigation, Software, Methodology, Writing-original draft. \textbf{Oswin So}: Conceptualization, Formal analysis, Investigation, Methodology, Software, Writing-review \& editing.
\textbf{Augustinos Saravanos}: Investigation, Methodology, Writing-review \& editing. 
\textbf{Evangelos Theodorou}: Conceptualization, Supervision, Project administration, Writing-review \& editing.
\section*{Statements and Declarations} 
\subsection*{Ethical considerations}
This article does not contain any studies with human or animal participants.
\subsection*{Consent to participate}
This article does not contain any studies with human or animal participants.
\subsection*{Consent for publication}
Not applicable.
\begin{dci}
The author(s) declared no potential conflicts of interest with respect to the research, authorship, and/or publication of this article.
\end{dci}
\begin{funding}
The author(s) disclosed receipt of the following financial support for the research, authorship, and/or publication of this article:  Augustinos D. Saravanos and Evangelos A. Theodorou were supported by ARO Award [grant number W911NF2010151]; 
Augustinos D. Saravanos was supported by the A. Onassis Foundation Scholarship. 
\end{funding}
\appendix
\section*{Appendix}
\section{Index to multimedia Extensions}
\begin{table}[h]
\centering
\resizebox{\columnwidth}{!}{
\begin{tabular}{ccc}
Extension & Media type & Description                                                 \\ \hline
1         & Video      & Movies of the trajectories of Panda.
\end{tabular}
}
\end{table}

\section{Detail of SQP for Dynamical Systems}\label{sup_sec:SQP}
This section is for detailed derivation of SQP and SQP for dynamical systems.
\subsection{General form of SQP}
The derivation of this section is based on those of the SQP solver SNOPT \citep{GillSNOPT2002} and NPSOL \citep{GILL1986NPSOL}. Nevertheless, we believe this section is important because some of the important techniques, e.g., solutions for problems to obtain parameters during optimization, were not explicitly explained.
\subsubsection{Problem formulation of SQP}
Let us revisit a SQP subproblem in \eqref{eq:SQP_QP}. Notice that the objective in has a gradient of $f_{0}$ and Hessian of $\mathcal{L}$, not Hessian of $f_{0}$. This discrepancy can be explained using a modified Lagrangian $\mathcal{L}_{\rm{m}}$ \citep{VanDerHoek1982, Robinson1972AQA}. Here we write
$x = x_{k} + \delta x_{k}$ with a subscript for iteration. 
\begin{align}\label{eq:SQP_modified_lagrangean}
    \mathcal{L}_{\rm{m}}(x, x_{k},\lambda_{k}) = f_{0}(x) + \lambda_{k}^{\tr}d_{\rm{L}}(x, x_{k}),
\end{align}
where $d_{\rm{L}}$ is the difference between the original and linearized constraints, which are written as
\begin{align*}
    d_{\rm{L}}(x,x_{k}) &= g(x) - g_{\rm{L}}(x, x_{k}),\\
    g_{\rm{L}}(x,x_{k}) &= g(x_{k}) + \nabla g(x_{k})(x-x_{k}).
\end{align*}
Since the gradient of $d_{\rm{L}}$ is given by
\begin{align*}
    \nabla d_{\rm{L}}(x, x_{k}) = \nabla g(x) - \nabla[g_{\rm{L}}(x,x_{k})] = \nabla g(x) - \nabla g(x_{k}),
\end{align*}
the gradient and Hessian of the modified Lagrangian are 
\begin{align*}
    \nabla \mathcal{L}_{\rm{m}}(x,x
    _{k},\lambda_{k}) &= \nabla f_{0}(x) + [\nabla g(x) - \nabla g(x_{k})]^{\tr}\lambda_{k},\\
    \nabla_{xx} \mathcal{L}_{\rm{m}}(x,x
    _{k},\lambda_{k}) &= \nabla_{xx} f_{0}(x) + \sum_{i=1}^{w}\big[[\lambda_{k}]_{i}\nabla_{xx} g_{i}(x)\big].
\end{align*}
Evaluated at $x=x_{k}$, $\mathcal{L}_{\rm{m}}$ and its gradient are equal to those of thr objective as 
\begin{align*}
    \mathcal{L}_{\rm{m}}(x_{k},x_{k},\lambda_{k}) = f_{0}(x_{k}), \
\nabla \mathcal{L}_{\rm{m}}(x_{k},x_{k},\lambda_{k}) = \nabla f_{0}(x_{k}).
\end{align*}
The Hessian of $\mathcal{L}_{\rm{m}}$, however, is equal to that of ${\mathcal{L}}$ at $x=x_{k}$,
\begin{align*} 
\nabla_{xx} \mathcal{L}_{\rm{m}}(x_{k},x_{k},\lambda_{k}) = \nabla_{xx} f_{0}(x_{k}) + \sum_{i=1}^{w}\big[[\lambda_{k}]_{i}\nabla_{xx} g_{i}(x_{k})\big].
\end{align*}
Therefore, quadratic approximation of $\mathcal{L}_{\rm{m}}$,  at $x = x_{k}$, denoted by $\mathcal{L}_{\rm{mq}}$, is obtained as
\begin{align*}
    &\mathcal{L}_{\rm{mq}}(x,x_{k},\lambda_{k})\\
    =& f_{0}(x_{k}) + [\nabla f_{0}(x_{k})]^{\tr}\delta x_{k} +\frac{1}{2}\delta x_{k}^{\tr}\nabla_{xx}\mathcal{L}(x_{k})\delta x_{k}.
\end{align*}
Given $x_{k}$, minimization of $\mathcal{L}_{\rm{mq}}$ is equivalent to that of 
\begin{align}\label{sup_eq:SQP_QP_subprob_simple}
     [\nabla f_{0}(x_{k})]^{\tr}\delta x_{k} +\frac{1}{2}\delta x_{k}^{\tr}\nabla_{xx}\mathcal{L}(x_{k})\delta x_{k},
\end{align}
which recovers \eqref{eq:SQP_QP}. More explanation is found in \cite{boggs_tolle_1995_SQP}. After solving the QP subproblem, $x_{k}$ is updated by the solution of \eqref{eq:SQP_QP} denoted by $\delta x^{\ast}_{k}$, giving $x_{k+1}$ for a new subproblem.

\subsubsection{Merit function and line search in SQP}\label{sec:sup_SQP_merit}
To determine an appropriate step size $\alpha$, we use Augmented Lagrangian (AL) merit function that achieves constraint satisfaction and cost reduction \citep{Gill1986SomeTP}. With the AL merit function, the appropriate step size is found by line search after setting the nonnegative penalty parameter $\rho \in \mathbb{R}^{w}$ so that it ensures the existence of good $\alpha$. Let $\phi$ be the AL merit function. $\phi$ is defined as a function of the penalty parameter parameterized by $\alpha$ as
\begin{align}{\label{eq:SQP_AL_merit}}
    \phi(\rho;\alpha) &= f(x_{k} + \alpha \delta x_{k})\\\notag 
    &\quad  +(\lambda_{k} + \alpha \delta \lambda_{k})^{\tr}[g(x_{k}+\alpha \delta x_{k})
     +(s+\alpha \delta s)] \\\notag &\quad\quad  + \sum_{i=1}^{w}\frac{1}{2}\rho_{i}\big[g_{i}(x_{k}+\alpha \delta x_{k})+[s+\alpha \delta s]_{i}\big]^2,
\end{align}
where $s\in \mathbb{R}^{w}$ is a slack variable that is introduced only for AL and line search. $\delta \lambda_{k}$ here is defined as a difference between optimal $\lambda_{k}$ denoted by $\lambda^{\ast}_{k}$ and current $\lambda_{k}$, that is
\begin{align*}
    \delta \lambda_{k} = \lambda^{\ast}_{k} - \lambda_{k}.
\end{align*}
$\lambda^{\ast}_{k}$ is obtained as byproducts of the solution of \eqref{eq:SQP_QP}. Here, we respect the original work and have negative sign for $\lambda$. $\delta x_{k}$ is the solution for \eqref{eq:SQP_QP}, but we drop $\ast$ for readability. The slack variable $s$ is initialized by
\begin{align*}
    s_{i} = \begin{cases}
    \max \{ 0, -g_{i}(x_{k}) \}, & \rho_{i} = 0, \\
    \max \{ 0, -g_{i}(x_{k})-\lambda_{i}/\rho_{i} \}, & \text{otherwise}. \\
    \end{cases}
\end{align*}
Its search direction is defined to satisfy
\begin{align}\label{eq:SQP_slack}
    g(x_{k}) + \nabla g(x_{k})\delta x_{k} = -(s + \delta s).
\end{align}
Since only the current iteration's variables matter in the merit function, we drop subscript $k$ for iteration hereafter. 
To ensure the existence of $\alpha$, $\rho$ needs to make $\phi'(\rho_{k};\alpha = 0)$ sufficiently negative, where $\phi' = \frac{\partial\phi}{\partial \alpha}$. This condition is typically given by
\begin{align}\label{eq:SQP_requirement_phi'(0)}
    \phi'(\rho;0) \leq -\frac{1}{2}\delta x^{\tr}H\delta x.
\end{align}
In SQP solver SNOPT \citep{GillSNOPT2002}, a minimum norm $\rho$, that achieves \eqref{eq:SQP_requirement_phi'(0)}, is used in the merit function. We follow the same strategy in our implementation. This $\rho$ has the following form.
\begin{align}\label{eq:SQP_min_norm_rho_sol}
    \rho^{\ast}_{i} &=
    \begin{cases}
     0, \quad   b \leq 0,\\
     \frac{b}{||a||^2}a_{i}, \ b > 0,
    \end{cases}\\
    \text{with} \ \quad 
    \notag
    a_{i} &= [g_{i}(x) + s_{i}]^{2}, \\ 
    \notag
    b = [\nabla f_{0}(x)]^{\tr}\delta x &+ \frac{1}{2}\delta x^{\tr}H\delta x + [\delta \lambda -\lambda]^{\tr}
[g(x)+s].
\end{align}
\proof{
The problem solved to find $\rho$ is given by
\begin{align}{\label{eq:SQP_min_norm_rho_original}}
    \min_{\rho} \frac{1}{2}||\rho||^2, \quad \text {s.t.} \quad \phi'(\rho;0) \leq -\frac{1}{2}\delta x^{\tr}H\delta x.
\end{align}
The LHS of the constraint is computed as
\begin{align}\label{eq:SQP_phi_0_original}
    \phi'(\rho;0) &= [\nabla f_{0}]^{\tr}\delta x \\ \notag 
    &  + \delta \lambda^{\tr}\big[g(x) + s\big] +\lambda^{\tr}[\nabla g(x)\delta x +s] \\\notag
    & + \sum_{i}^{w}\rho_{i}[g_{i}(x)+ s_{i}][\nabla g_{i}(x)\delta x +\delta s_{i}].
\end{align}
Plugging \eqref{eq:SQP_slack} and \eqref{eq:SQP_phi_0_original} back into \eqref{eq:SQP_requirement_phi'(0)}, we obtain
\begin{align}\label{eq:SQP_rho_cnst}
[\nabla f_{0}(x)]^{\tr}\delta x + \frac{1}{2}\delta x^{\tr}H\delta x
+ [\delta \lambda -\lambda]^{\tr}
[g(x)+s] \\ \leq \notag \sum_{i=1}^{w}\rho_{i}[g_{i}(x)+s_{i}]^2.
\end{align}
Notice that the first two terms are the optimal value of QP subproblem \eqref{sup_eq:SQP_QP_subprob_simple}.
Now, the constraint of the problem \eqref{eq:SQP_min_norm_rho_original} is equivalent to \eqref{eq:SQP_rho_cnst}.
When the LHS of \eqref{eq:SQP_rho_cnst} is not positive, the solution for \eqref{eq:SQP_min_norm_rho_original} is $\rho_{i} = 0,   \forall i = 1, \cdots w$. 
We consider the case where the LHS of the constraint is positive. To simplify the expression, we use $a$, $b$ in \eqref{eq:SQP_min_norm_rho_sol} and write \eqref{eq:SQP_min_norm_rho_original} as 
\begin{align}\label{eq:SQP_min_norm_rho}
\min_{\rho}\frac{1}{2}||\rho||^{2}, \quad \text{s.t.} \quad a^{\tr}\rho \geq b>0, \quad \rho \geq 0,
\end{align}
where $b>0$ is because we are considering the case where the LHS of the constraints in \eqref{eq:SQP_rho_cnst} is positive.
Lagrangian of this problem is given by
\begin{align*}
    \mathcal{L} = \frac{1}{2}\rho^{\tr}\rho- {\lambda_{1}}(a^{\tr}\rho-b) - {\lambda_{2}}^{\tr}\rho, 
\end{align*}
with corresponding multipliers $\lambda_{1}\in\mathbb{R}$ and $\lambda_{2}\in\mathbb{R}^{w}$. KKT condition yields
\begin{subequations}
\begin{align}
\rho - {\lambda_{1}} a - \lambda_{2} &= 0,\label{eq:SQP_minnorm_1}\\ 
{\lambda_1} \geq 0, \quad 
\lambda_{1}(a^{\tr}\rho-b) &= 0,\label{eq:SQP_minnorm_2}\\
\lambda_{2,i}\geq0, \quad 
\lambda_{2,i} \rho_{i} &= 0,  \quad i =1\cdots w, \label{eq:SQP_minnorm_3}
\end{align}
\end{subequations}
From \eqref{eq:SQP_minnorm_1},
\begin{align*}
    \lambda_2 = \rho-\lambda_{1}a. 
\end{align*}
Plugging this into the second equation of \eqref{eq:SQP_minnorm_3} yields
\begin{align}\label{eq:SQP_minnorm_4}
    [\rho - {\lambda_{1}}a]_{i}\rho_{i} = 0.
\end{align}
If $\rho_{i} = 0, \forall i$, then $a^{\tr}\rho=0$, which contradicts with $a^{\tr}\rho = b >0$ in the constraints in \eqref{eq:SQP_min_norm_rho}. Therefore, we need $\rho_{i} = \lambda_{1} a_{i} \neq0$ for some $i$. Plugging this in the second equation of \eqref{eq:SQP_minnorm_2} leads to 
\begin{align}
    {\lambda_{1}}(\lambda_{1}||a||^{2} - b) = 0.
\end{align}
If $\lambda_{1} = 0$, $\rho = 0$, from \eqref{eq:SQP_minnorm_4}, we have $a^{\tr}\rho=0$, again. Therefore, $\lambda_{1} \neq 0$, which gives 
\begin{align*}
    \lambda_{1} = b /{||a||^{2}}.
\end{align*}
Thus, the solution $\rho^{\ast}$ is obtained as
\begin{align*}
    \rho^{\ast}_{i} = \lambda_{1}a_{i} = \frac{b}{||a||^2}a_{i}.
\end{align*}
}
In addition to the minimum norm problem, a damping mechanism is implemented to allow $\rho$ to oscillate only finite times over iterations, which gives $\bar{\rho}$    
\begin{align}\label{eq:SQP_dumping_rho}
    &\bar{\rho_{i}} = \max\{\rho_{i}^{\ast}, \hat{\rho_{i}}\},\\\notag 
    &{\text{where}}\
    \hat{\rho}_{i} = \begin{cases}
    \rho_{i}, & \text{if} \quad \rho_{i} < 4(\rho_{i}^{\ast} + \Delta_\rho),\\
    [\rho_{i}(\rho_{i}^{\ast}+\Delta_\rho)]^{1/2}, & \text{otherwise}.
    \end{cases}
\end{align}
with initial damping parameter $\Delta_{\rho} = 1$. With this $\bar{\rho}$, line search is performed to find $\alpha$. $\alpha = 1$ if the following conditions are satisfied.
\begin{subequations}
\begin{align*}
&\phi(1)-\phi(0) \leq \sigma \phi'(0), \ 
\text{and} \\ \notag &\quad \phi'(1) \leq \eta \phi'(0) \quad \text{or} \quad |\phi'(1)| \leq - \eta \phi'(0).
\end{align*}
\end{subequations}
with $0<\sigma\leq\eta<\frac{1}{2}$. Otherwise, $\alpha$ that satisfies the following conditions is found by by backtracking $\alpha$
\begin{subequations}
\begin{align*}
\phi(\alpha)-\phi(0) \leq \sigma \alpha \phi'(0),\quad 
\text{and} \quad |\phi'(\alpha)| \leq - \eta \phi'(0).
\end{align*}
\end{subequations}

\subsubsection{Hessian Update}
After updating $x$ with appropriate $\alpha$, the new subproblem is obtained, whose gradient is computed with the new $x$. For the Hessian, since its exact computation is expensive and needs to be PD, BFGS quasi-Newton update \citep{Broyden1970BFGS, Fletcher1970BFGS, Goldfarb1970BFGS, Shanno19070BFGS} is used, which is give as follows.
\begin{align}\label{eq:SQP_BFGS_update}
    & \hspace{8mm} H_{k+1} = H_{k} + \theta_{k}y_{k}y_{k}^{\tr}-\psi_{k}q_{k}q_{k}^{\tr},\\\notag
\text{with} \\
\notag
y_{k} &= \nabla \mathcal{L}_{\rm{m}}(x_{k+1},x_{k},\lambda_{k+1})
    - \nabla  \mathcal{L}_{\rm{m}}(x_{k},x_{k},\lambda_{k+1})\\\notag
    &= \nabla f_{0}(x_{k+1}) - \nabla f_{0}(x_{k}) \\\notag 
    & \hspace{8mm} + [\nabla g(x_{k+1})-\nabla g(x_{k})]^{\tr}\lambda_{k+1},\\\notag
    \delta_{k} &= x_{k+1} - x_{k}, \
    q_{k} = H_{k}\delta_{k}, \
    \theta_{k} = \frac{1}{y_{k}^{\tr}\delta_{k}}, \ \psi_{k} = \frac{1}{q_{k}^{\tr}\delta x_{k}}.
\end{align}
This update law can keep $H_{k+1}$ PD if $H_{k}$ is PD and the approximate curvature $y_{k}^{\tr}\delta_{k}$ is positive. The term curvature comes from the curvature condition, which states that the step size $\alpha$ must satisfy. 
\begin{align*}
    \nabla f_{0}(x_{k}+\alpha p_{k})^{\tr} p_{k}\geq c \nabla f_{0}(x_{k})^{\tr}p_{k}
\end{align*}
for an objective $f_{0}(x)\in \mathbb{R}$, constant $c$ and a search direction $p_{k}$ to make a successful progress  
\citep{Nocedal2006numerical}. In our setting, the LH side is $\phi'(\alpha)$ and the RH side is a constant times $\phi'(0)$. This condition was originally used to decide when to terminate the line search. Setting $c=1$ gives the condition used here. When $y_{k}^{\tr}\delta x_{k}$ is not positive or has a very small positive value, the update law is modified to ensure that $H_{k+1}$ is PD. The modification is invoked when
\begin{align}\label{eq:SQP_curve_small}
    y_{k}^{\tr}\delta_{k} < \sigma_{k}, \quad \sigma_{k} = \alpha (1-\eta)\delta x_{k}H_{k}\delta x_{k},
\end{align}
 where $\eta$ is predefined constant $0<\eta<1$. 
Although SNOPT uses two modification techniques, we only tried the second one because the first one requires the intermediate results of the QP subproblems, which are unavailable to us because we use MATLAB {\texttt{quadprog}} function to solve the problem. The modification uses modified AL $\mathcal{L}_{\rm{mA}}$ rather than modified Lagrangian when computing $y_{k}$ as
\begin{align*}
    \mathcal{L}_{\rm{mA}}(x, x_{k},\lambda_{k};\omega) &= f_{0}(x) + \lambda_{k}^{\tr}d_{\rm{L}}(x,x_{k}) \\
    &\hspace{8mm}+ \frac{1}{2}d_{\rm{L}}(x,x_{k})^{\tr}\Omega d_{\rm{L}}(x,x_{k}), \\\notag
    {\text{with}} \quad \Omega &= {\rm{diag}}[\omega_{i}], \quad \omega_{i} \geq 0.
\end{align*}
The third penalty term with the parameter $\omega$ is added to the modified Lagrangian in \eqref{eq:SQP_modified_lagrangean}. 
Although the names are similar, this $\mathcal{L}_{\rm{mA}}$ has nothing to do with the AL merit function $\phi$ in \eqref{eq:SQP_AL_merit} used to find the step size $\alpha$. $\mathcal{L}_{\rm{mA}}$ is defined only to modify $y_{k}$ to update Hessian.
Here, new $y_{k}$ is obtained as a sum of current $ y_{k}$ and difference $\delta y_{k}$ as
\begin{align*}
    y_{k} + \delta y_{k} &= \nabla \mathcal{L}_{\rm{mA}}(x_{k+1},x_{k},\lambda_{k+1})-
    \nabla \mathcal{L}_{\rm{mA}}(x_{k},x_{k},\lambda_{k+1})\\
    &= \nabla \mathcal{L}_{\rm{m}}(x_{k+1},x_{k},\lambda_{k+1})-
    \nabla \mathcal{L}_{\rm{m}}(x_{k},x_{k},\lambda_{k+1})\\\notag
    &\hspace{12mm} +\underbrace{[\nabla g(x_{k+1}) - \nabla g(x_{k})]^{\tr}\Omega d_{\rm{L}}(x_{k+1},x_{k})}_{\delta y_{k}}.
\end{align*}
The modification uses minimum norm $\omega$ which satisfies 
$$[y_{k}+\delta y_{k}]^{\tr}\delta_{k} = \sigma_{k}.$$
The solution is given as 
\begin{align}\label{eq:SQP_omega_final}
\omega_{i} = 
    \begin{cases}
    0, & \bar{a}_{i}\leq 0\\
    \frac{\bar{b}}{||\bar{a}||^{2}}\bar{a}_{i}, &\bar{a}_{i}> 0.
    \end{cases}
\end{align}
\begin{proof}
The condition above is equivalent to
\begin{align}\label{eq:SQP_min_omega_condition}
    \delta y_{k}^{\tr} \delta_{k} &= \sigma_{k}-y_{k}^{\tr}\delta_{k}\\\notag
    & = \big[ [\nabla g(x_{k+1})-\nabla g(x_{k})]^{\tr}\Omega d_{\rm{L}}(x_{k+1},x_{k})\big]^{\tr} \delta_{k}\\\notag
    & = 
    \big[ \mathrm{diag}\{d_{\rm{L}}(x_{k+1},x_{k})\}[\nabla g(x_{k+1})-\nabla g(x_{k})]\delta_{k}\big]^{\tr}\omega.
\end{align}
To set up a problem to find $\omega$, define $\bar{a}\in\mathbb{R}^{w}$ and $\bar{b}\in\mathbb{R}$ as
\begin{align*}
    \bar{a} &= \mathrm{diag}\{d_{\rm{L}}(x_{k+1},x_{k})\}[\nabla g(x_{k+1})-\nabla g(x_{k})]\delta_{k},\\
    \bar{b} &= \sigma_{k}-y_{k}\delta_{k}>0.
\end{align*}
The inequality of $\bar{b}>0$ comes from \eqref{eq:SQP_curve_small}. Now, the problem of finding $\omega$ is
\begin{align}\label{eq:SQP_min_omega}
    \min_{\omega} \frac{1}{2}\omega^{\tr}\omega, \quad \text{subject to} \quad \bar{a}^{\tr}\omega = \bar{b} >0.
\end{align}
Lagrangian of this problem is
\begin{align*}
    \mathcal{L} = \frac{1}{2}\omega^{\tr}\omega - \lambda_{2}(\bar{a}^{\tr}\omega -\bar{b}),
\end{align*}
with a nonnegative Lagrangian multiplier $\lambda_{2}\in\mathbb{R}$. KKT condition yields
\begin{subequations}
\begin{align}
    \nabla \mathcal{L} = \omega - \lambda_{2}\bar{a} &= 0,\label{eq:SQP_min_omega_KKT}\\
    \lambda_{2}\geq0, \quad \bar{a}^{\tr}\omega &= b>0,\label{eq:SQP_min_omega_KKT2}\\
    \quad \lambda_{2}(\bar{a}^{\tr}\omega-\bar{b}) &=0. \label{eq:SQP_min_omega_KKT3}
\end{align}
\end{subequations}
If $\bar{a}_{i}\leq0, \forall i$, from \eqref{eq:SQP_min_omega_KKT}, $\omega=0$, which violates the constraints in \eqref{eq:SQP_min_omega_KKT2}. Therefore, we need $a_{i}>0$ for some $i$.
Then, from \eqref{eq:SQP_min_omega_KKT} we have
\begin{align}\label{eq:SQP_min_omega_2}
    \omega_{i} = \lambda_{2}\bar{a}_{i}.
\end{align}
Plugging this back into the constraints of \eqref{eq:SQP_min_omega} yields
\begin{align*}
\lambda_{2} ||\bar{a}||^{2} = \bar{b} \Leftrightarrow \lambda_{2} = \bar{b}/||\bar{a}||^{2},    
\end{align*}
which from \eqref{eq:SQP_min_omega_2} gives
\begin{align*}
\omega_{i}=\frac{\bar{b}}{||\bar{a}||^2}\bar{a}_{i}.
\end{align*}
For $\bar{a}_{i} \leq 0$, we put $\omega_{i} =0$, achieving minimum $\norm{\omega}$. 
\end{proof}

Using this $\omega$, Hessian for the next iteration $H_{k+1}$ is computed again by substituting $y_{k} + \delta y_{k}$ for $y_{k}$ in \eqref{eq:SQP_BFGS_update}. When $\omega^{\tr}\omega$ is too large, (say $10^5$), or $\omega$ does not exist, the modification is not performed, leaving the Hessian as it is for the next iteration.

\subsubsection{Merit function in SQP for dynamical systems}
This section show how the AL merit function can determine the step size $\alpha$ in SQP for dynamical systems. We take the multiple-shooting formulation as an example. 
Consider an optimization problem in \eqref{eq:SQP_multi}.

We define the penalty parameters for inequality and inequality constraints as $\rho_{\rm{I}}\in\mathbb{R}^{wN}$ and $\rho_{\rm{E}}\in\mathbb{R}^{nN}$. With these, define diagonal matrices
$
P_{\rm{E}} = {\rm{diag}}[\rho_{\rm{E}}], \quad P_{\rm{I}} = {\rm{diag}}[\rho_{\rm{I}}]$.
AL merit function $\phi$ for multiple shooting SQP and its derivative $\phi' = \frac{\partial{\phi}}{\partial \alpha}$ is obtained using a slack variable $s$ and Lagrangian multipliers $\lambda_{\rm{v}}$ for inequality and $\nu_{\rm{v}}$ for equality constraints as
\begin{align*}
&\phi(\rho_{\rm{I}},\rho_{\rm{E}};\alpha) \\
=& J(\bar{Y}+\alpha\delta Y)- (\lambda_{\rm{v}}+\alpha \delta \lambda_{\rm{v}})^{\tr}[G_{\alpha}(\delta Y)-(s+\alpha \delta s)]\\&+ \frac{1}{2}[G_{\alpha}(\delta Y)-(s+\alpha \delta s)]^{\tr}P_{\rm{I}}[G_{\alpha}(\delta Y)-(s+\alpha \delta s)]\\&
    -(\nu_{\rm{v}}+\alpha \delta \nu_{\rm{v}})^{\tr}\hat{F}_{\alpha}(\delta Y) +
    \frac{1}{2}\hat{F}_{\alpha}(\delta{Y})^{\tr}P_{\rm{E}}\hat{F}_{\alpha}(\delta Y),\\
&\text{with} \ 
G_{\alpha}(\delta Y) = G(\bar{Y} + \alpha \delta Y), \ \hat{F}_{\alpha}(\delta Y) = \hat{F}(\bar{Y} + \alpha \delta Y).
\end{align*}
Following the same procedure as in the derivation of \eqref{eq:SQP_phi_0_original}, we evaluate the merit function at $\alpha=0$
\begin{align*}
&\phi(\rho_{\rm{I}},\rho_{\rm{E}};0)'
%
= \hat{b} - \hat{a}_{2}^{\tr}\rho_{\rm{v}},\\\notag
&\text{where} \
     \rho_{\rm{v}} = [\rho_{\rm{I}}^{\tr}, \ \rho_{\rm{E}}^{\tr}]^{\tr}, \quad \hat{a}_{2} = \hat{a}_{1}\odot \hat{a}_{1}, \\
    &\hat{a}_{1} = [(G(\bar{Y})-s)^{\tr}, \hat{F}(\bar{Y})^{\tr}]^{\tr},\\\notag
    &\hat{b} = J_{Y}^{\tr}\delta  Y-[\delta  \lambda_{\rm{v}}-\lambda_{\rm{v}}]^{\tr}[G(\bar{Y})-s]
    -[\delta \nu_{\rm{v}}-\nu_{\rm{v}}]^{\tr}\hat{F}(\bar{Y}).
\end{align*}
 The minimum norm $\rho_{\rm{v}}$ that satisfies $\phi(\rho_{\rm{v}},0)'\leq -\frac{1}{2}\delta Y^{\tr} [\nabla_{YY} \mathcal{L}] \delta Y$ is obtained by substituting $\hat{a}_{2}$ and $\hat{b}$ for $a$ and $b$ in \eqref{eq:SQP_min_norm_rho_sol}, followed by the dumping mechanism in \eqref{eq:SQP_dumping_rho}.
With this penalty parameter, the line search is performed as in section \ref{sec:sup_SQP_merit}. When a reasonable step size is not found with a small $\alpha$, we discard the solution and regularize the Hessian as in \eqref{eq:DDP_regularization}, solving the QP subproblem again with the regularized Hessian.  

\subsubsection{Hessian Update}
In this section, we will discuss the Hessian update technique in multiple-shooting SQP for dynamical systems. As before, we first use \eqref{eq:SQP_BFGS_update} with modified Lagrangian 
\begin{align*}
    \mathcal{L}_{\rm{m}}(Y,Y_{k},\lambda_{\rm{v}}, \nu_{\rm{v}}) &= J(Y) + \lambda_{\rm{v}}^{\tr}d_{\rm{LI}}(Y, Y_{k})  +\nu_{\rm{v}}^{\tr}d_{\rm{LE}}(Y, Y_{k}),\\
    \text{with} \ d_{\rm{LI}}(Y,Y_{k}) &= G(Y)-G_{\rm{L}}(Y, Y_{k}), \\
    G_{\rm{L}}(Y,Y_{k}) &=  G(Y_{k})+G_{Y}(Y-Y_{k}),\\
    d_{\rm{LE}}(Y,Y_{k}) &= \hat{F}(Y)-\hat{F}_{\rm{L}}(Y, Y_{k}), \\
    \hat{F}_{\rm{L}}(Y,Y_{k}) &=  \hat{F}(Y_{k})+\hat{F}_{Y}(Y-Y_{k}).
\end{align*}
Here, the subscript $k$ in $Y_{k}$ stands for $Y$ in $k$ th iteration. Using this $\mathcal{L_{\rm{m}}}$, approximate curvature \eqref{eq:SQP_curve_small} is computed with 
\begin{align*}
    y_{k} &= \nabla \mathcal{L}_{\rm{m}}(Y_{k+1},Y_{k},\lambda_{{\rm{v}},k+1}, \nu_{{\rm{v}},k+1}) \\
    &\quad\quad  - \nabla \mathcal{L}_{\rm{m}}(Y_{k},Y_{k},\lambda_{{\rm{v}},k+1}, \nu_{{\rm{v}},k+1}).
\end{align*}
If the curvature is not large enough, we use modified AL
\begin{align*}
    &\mathcal{L}_{\rm{mA}}(Y,Y_{k},\lambda_{\rm{v}}, \nu_{\rm{v}}) \\
    =& J(Y) + \lambda_{\rm{v}}^{\tr}d_{\rm{LI}}(Y, Y_{k}) + \frac{1}{2}d_{\rm{LI}}(Y,Y_{k})^{\tr}\Omega_{\rm{I}}d_{\rm{LI}}(Y,Y_{k})\\
    &\hspace{8mm} + \nu_{\rm{v}}^{\tr}d_{\rm{LE}}(Y, Y_{k}) + \frac{1}{2}d_{\rm{LE}}(Y,Y_{k})^{\tr}\Omega_{\rm{E}}d_{\rm{LE}}(Y,Y_{k}),\\
    \text{with}
    &\quad \Omega_{I} = {\rm{diag}}[\omega_{I}], \quad \Omega_{E} = {\rm{diag}}[\omega_{E}], \quad \omega_{\rm{I}}, \omega_{\rm{E}} \geq0,
\end{align*}
which yields modification of $y_{k}$
\begin{align*}
    &\nabla \mathcal{L}_{\rm{mA}}(Y_{k+1},Y_{k},\lambda_{\rm{v}}, \nu_{\rm{v}}) - \nabla \mathcal{L}_{\rm{mA}}(Y_{k},Y_{k},\lambda_{\rm{v}}, \nu_{\rm{v}})\\
    =&y_{k}
    +[G_{Y}(Y_{k+1})-G_{Y}(Y_{k})]^{\tr}\Omega_{\rm{I}}d_{\rm{LI}}(Y_{k+1},Y_{k})\\\notag
    &\hspace{8mm} +[\hat{F}_{Y}(Y_{k+1})-\hat{F}_{Y}(Y_{k})]^{\tr}\Omega_{\rm{E}}d_{\rm{LE}}(Y_{k+1},Y_{k})\\
    =&y_{k} + \delta y_{k},\\
    \text{where} & \quad
    \delta y_{k} = G_{\rm{IE}}^{\tr}\Omega_{\rm{LE}}D_{\rm{L}},\\
    \text{with} & \hspace{2mm}
    G_{\rm{IE}} = \begin{bmatrix}
    G_{Y}(Y_{K+1}) - G_{Y}(Y_{k}) \\\notag
    \hat{F}_{Y}(Y_{k+1})-\hat{F}_{Y}(Y_{k})
    \end{bmatrix}, \\
    & \ \ \ D_{\rm{L}} = \mathrm{diag}\{
    d_{\rm{LI}}, d_{\rm{LE}} \},\\
    & \ \ \Omega_{\rm{LE}} = \begin{bmatrix}
    \Omega_{\rm{I}} & O_{nw, nN} \\ O_{nN, nw} & \Omega_{\rm{E}}
    \end{bmatrix}.
\end{align*}
Defining
\begin{align*}
    a^{\dagger} &= D_{\rm{L}} [\nabla G(Y_{k+1})-\nabla G(Y_{k})]^{\tr}(Y_{k+1}-Y_{k}), \\
    b^{\dagger} &= \sigma_{k} - y_{k}^{\tr}(Y_{k+1}-Y_{k}),
\end{align*}
and substituting $a^{\dagger}$ and $b^{\dagger}$ for $\bar{a}$ and $\bar{b}$ in \eqref{eq:SQP_omega_final} gives minimum norm penalty parameters $[\omega_{I}^T, \omega_{E}^{\tr}]^{\tr}$, which is used to compute $\delta y_{k}$.

\section{Appendix for IP DDP}\label{sec:app_IP}
In the main article, we have shown how $Q$ function and the propagation of the value function are modified from normal DDP to IP DDP. In this section, we complete the value function recursion by explaining its terminal condition. Then, we explain the difference between this recursion and that of the original work \cite{Pavlov2021IPDDP}. 

\subsection{Terminal condition with constraints.}
The terminal condition of the recursion of the value function is given by considering the cost, constraints, complimentary slackness at terminal time step as
\begin{align*}
    \hat{\Phi}(x_{N},\lambda_{N}) = \lambda^{\tr}_{N}g(x_{N}) + \Phi(x_{N}), \\
    g(x_{N}) + g_{x}\delta x_{N} + s_{N} + \delta s_{N} = 0,\\
    \bar{\Lambda}\delta s_{N} + S\delta\lambda_{N}  = - \bar{\Lambda}s_{N}+\mu e.
\end{align*}
From the second and the third equations, gains of $s$ and $\lambda$ at the terminal time step, which gives $\delta s_{N} = \kappa_{s} + K_{s} \delta x_{N}$ and $\delta \lambda_{N} = r + R \delta x_{N}$, are computed as
\begin{align*}
\kappa_{s} &= -(g(x_{N}) + s_{N}), \quad K_{s} = -g_{x}, \\
r &= S^{-1}[\bar{\Lambda} g(x_{N})+ \mu e], \quad R = S^{-1}\bar{\Lambda} g_{x}.
\end{align*}
Plugging $V_{N} = \hat{\Phi}$ and $\delta \lambda_{N} = r + R \delta x_{N}$ in quadratically expanded $V_{N}$ and mapping $\delta x$ terms with that of $V_{N}$ gives 
\begin{align*}
    V_{x,N} &= \hat{\Phi}_{x} + R^{\tr}\hat{\Phi}_{\lambda} + \hat{\Phi}_{x\lambda}r, \\
    V_{xx,N} &= \hat{\Phi}_{xx} + \hat{\Phi}_{x\lambda}R + R^{\tr}\hat{\Phi}_{\lambda x}.
\end{align*}
These two equations with $V_{N} = \hat{\Phi}$ are used as the terminal condition of the recursion in \eqref{eq:IPDDP_value_recursion_full}.
\subsection{Multipliers in the recursion}
In the main paper, we presented how $\delta \lambda$ affects the value function in \eqref{eq:IPDDP_value_recursion_full}. However, the authors of \cite{Pavlov2021IPDDP}, take a different approach. They claim that since the first and second equations of \eqref{eq:IP_DDP_KKT} lead to
\begin{align*} 
[\hat{Q}_{uu} + \hat{Q}_{u\lambda}S^{-1}\bar{\Lambda} g_{u}]\delta u_{k}  &= - \hat{Q}_{u} - Q_{u\lambda}S^{-1}[\bar{\Lambda} g - \mu] \\
&\hspace{5mm}
-[\hat{Q}_{u\bar{\lambda}}S^{-1}\bar{\Lambda} g_{x} + \hat{Q}_{u\lambda}]\delta x_{k},
\end{align*}
comparing this with \eqref{eq:delta-u-star} yields
\begin{align}\label{eq:IP_DDP_dQ_dagger}
    Q^{\dagger}_{u} &=  \hat{Q}_{u} + \hat{Q}_{u\lambda}S^{-1}[\bar{\Lambda} g(x_{k}, u_{k}) - \mu]\\\notag
    Q^{\dagger}_{uu} &= \hat{Q}_{uu} + \hat{Q}_{u\lambda}S^{-1}\bar{\Lambda} g_{u},\quad
    Q^{\dagger}_{ux} = \hat{Q}_{ux} \hat{Q}_{u\lambda}S^{-1}\bar{\Lambda} g_{x}.
\end{align}
Using these $\hat{Q}^{\dagger}$, they update $V$ with gains of $u_{k}$, i.e., $\kappa$ and $K$ in \eqref{eq:value_update_with_gain}. Observe that the gains for $\lambda_{k}$, that is, $r$ and $R$ are not used. This can be seen as taking $\delta \lambda_{k} \approx 0$ during the recursion of the value function, which may not be true. 
The authors also remove the constraints at the terminal time step. In most cases, the optimal solution is located within the feasible region. Hence, the modification would not significantly affect the problem. However, we have observed that the feedback gain of $\lambda_{N}$ increases when the corresponding $s_{N}$ approaches zero. As a result, the fraction to the boundary rule for DDP in \eqref{eq:DDP_fraction_to_boundary} is not satisfied even with a small step size $\alpha$ because it can only affect the feedforward part during the line search. The numerical instability of $s^{-1}$ when $s$ becomes close to zero is a typical problem with the IP method. One measure is the change of variables by multiplying $S$ to remove the term $s^{-1}$ \citep{Nocedal2006numerical}, but this is not implemented in IP DDP. In the experiment in section \ref{sec:experiments}, we used the same formulation as in the original work, i.e., $\delta \lambda \approx 0$ in the recursion of the value function and remove the constraints at the final time step.

\section{Matrix Transformation in PDAL DDP}\label{sec:app_pdal_transform}
This section provides matrix transformation for single- and multiple-shooting DDP introduced in section \ref{subsec:PDAL DDP}. Here, we first explain how we keep the Hessian matrices PD and then present the transformation.
\subsection{Single-shooting PDAL}
We first examine the Hessian matrix on the left-hand side of \eqref{eq:PDAL_DDP_single_opt_matrix} is PD. 
We introduce a nonsingular matrix $N_{\rm{s}}$ and denote the Hessian matrix, i.e., the matrix on the left-hand side of \eqref{eq:PDAL_DDP_single_opt_matrix}, as $H_{\hat{Q}}$
\begin{align*}
N_{\rm{s}} &= \begin{bmatrix}
I_{m} & [g_{u}]^{\tr}P_{\rm{I}}^{-1} \\ O_{w,m}  & I_{w}
\end{bmatrix}, \\
H_{\hat{Q}} &= \begin{bmatrix}
H_{\rm{s}} + 2[g_{u}]_{\mathcal{A}}^{\tr}P_{\rm{I}}[g_{u}]_{\mathcal{A}} &-[g_{u}]_{\mathcal{A}}^{\tr} \\
-[g_{u}]_{\mathcal{A}} & \frac{P^{-1}}{2} [[I_{w}]_{\mathcal{A}}+I_{w}]
\end{bmatrix},\\ \notag
&\text{with} \quad H_{\rm{s}} = Q_{uu} + [2\pi_{\rm{I}}-\lambda]_{+}[g_{uu}]_{\mathcal{A}}.
\end{align*}
$N_{\rm{s}}$ gives the following transformation.
\begin{align}\label{eq:sup_inertia}
N_{\rm{s}}^{\tr}H_{\hat{Q}}
N_{\rm{s}}
 = 
\begin{bmatrix}
H_{\rm{s}} + [g_{u}]_{\mathcal{A}}^{\tr}P_{\rm{I}}[g_{u}]_{\mathcal{A}} & O_{m,w} \\ O_{w,m} & \frac{P_{\rm{I}}^{-1}}{2}[[I_{w}]_{\mathcal{A}}+I_{w}]
\end{bmatrix},
\end{align}
the $(2,2)$ block is a diagonal matrix with positive elements. This block matrix has positive eigenvalues and is shared by \eqref{eq:sup_inertia}. From Sylvester's law of inertia, we know that the number of positive and negative eigenvalues is preserved under the transformation performed above \citep{Sylvester1852inertia,Ostrowski1959inertia}. Since \eqref{eq:sup_inertia} has $w$ nonnegative elements from its $(2, 2)$ block, the entire matrix is PD if the $(1, 1)$ block $H_{\rm{s}}+[g_{u}]^{\tr}_{\mathcal{A}}P[g_{u}]_{\mathcal{A}}$ is PD. If not, this block is regularized to ensure PSD by \eqref{eq:DDP_regularization} as other DDP methods.

Next, we explain the transformation applied to \eqref{eq:PDAL_DDP_single_opt_matrix}. The transformation is performed by a non-singular transformation matrix
\begin{align*}
    M_{\rm{s}} = \begin{bmatrix}
    I_{m} & {2}[P_{\rm{I}}^{-1}[g_{u}]_{\mathcal{A}}]^{\tr} \\
    0_{w, m} & I_{w}
    \end{bmatrix}.
\end{align*}
Multiplying $M_{\rm{s}}$ from the right side of \eqref{eq:PDAL_DDP_single_opt_matrix}, we have
\begin{align}\label{sup_eq:PDAL_DDP_lag_transformed}
    \begin{bmatrix}
    H_{\rm{s}} & [g_{u}]_{\mathcal{A}}^{\tr}\\ -[g_{u}]_{\mathcal{A}} &  \frac{{\rm{diag}}[\mu]}{2}\big[[I_{w}]_{\mathcal{A}} + I_{w}\big]
    \end{bmatrix}.
\end{align}
Note that $\rho$, a source of instability when it becomes large, disappears after transformation. $M_{\rm{s}}$ also transforms the right-hand side vectors as
\begin{align*}
    M_{\rm{s}}\begin{bmatrix}
    \hat{Q}_{u}\\ \hat{Q}_{\lambda}
    \end{bmatrix}
    &= 
    \begin{bmatrix}
    Q_{u} + [g_{u}]_{\mathcal{A}}^{\tr}\lambda \\ \hat{Q}_{\lambda}
    \end{bmatrix}\\
    M_{\rm{s}}\begin{bmatrix}
    \hat{Q}_{ux}\\ \hat{Q}_{\lambda}
    \end{bmatrix}
    &=
    \begin{bmatrix}
    Q_{ux} + [2\pi_{\rm{I}}-\lambda_{k}]_{+}[g_{ux}]_{\mathcal{A}} \\ \hat{Q}_{\lambda x}
    \end{bmatrix}.
\end{align*}
These transformations lead to the symmetric system presented in the main article.
\subsection{Multiple-shooting PDAL DDP}
Consider optimality condition of quadratically expanded $Q(x,u,\lambda, \nu)$ in multiple-shooting PDAL DDP, which gives 
\begin{align*}
    \begin{bmatrix}
    \hat{Q}_{\tilde{u}\tilde{u}} & \hat{Q}_{\tilde{u}\lambda} &\hat{Q}_{\tilde{u} \nu}\\
    \hat{Q}_{\lambda \tilde{u}} & \hat{Q}_{\lambda \lambda} & O_{w,n} \\
    \hat{Q}_{\nu \tilde{u}} & O_{n,w} & \hat{Q}_{\nu \nu}
    \end{bmatrix}
    \begin{bmatrix}
    \delta \tilde{u}_{k} \\ \delta \lambda_{k} \\ \delta \nu_{k}
    \end{bmatrix}
    = 
    -\begin{bmatrix}
    \hat{Q}_{\tilde{u}}\\ \hat{Q}_{\lambda} \\ \hat{Q}_{\nu}
    \end{bmatrix}
    - \begin{bmatrix}
    \hat{Q}_{\tilde{u}x} \\ \hat{Q}_{\lambda x} \\ \hat{Q}_{\nu x}
    \end{bmatrix} \delta x_{k}.
\end{align*}
Following the same procedure as in the single-shooting, the matrix on the left-hand side is PD if its $(1, 1)$ block is PD. Using a transformation matrix $M_{\rm{m}}$
\begin{align*}
    M_{\rm{m}} = \begin{bmatrix}
    I_{m} & 2[P_{\rm{I}}[g_{\tilde{u}}]_\mathcal{A}]^{\tr} & 2[P_{\rm{E}} [h_{\tilde{u}}]_\mathcal{A}]^{\tr} \\
    O_{w,m} & I_{w} & O_{w,n} \\
    O_{n,m} & O_{n,w} & I_{n}\\
    \end{bmatrix},
\end{align*}
we have a symmetric system
\begin{align*}
&
    \begin{bmatrix}
    H_{m} & -[g_{\tilde{u}}]_{\mathcal{A}}^{\tr} & -[g_{\tilde{u}}]_{\mathcal{A}}^{\tr} \\  -[g_{\tilde{u}}]_{\mathcal{A}} &  
    -M_{\mu}
    & O_{w,n} \\
    -h_{\tilde{u}} & O_{n,w} & -{\rm{diag}}[\mu_{\rm{E}}]
    \end{bmatrix}
    \begin{bmatrix}
    \delta u_{k} \\ -\delta \lambda_{k} \\ -\delta \nu_{k}
    \end{bmatrix} \\
    =&
    -\begin{bmatrix}
    Q_{\tilde{u}} + h_{\tilde{u}}^{\tr}\nu + [g_{\tilde{u}}]_{\mathcal{A}}^{\tr}\lambda \\ \hat{Q}_{\lambda}\\
    \hat{Q}_{\nu}
    \end{bmatrix}\\\notag
    &\hspace{6mm} - 
    \begin{bmatrix}
    Q_{\tilde{u}x} + [2\pi_{\rm{E}}-\nu][h_{\tilde{u}x}]+[2\pi_{\rm{I}}-\lambda]_{+}[g_{\tilde{u}x}]_{\mathcal{A}} \\ \hat{Q}_{\lambda x} \\\hat{Q}_{\nu x}
    \end{bmatrix}
    \delta x, \\
\text{with}& \
M_{\mu} = \frac{{\rm{diag}}[\mu_{{\rm{I}}}]}{2}\big[[I_{w}]_{\mathcal{A}} + I_{w}\big],\\
&H_{\rm{m}} 
= \underbrace{l_{\tilde{u}\tilde{u}} + \frac{\partial^2 V(\Pi(\tilde{u}))}{\partial \tilde{u}^2}}_{{Q}_{\tilde{u}\tilde{u}}}+ [2\pi_{\rm{E}}-\nu]h_{\tilde{u}\tilde{u}} \\ \notag
&\hspace{40mm} +[2\pi_{\rm{I}}-\lambda]_{+}[g_{\tilde{u}\tilde{u}}]_\mathcal{A},
\end{align*}
As in the previous case, $H_{\rm{m}}$ can be seen as a sum of Hessian of the problem's original objective and those of constraints multiplied by Lagrangian multipliers, and thus seen as Lagrangian of the constraint optimization problem.

\section{Details of Numerical Experiments}\label{sec:app_numerical}
In this section, we show details of experiments, including system dynamics, constraints, and cost structure. The dynamics include the inverted pendulum, quadpend, and swimmer. Finally, we present an additional experimental result of ADMM DDP of inverted pendulum with loosened control limits.  

\subsection{Cost Structure}
In the experiments, we have used the same cost structure as in \eqref{eq:unconstrained_optimalcontrol}. The running and terminal costs are given by
\begin{align}\label{sup_eq:numerical_example_quad_cost}
    l(y_{k}) &= 0.5 [u_{k}^{\tr}R_{1}u_{k} + (x_{k}-x_{\rm{g}})^{\tr}  R_{2}(x_{k}-x_{\rm{g}})], \\\notag \Phi(x_{N}) &= 0.5(x_{N}-x_{\rm{g}})^{\tr}Q(x_{N}-x_{\rm{g}}),
\end{align}
respectively. $R_{1} \in \mathbb{R}^{m}$, $R_{2}\in\mathbb{R}^{n}$ are the weight matrices for the running cost and $Q\in\mathbb{R}^{n}$ is that for the terminal cost. $x_{\rm{g}}$ is a desired state.

\subsection{Inverted Pendulum}
Consider the inverted pendulum. Let $l$, $m_{\rm{p}}$ be the length and mass of a pendulum, respectively. We assume that the mass is concentrated on the tip of the pendulum. The angle between a vertical line and the pendulum is $\theta$. The control of the system to be the torque $u$ applied to the joint. The dynamics of the system is given as follows.
\begin{align*}
    \begin{bmatrix}
    \theta_{k+1}\\
    \dot{\theta}_{k+1}
    \end{bmatrix}
     =
     \begin{bmatrix}
     \theta_{k}\\
     \dot{\theta}_{k}
     \end{bmatrix}
     + 
     \begin{bmatrix}
     \dot{\theta}_{k}\\
     \frac{u_{k}}{ml^2} - \frac{g_{0}}{l}\sin{\theta_{k}}
    \end{bmatrix}\Delta t,
\end{align*}
where $\Delta t$ is a discretization time interval.
We use $l=0.5$
, $m=0.2$
, $g = 9.81$
.  

For cost, we set the weight matrices as follows.
\begin{align*}
R_{1} =  0.001, \quad                             
R_{2} = {\rm{diag}}([100, 100]),\quad   
Q = 0.005 I_{2}.
\end{align*}
The time horizon is $N = 100$, the discretization time interval is $\Delta t = 0.02$. 



\begin{figure}[htbp]
  \centering
 \includegraphics[trim={0cm 20.5cm 10cm 3cm},clip,width=0.6\linewidth]{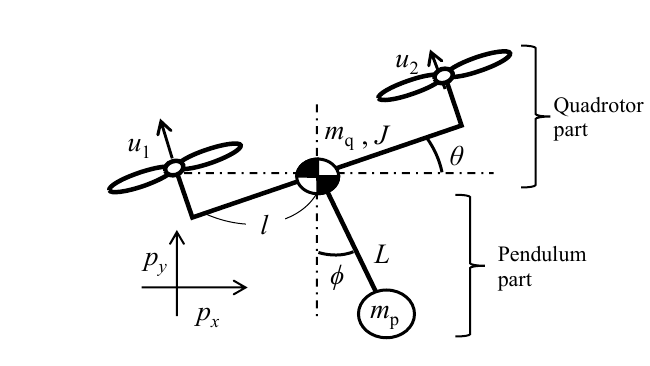}
  \caption{Schematic of a quadpend.}
  \label{fig:sup_quadpend_schematic}
\end{figure}

\subsection{Quadpend}\label{sup_sec:quad}
This section provides the dynamics and constraints of quadpend. 
\subsubsection{Dynamics and Parameters}
In this section, we derive the state-space representation of a quadpend based on the Lagrangian formulation given in \cite{Singh2022}. We define the position of the center of gravity of the quadrotor, its angle relative to the horizontal line, and the angle of the pendulum as $p_x$, $p_{y}$, $\theta$, and $\phi$, respectively, as shown in Fig.  \ref{fig:sup_quadpend_schematic}. The state of the system is given by $x = [p_{x}, p_{y}, \theta, \phi, \dot{p}_{x}, \dot{p}_{y}, \dot{\theta}, \dot{\phi}]^{\tr}$. The mass matrix of the system is 
\begin{align*}
    M = \begin{bmatrix}
    m_{\rm{q}} + m_{\rm{p}} & 0 & 0 & m_{\rm{p}}L\cos{\phi} \\
    0 & m_{\rm{q}} + m_{\rm{p}} & 0 &m_{\rm{p}}L\sin{\phi} \\
    0 & 0 &  J_{\rm{q}} & 0\\
    m_{\rm{p}}L\cos{\phi} &m_{\rm{p}}L\sin{\phi} & 0 & m_{\rm{p}}L^{2}
    \end{bmatrix},
\end{align*}
where $m_{\rm{q}}$ is mass of the quadrotor, $m_{p}$ is mass of the pendulum, $L$ is length of the pendulum, and $J_{\rm{q}}$ is inertia of quadrotor. We assume that the mass of the pendulum is concentrated on the tip. Lagrangian of the system $\mathcal{L}$ is given as a function of position and velocity in generalized coordinates as 
\begin{align*}
    \mathcal{L} = T(q,\dot{q})-V(q),
\end{align*}
with generalized position
\begin{align}\label{eq:sup_gen_pos}
q = [p_{x}, p_{y}, \theta, \phi]^{\tr},
\end{align}
kinetic energy $T(q, \dot{q})$ and potential energy $V(q)$. The energy terms are given by  
\begin{align*}
    T(q,\dot{q}) &= \frac{1}{2}\dot{q}^{\tr}M(q)\dot{q},\\
    V(q) &= m_{q}g_{0}p_{y} + m_{p}g_{0}(p_{y}-L\cos\phi),
\end{align*}
where we write the mass matrix $M$ as a function of $q$, explicitly.
Lagrangian formulation of the dynamics of the system is expressed as    
\begin{align}\label{eq:quadpend_lagrangian_dyn}
\odv{}{t}\frac{\partial{\mathcal{L}}}{\partial{\dot{q}}} - \frac{\partial \mathcal{L}}{\partial q} = F
\end{align}
with generalized force
\begin{align*}
F =
    [
    -(u_{1} + u_{2})\sin{\theta}, \ (u_{1}+u_{2})\cos\theta, \\ 
    \quad (u_{1}-u_{2})l-\tau_{f}, \
    \tau_{f}
    ]^{\tr},
\end{align*}
where control $u$ is thrust generated by two rotors and $\tau_{f} = -K_{\nu}(\dot{\phi}-\dot{\theta})$ is frictional torque with constant $K_{\nu}$. $l$ is distance from the center of the quadrotor to the rotors.
We plug the energy term into the Lagrangian to derive a state-space representation.
From the first term of \eqref{eq:quadpend_lagrangian_dyn} we get,
\begin{align}\label{sup_eq:quadpend_lagran_first_term}
    \odv{}{t}\frac{\partial \mathcal{L}}{\partial \dot{q}} &= \odv{}{t}M(q)\dot{q} = \dot{M}(q)\dot{q} + M(q)\ddot{q}.
\end{align}
$\dot{M}$ is computed as
\begin{align*}    
[\dot{M}(q)]_{i,j} &=     
\begin{cases}
-m_{\rm{p}}L\dot{\phi}\sin{\phi}, \ &\text{if} \ (i,j) = (1,4),(4,1),\\
m_{\rm{p}}L\dot{\phi}\cos{\phi}, \ &\text{if} \ (i,j) = (2,4),(4,2),\\
0, \ &\text{otherwise},
\end{cases}\\
\end{align*}
which gives
\begin{align}\label{eq:quadpend_lagrangian_dyn_1}
\dot{M}(q)\dot{q} &= [
-m_{p}L\dot{\phi}^{2}\sin{\phi}, \ m_{p}L\dot{\phi}^{2}\cos{\phi}, \\\notag
& \hspace{20mm}0, \
m_{p}L\dot{\phi}(-\dot{x}\sin{\phi} + \dot{y}\cos{\phi})
]^{\tr}.
\end{align}
For the second term of  \eqref{eq:quadpend_lagrangian_dyn}, we have
\begin{align}\label{sup_eq:quadpend_lagran_second_term}
\frac{\partial\mathcal{L}}{\partial q} &= \frac{\partial}{\partial q}\Big[\frac{1}{2}\dot{q}^{\tr}M(q)\dot{q}-V(q)\Big].
\end{align}
Since $M(q)$ depends only on $\phi$ in $q$, non-zero elements of $\frac{\partial M(q)}{\partial q}$ are given as follows.
\begin{align*}
\Big[\frac{\partial M(q)}{\partial q}\Big]_{1,4} &= \Big[\frac{\partial M(q)}{\partial q}\Big]_{4,1} = \frac{\partial m_{\rm{p}}L\cos{\phi}}{\partial \phi} = -m_{\rm{p}}L\sin\phi, \\
\Big[\frac{\partial M(q)}{\partial q}\Big]_{2,4} &= \Big[\frac{\partial M(q)}{\partial q}\Big]_{4,2} = \frac{\partial m_{\rm{p}}L\sin{\phi}}{\partial \phi} = m_{\rm{p}}L\cos\phi.
\end{align*}
Therefore,
\begin{align}\label{eq:quadpend_lagrangian_dyn_2}
    &\frac{\partial}{\partial q}\Big[\frac{1}{2}\dot{q}^{\tr}M(q)\dot{q}\Big]_{i} \\ =&
    \begin{cases}
    \frac{1}{2}\dot{q}^{\tr} \frac{\partial M(q)}{\partial \phi}\dot{q} = m_{p}L\dot{\phi}(-\dot{p}_{x}\sin{\phi} + \dot{p}_{y}\cos\phi), \quad i=4\\\notag
    0, \quad i = 1, 2, 3.
    \end{cases}
\end{align}
Differentiating $V(q)$ by $q$, we have
\begin{align}\label{eq:quadpend_lagrangian_dyn_2_V}
    \frac{\partial V(q)}{\partial q} = \begin{bmatrix}
    0, & (m_{q} + m_{p})g_{0}, & 0, & m_{p}Lg_{0}\sin\phi
    \end{bmatrix}^{\tr}.
\end{align}
By plugging \eqref{eq:quadpend_lagrangian_dyn_1} back into \eqref{sup_eq:quadpend_lagran_first_term}, and \eqref{eq:quadpend_lagrangian_dyn_2}, into \eqref{sup_eq:quadpend_lagran_second_term}, we have the two terms in \eqref{eq:quadpend_lagrangian_dyn}, and finally we have the state-space representation.
\begin{align*}
&M(q){\ddot{q}} = F + 
[
m_{p}L\dot{\phi}^2\sin{\phi}, \\& \quad  \ -(m_{p} + m_{q})g_{0}-m_{p}L\dot{\phi}^{2}\cos{\phi},  0,  -m_{p}L g_{0}\sin{\phi} 
]^{\tr},
\end{align*}
with generalized position in \eqref{eq:sup_gen_pos}. Using the Euler integration scheme, we get the dynamics of the quadpend.

The parameters are set to the following values, i.e., 
$m_p = 0.468$, $m_q = 0.2m_q$, $l = 0.25$, $L = 2l$, $g_{0}=9.81$, $J = 3.83\times10^{-3}$, $\nu = 0.01$. Weight matrices for cost function are
\begin{align*}
R_{1} &= 0.01I_{2}, \ R_{2} = 0.001{\rm{diag}}{[5, 5, 100, 10, 5, 5, 10, 10]} \\
Q &= 2{\rm{diag}}([100, 100, 10, 100, 50, 50, 10, 50]).
\end{align*}
The time horizon is $N = 100$, and the time discretization interval is $\Delta t = 0.02$. Control constraints are as follows
\begin{equation*}
    u_{\rm{l}} \leq u \leq u_{\rm{u}}, \ \text{with}\
    u_{\rm{l}} = 3 m_{\rm{q}}g_{0}[1, 1]^{\tr},\ u_{\rm{u}} = 0.1 m_{\rm{q}}g_{0}[1,1]^{\tr}. 
\end{equation*}

\subsection{Swimmer}
This section shows the dynamics of the swimmer, the details of which are given in the supplementary material of \cite{Tassa2007MPCDDP} in an unconstrained setting. In this work, the dynamics is provided by the center of mass $[x_{\rm{cm}}, y_{\rm{cm}}]$ and angles of the joints in a global reference frame. Here, we describe how to recover the positions of the links from them. We consider a swimmer with three links for simplicity, but the same derivation is used for an arbitrary number of links.
Let the position of the centers of the links $\bm{r}_{i}$, $i=1,\cdots,3$ and unit vectors along the $i$ th link to the nose link $\bm{t}_{i}$. We assume that the mass of the link ($m_{i}$) is evenly distributed. We use bold letters to emphasize that these are vectors. From the geometric relation and the definition of the center of mass, we have the following relation.
\begin{align}
\begin{cases}
&\bm{r}_{3} - \bm{r}_{2} = (l_{3}/2) \bm{t}_{3} + (l_{2}/2) \bm{t}_{2}, \\
&\bm{r}_{2} - \bm{r}_{1} = (l_{2}/2) \bm{t}_{2} + (l_{1}/2) \bm{t}_{1}, \\
&\sum_{i=1}^{3} m_{i} \bm{r}_{i} = 0, 
\end{cases}
\end{align}
with $\bm{t}_{i} = [\cos{\theta_{i}}, \sin{\theta_{i}}]^{\tr}$. These equations are written as systems of vector equations as follows.
\begin{align*}
\underbrace{
\begin{bmatrix}
1 & -1 & 0 \\
0 & 1 & -1 \\
m_{3} & m_{2} & m_{1}
\end{bmatrix}
}_{Q_{\rm{s}}}
\begin{bmatrix}
\bm{r}_{3} \\ \bm{r}_{2} \\ \bm{r}_{1} \\
\end{bmatrix}
=\frac{1}{2}
\underbrace{
\begin{bmatrix}
1 & 1 & 0 \\
0 & 1 & 1 \\
0 & 0 & 0
\end{bmatrix}
}_{A_{\rm{s}}}
\underbrace{
\begin{bmatrix}
l_{3} & 0 & 0 \\
0 & l_{2} & 0 \\
0 & 0 & l_{1}
\end{bmatrix}
}_{L_{\rm{s}}}
\begin{bmatrix}
\bm{t}_{3} \\ \bm{t}_{2} \\ \bm{t}_{1} 
\end{bmatrix}
\end{align*}
From the above equation, the $x$ coordinate of the tip of the links other than the nose is computed as 
\begin{align*}
 \bm{x}_{l} = x_{\rm{cm}}\bm{e} + \bm{r}_{x} - (1/2)Q_{\rm{s}}^{-1}A_{\rm{s}}L_{\rm{s}} \bm{t}_{x},
\end{align*}
with $\bm{e} = [1,1,1]^{\tr}$, $\bm{t}_{x} = [t_{3, x}, t_{2, x}, t_{1, x}]^{\tr}$ and $\bm{r}_{x} = [r_{3, x}, r_{2, x}, r_{1, x}]^{\tr}$ .
The position of the nose is 
\begin{align*}
x_{\rm{n}} = x_{\rm{cm}} + r_{3,x} + (1/2)l_{3}\cos{\theta_{3}}.
\end{align*}
For $y$ components, we have 
\begin{align*}
 \bm{y}_{l} = y_{\rm{cm}}\bm{e} + \bm{r}_{y} - (1/2)Q_{\rm{s}}^{-1}A_{\rm{s}}L_{\rm{s}} \bm{t}_{y},\\
 y_{\rm{n}} = y_{\rm{cm}} + r_{3,y} + (1/2)l_{3}\sin{\theta_{3}},
\end{align*}
where $\bm{t}_{y}$ and $\bm{r}_{y}$ are defined in the similar manner as those of $x$ elements. We set $l=1$ and $m=1$. For time, $N=80$, and $\Delta t$ = 0.02. Cost matrices are as follows. 
\begin{align*}
R_{1} &= \num{1e-4}I_{2}, \ R_{2} = 0.01Q \\
Q &= {\rm{diag}}([500, 500, 0.001, \cdots, 0.001]).
\end{align*}
Finally, the control constraint is $-100 \leq u_{i} \leq 100$.

\subsection{Panda}
This section provides the details of the experiment with Panda. We denote the position and velocity of the end effector as
\begin{equation*}
x_{\mathrm{e}} = [x_{\mathrm{ep}}^{\tr}, v_{\mathrm{e}}^{\tr}, \omega_{\mathrm{e}}^{\tr}]^{\tr} \in \Rb^{12},
\end{equation*}
where $x_{\mathrm{ep}}$,  $v_{\mathrm{e}}$, and $\omega_{\mathrm{e}}$ are position, velocity and angular velocity, respectively. We set the targets for the tasks by specifying the position $x_{\mathrm{ep}}$ and setting $v_{\mathrm{e}}=\omega_{\mathrm{e}} = 0$. The quantity $x_{\mathrm{e}}$ is computed from the state $x$. The limits of the joints are based on the suggested limits, and the command limits are from the \texttt{xml} file of the model \cite{menagerie2022github}. To encode obstacle constraints, we place spheres along the links of the arm so that the spheres cover the link. Then, we define constraints between the spheres and obstacles. The time horizon and discretization interval of the problem are $N=100$ and $\Delta t$ = 0.02. Cost matrices are as follows: 
\begin{align*}
R_{1} &= 0.01I_{7}, \ R_{2} = \mathrm{blkdiag}[0.8I_{3}, 0.2I_{3}, 0.02I_{3}] \\
Q &= \mathrm{blkdiag}[100I_{3}, I_{3}, 0.1I_{3}],
\end{align*}
where $\mathrm{blkdiag}$ stands for block diagonal matrix. 

\subsection{Additional Results for ADMM DDP}
\label{sec:sup_ADMM_loose_ctrl}
Here, we provide an additional example of ADMM DDP with the inverted pendulum with relaxed control bounds. In the results presented in section \ref{subsec:result_one_guess}, ADMM DDP shows slow progress compared to other methods. This is because the control limits are so tight that it cannot achieve the constraint satisfaction and completion of task simultaneously as mentioned in \ref{sec:AL_based_DDP_analysis}. Indeed, by loosening the control bounds, the performance of ADMM DDP is improved. 
\begin{figure}[!ht]
\centering
\includegraphics[trim={0.8cm 0.5cm 1.6cm 0.3cm},clip,width=\linewidth]{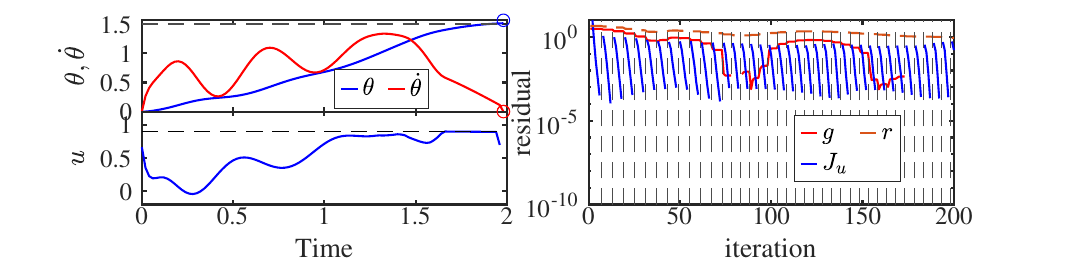}
\caption{ADMM DDP with loosened control bounds.}
\label{fig:pend/pend_admm_loose.pdf}
\end{figure}

\bibliographystyle{SageH}






\end{document}